\documentclass[invmat,final,envcountsect]{svjour}

\usepackage{amssymb}
\usepackage{amsmath}
\newcommand{\cst}{\mbox{\textnormal{Cst }}}
\newcommand{\eps}{\varepsilon}
\newcommand{\epsi}{\epsilon}
\newcommand{\jap}{\langle \xi\rangle}
\newcommand{\dsp}{\displaystyle}

\newcommand{\mafrS}{\mathfrak{S}}
\newcommand{\mfa}{\underline{\mathfrak{a}}}
\newcommand{\mfb}{\underline{\mathfrak{b}}}
\newcommand{\R}{{\mathbb R}}
\newcommand{\N}{{\mathbb N}}

\newcommand{\tX}{\widetilde{X}}
\newcommand{\tV}{\widetilde{V}}
\newcommand{\tz}{\widetilde{\zeta}}

\newcommand{\us}{\underline{s}}
\newcommand{\bv}{\underline{\bf v}}
\newcommand{\bbv}{{\bf v}}

\newcommand{\up}{\underline{\psi}}
\newcommand{\uz}{\underline{\zeta}}
\newcommand{\uZ}{\underline{Z}}
\newcommand{\uU}{\underline{U}}

\newcommand{\uC}{\underline{C}}
\newcommand{\uV}{\underline{V}}
\newcommand{\usi}{\underline{\sigma}}

\newcommand{\cL}{{\mathcal L}}

\newcommand{\cF}{{\mathcal F}}
\newcommand{\cA}{{\mathcal A}}
\newcommand{\cT}{{\mathcal T}}

\newcommand{\cQ}{{\mathcal Q}}
\newcommand{\cS}{{\mathcal S}}
\newcommand{\cE}{{\mathcal E}}
\newcommand{\cLu}{{\mathfrak L}_{(\uz,\up)}}
\newcommand{\cMu}{{\mathfrak M}_{(\uz,\up)}}
\newcommand{\cG}{{\mathcal G}}
\newcommand{\cZ}{{\mathcal Z}}
\newcommand{\Mu}{{M}_{(\uz,\up)}}

\newcommand{\dt}{\partial_t}
\newcommand{\dx}{\partial_x}
\newcommand{\dy}{\partial_y}
\newcommand{\dz}{\partial_z}

\newcommand{\Pig}{{\mathfrak P}}
\newcommand{\init}{{\vert_{t=0}}}
\newcommand{\surf}{{\vert_{z=0}}}
\newcommand{\fond}{{\vert_{z=-1}}}

\newcommand{\nag}{\nabla^{\gamma}}
\newcommand{\divg}{\mbox{\textnormal{div}}_\gamma}
\newcommand{\nagm}{\nabla^{\mu,\gamma}}

\newcommand{\Dg}{\vert D^{\gamma}\vert}
\newcommand{\xig}{\vert \xi^{\gamma}\vert}
\newcommand{\G}{\cG_{\mu,\gamma}[\eps\zeta,\beta b]}
\newcommand{\GG}{\cG[\eps\zeta]}

\newcommand{\GGu}{\cG[\eps\uz]}

\newcommand{\GGuo}{\mathcal{G}[0]}

\newcommand{\matr}{(1+Q[\sigma])}

\spnewtheorem{theo}{Theorem}[section]{\bf}{\rm}
\spnewtheorem{defi}{Definition}[section]{\bf}{\rm}
\spnewtheorem{assu}{Assumption}[section]{\bf}{\rm}

\spnewtheorem{nota}{Notation}[section]{\bf}{\rm}
\spnewtheorem{prop}{Proposition}[section]{\bf}{\rm}
\spnewtheorem{lemm}{Lemma}[section]{\bf}{\rm}
\spnewtheorem{coro}{Corollary}[section]{\bf}{\rm}
\spnewtheorem{rema}{Remark}[section]{\bf}{\rm}
\spnewtheorem{exam}{Example}[section]{\bf}{\rm}

\numberwithin{equation}{section}
\smartqed

\title{Large time existence for $3D$ water-waves and asymptotics}
\author{Borys Alvarez-Samaniego \and David Lannes}
\institute{Universit\'e Bordeaux I; IMB and CNRS UMR 5251, 
	  351 Cours de la Lib\'eration,
	  33405 Talence Cedex, France \\
          (\email{balvarez@math.uic.edu; David.Lannes@math.u-bordeaux1.fr})} 
\date{}
\begin{document}
\authorrunning{B. Alvarez-Samaniego \and D. Lannes}
\titlerunning{Large time existence for $3D$ water-waves and asymptotics}

\maketitle

\begin{abstract}
We rigorously justify in $3D$ the main asymptotic models 
used in coastal oceanography, including: shallow-water equations,
Boussinesq systems, Kadomtsev-Petviashvili (KP) approximation, \\
Green-Naghdi equations, Serre approximation and full-dispersion \\
model.  
We first introduce a ``variable'' nondimensionalized version of 
the water-waves equations which vary from shallow to deep water, 
and which involves four dimensionless parameters. Using
a nonlocal energy adapted to the equations, we can prove a
well-posedness theorem, uniformly with respect to all the parameters.
Its validity ranges therefore from shallow to deep-water, from small to large
surface and bottom variations, and from fully to weakly transverse waves.\\
The physical regimes corresponding to the aforementioned models can 
therefore be studied as particular cases; it turns out that the existence
time and the energy bounds given by the theorem are always those needed
to justify the asymptotic models. We can therefore derive and justify them
in a systematic way.
\end{abstract}
\section{Introduction}

\subsection{General setting}

The motion of a perfect, incompressible and irrotational fluid under the 
influence of gravity is described by the free surface Euler (or water-waves)
equations. Their complexity led physicists and mathematicians to
derive simpler sets of equations likely to describe the dynamics 
of the water-waves equations in some specific physical regimes. In
fact, many of the most famous equations of mathematical physics
were historically obtained as formal asymptotic limits of the
water-waves equations: the shallow-water equations, 
the Korteweg-de Vries (KdV) and 
Kadomtsev-Petviashvili (KP) equations, the Boussinesq
systems, etc. Each of these
asymptotic limits corresponds to a very specific physical regime whose
range of validity is determined in terms of the characteristics of the
flow (amplitude, wavelength, anisotropy, bottom topography, depth, ...).\\
The derivation of these models goes back to the XIXth century, but the
rigorous analysis of their relevance as approximate models for the
water-waves equations only began three decades ago with the works of
Ovsjannikov \cite{Ovsjannikov0,Ovsjannikov}, Craig \cite{Craig}, and 
Kano and Nishida \cite{KanoNishida1,KanoNishida2,Kano}
who first addressed the problem of justifying the formal asymptotics.
For all the different asymptotic models, the problem can be formulated
as follows: 1) do the 
water-waves equations have a solution on the time scale relevant for the 
asymptotic model and 2) does this model furnish a good approximation
of the solution? Answering the first question requires a large-time existence
theorem for the water-waves equations, while the second one requires
a rigorous derivation of the asymptotic models and a precise control
of the approximation error.

Following the pioneer works for one-dimensional surfaces ($1DH$) of 
Ovsjannikov \cite{Ovsjannikov} and Nalimov \cite{Nalimov} (see also
Yosihara \cite{Yosihara1,Yosihara2}), Craig \cite{Craig}, 
and Kano and Nishida \cite{KanoNishida1} provided the first
justification of the KdV and $1DH$ Boussinesq and shallow water 
approximations.  However, the comprehension
of the well-posedness theory for the water-waves equations hindered the
perspective of justifying the other asymptotic regimes until the
breakthroughs of S. Wu (\cite{Wu1} and \cite{Wu2} respectively for the 
$1DH$ and $2DH$ case, in infinite depth, and without restrictive assumptions).
Since then, the literature on free surface Euler equations has been
very active: the case of finite depth was proved in \cite{LannesJAMS}, and
in the related case of the study of the free surface
of a liquid in vacuum with zero gravity, Lindblad \cite{Lindblad1,Lindblad2}.
More recently Coutand and Shkoller \cite{CoutandShkoller} and
Shatah and Zeng \cite{ShatahZeng} managed to remove the irrotationality
condition and/or took into account surface tension effects (see also
\cite{AmbroseMasmoudi} for $1DH$ water-waves with surface tension).

In order to review the existing results of the rigorous justification
of asymptotic models for water-waves, it is suitable to classify
the different physical regimes using two dimensionless numbers: the amplitude 
parameter $\eps$
and the shallowness parameter $\mu$ (defined below in (\ref{KPreg})):
\begin{itemize}
	\item Shallow-water, large amplitude ($\mu\ll 1$, $\eps\sim 1$). 
Formally, this
regime leads at first order to the well-known ``shallow-water equations'' 
(or Saint-Venant) and at second order to the
so-called ``Green-Naghdi'' model, often used in coastal oceanography because it
takes into account the dispersive effects neglected by the shallow-water 
equations. The first rigorous justification of the shallow-water model
goes back to Ovsjannikov \cite{Ovsjannikov0,Ovsjannikov} and 
Kano and Nishida \cite{KanoNishida1} who proved the convergence
of the solutions of the shallow-water equations to solutions of the 
water-waves equations as $\mu\to 0$ in $1DH$, and under some restrictive 
assumptions
(small and analytic data). More recently, Y. A. Li \cite{LiCPAM}
removed these assumptions and rigorously justified the shallow-water and
Green-Naghdi equations, in $1DH$ for flat bottoms. Finally, the 
first and so far only
rigorous work on a $2DH$ asymptotic model is due
to a very recent work by T. Iguchi \cite{Iguchi2} in which he justified
the $2DH$ shallow-water equations, also allowing non-flat bottoms, but under
a restrictive zero mass assumption on the velocity.
	\item Shallow water, medium amplitude ($\mu\ll 1$, 
$\eps\sim\sqrt{\mu}$). This regime leads to the so-called Serre equations,
which are quite similar to the aforementioned Green-Naghdi equations and
are also often used in coastal oceanography. To our knowledge, no rigorous
result exists on that model.
	\item Shallow water, small amplitude ($\mu\ll 1$, 
$\eps\sim \mu$). This regime (also called long-waves regime) leads to
many mathematically interesting models due to the balance of nonlinear
and dispersive effects:
\begin{itemize}
	\item Boussinesq systems: since the first derivation by Boussinesq, 
many formally equivalent models (also named after Boussinesq) have 
been derived. W. Craig \cite{Craig} and Kano and Nishida \cite{KanoNishida2}
were the first to 
give a full justification of
these models, in $1DH$ (and for flat bottoms and small data). Note, however, 
that the convergence result given in \cite{KanoNishida2} is given on a time
scale too short to capture the nonlinear and dispersive 
effects specific to the Boussinesq systems; in \cite{Craig}, the correct
\emph{large time}
existence (and convergence) results for the water-waves equations are given.
The proof, of such a large time well-posedness result for the water-waves
equations, is the most delicate point in the justification process.  Furthermore, it
is the last step needed to fully justify the Boussinesq systems in $2DH$, 
owing to \cite{BCL} (flat bottoms) and \cite{Chazel} (general bottom topography),  where the convergence property is proved \emph{assuming}
that the large-time well-posedness theorem holds.
	\item Uncoupled models: at first order, the Boussinesq systems
reduce to a simple wave equation and,  in $1DH$, the motion of the free surface
can be described as the sum of two uncoupled counter-propagating waves, 
slightly modulated by a Korteweg-de Vries (KdV) equation. 
In $2DH$ and for weakly transverse
waves, a similar phenomenon occurs, but with the Kadomtsev-Petviashvili (KP) 
equation replacing the KdV equation. Many papers addressed the problem of
validating the KdV model 
(e.g. \cite{Craig,KanoNishida2,SchneiderWayne1,BCL,Wright,Iguchi1}) and
its justification is now complete. For the KP model, a first attempt
was done in \cite{Kano}, under restrictive assumption (small and analytic
data), but as in \cite{KanoNishida1}, the time scale considered is 
unfortunately too small for the relevant dynamics. A series
of works then proved the KP limit for simplified systems and toy models
\cite{GallaySchneider,BenyouLannes,Paumond}, while a different approach
was used in \cite{LannesSaut} where the KP limit is proved for the
full water-waves equations, \emph{assuming} a large-time well-posedness theorem and
a specific control of the solutions.
\end{itemize}
\item Deep-water, small steepness ($\mu\geq 1$, $\eps\sqrt{\mu}\ll1$). This
regime leads to the full-dispersion (or Matsuno) equations; 
to our knowledge, no rigorous result exists on this point. 
\end{itemize}

Instead of developing an existence/convergence theory for each 
physical scaling, 
we hereby propose a global method which allows one to
justify all the asymptotics mentioned above at once. In order to do
that, we nondimensionalize the water-waves equations, and keep track
of the five physical quantities which characterize the dynamics of the
water-waves: amplitude, depth, wavelength in the longitudinal direction,
wavelength in the transverse direction and amplitude of the
bottom variations.\\
Our main theorem gives an estimate of the existence time of the solution of the
water-waves equations which is \emph{uniform} with respect to \emph{all} these
parameters. In order to prove this theorem, 
we introduce an energy which involves the aforementioned parameters 
and
use it to construct our solution by an iterative scheme. Moreover, this energy 
provides some bounds on the solutions which appear to be exactly those needed
in the justification of the asymptotics regimes mentioned above. 

\subsection{Presentation of the results}

Parameterizing the free surface by $z=\zeta(t,X)$ (with $X=(x,y)\in\R^2$) 
and the bottom by
$z=-d+b(X)$ (with $d>0$ constant), one can use the incompressibility and 
irrotationality conditions to
write the water-waves equations under Bernouilli's formulation, in terms
of a velocity potential $\phi$ (i.e., the velocity field is given by 
${\mathbf v}=\nabla_{X,z}\phi$):
\begin{equation}
	\label{eqbern}
	\left\lbrace
	\begin{array}{lcl}
	\displaystyle \partial_{x }^2\phi 
	+\partial_{y }^2\phi 
	+\partial_{z }^2\phi =0,
	& & -d+b\leq z \leq \zeta ,\\
	\displaystyle \partial_n\phi=0, & &z =-d+b ,\\
	\displaystyle \partial_{t }\zeta +\nabla\zeta\cdot\nabla\phi
	=\partial_{z }\phi,
	& &  z = \zeta ,\\
	\displaystyle \partial_{t }\phi 
	+\frac{1}{2}\big(\vert\nabla\phi\vert^2+
	(\partial_{z }\phi )^2\big)+\zeta =0,
	& & z =\zeta,
	\end{array}
	\right.
\end{equation}
where $\nabla=(\dx,\dy)^T$ and $\partial_n\phi$ is the outward normal
derivative at the boundary of the fluid domain.\\
The qualitative study of the water-waves equations is made easier by
the introduction of dimensionless variables and unknowns. This requires
the introduction of various orders of magnitude linked to the physical
regime under consideration. More precisely, let us introduce the
following quantities:
$a$ is the order of amplitude of the waves;
$\lambda$ is the wave-length of the waves in the $x$ direction;
$\lambda/\gamma$ is the wave-length of the waves in the $y$ direction;
$B$ is the order of amplitude of the variations of the bottom topography.\\
We also introduce the following dimensionless parameters
\begin{equation}
	\label{KPreg}
	\frac{a}{d}=\eps,\quad
	\frac{d^2}{\lambda^2}=\mu,\quad
	\frac{B}{d}=\beta;
\end{equation}
the parameter $\eps$ is often called \emph{nonlinearity} parameter, while
$\mu$ is the \emph{shallowness} parameter. In total generality, one has
\begin{equation}\label{values}
	(\eps,\mu,\gamma,\beta)\in (0,1]\times (0,\infty)\times
	(0,1]\times [0,1]
\end{equation}
(the conditions $\eps\in (0,1]$ and $\beta \in [0,1]$ mean that the
the surface and bottom variations are at most of the order of 
depth ---$\beta=0$ corresponding to flat bottoms--- and the condition $\gamma\in (0,1]$ says
that the $x$ axis is chosen to be the longitudinal direction for weakly
transverse waves).\\
Zakharov \cite{Zakharov} remarked that the system (\ref{eqbern}) could
be written in Hamiltonian form in terms of the free surface elevation
$\zeta$ and of the trace of the velocity potential at the surface
$\psi=\phi_{\vert_{z=\zeta}}$ and Craig, Sulem and Sulem \cite{CSS1} 
and Craig, Schanz and Sulem \cite{CSS2}  used the
fact that (\ref{eqbern}) could be reduced to a system of two evolution 
equations on
$\zeta$ and $\psi$;  this formulation has commonly been used
since then. The dimensionless form
of this formulation involves the parameters introduced in (\ref{KPreg}), 
the transversity $\gamma$, and a parameter $\nu=(1+\sqrt{\mu})^{-1}$
whose presence is due to the fact that the nondimensionalization is not
the same in deep and shallow water. It is derived in Appendix \ref{appnd}:
\begin{equation}
	\label{nondimww}
	\left\lbrace
	\begin{array}{l}
	\dsp \dt \zeta-\frac{1}{\mu\nu}\cG_{\mu,\gamma}[\eps\zeta,\beta b]
	\psi=0,\\
	\dsp \dt \psi+\zeta+\frac{\eps}{2\nu}\vert\nag\psi\vert^2
	-\frac{\eps\mu}{\nu}\frac{(\frac{1}{\mu}\cG_{\mu,\gamma}
	[\eps\zeta,\beta b]\psi+\eps\nag\zeta\cdot\nag\psi)^2}
	{2(1+\eps^2\mu\vert\nag\zeta\vert^2)}=0,
	\end{array}\right.
\end{equation}
where $\nag=(\dx,\gamma\dy)^T$ and 
$\cG_{\mu,\gamma}[\eps\zeta,\beta b]$ is the Dirichlet-Neumann
operator defined by $\cG_{\mu,\gamma}[\eps\zeta,\beta b]\psi
=(1+\eps^2\vert\nabla\zeta\vert^2)^{1/2}\partial_n\Phi_{\vert_{z=\eps\zeta}}$,
with $\Phi$ solving
\begin{equation}\label{Laplintro}
	\left\lbrace
	\begin{array}{l}
	\dz^2\Phi+\mu \dx^2\Phi+\gamma^2\mu \dy^2\Phi=0,
	\qquad
	-1+\beta b <z <\eps \zeta,\\
	\Phi_{\vert_{z=\eps\zeta}}=\psi,\qquad
	\partial_n \Phi_{\vert_{z=-1+\beta b}}=0.
	\end{array}\right.
\end{equation}

In Section \ref{sectprel}, we give some preliminary results which will
be used throughout the paper: a few technical results (such as commutator
estimates) are given in \S \ref{sectcomm} and elliptic boundary value
problems directly linked to (\ref{Laplintro}) are studied in \S \ref{sectell}.\\
Section \ref{sectDN} is devoted to the study of various aspects of
the Dirichlet-Neumann operator $\cG_{\mu,\gamma}[\eps\zeta,\beta b]$. It is 
well-known that the Dirichlet-Neumann operator is a pseudo-differential
operator of order one; in particular, it acts continuously on Sobolev 
spaces and its operator norm, commutators with derivatives, etc., have
been extensively studied. The task here is more delicate because of the
presence of four parameters $(\eps,\mu,\gamma,\beta)$ in the operator
$\cG_{\mu,\gamma}[\eps\zeta,\beta b]$. Indeed, some of the classical 
estimates on the Dirichlet-Neumann are not uniform with respect to the
parameters and must be modified. But the main difficulty is due to the
fact that the energy introduced 
in this paper is not of Sobolev type; namely,
it is given by
\begin{equation}\label{introNRJ}
	\forall s\geq 0,\quad\forall U=(\zeta,\psi),\qquad
	\vert U\vert_{\tX^s}=
	\big\vert \zeta\vert_{H^s}
	+\vert \frac{\nu^{-1/2}\Dg}{(1+\sqrt{\mu}\Dg)^{1/2}}\psi\big\vert_{H^s},
\end{equation}
where $\Dg:=(-\partial_x^2-\gamma^2\partial_y^2)^{1/2}$. For high 
frequencies, this energy
is equivalent to the $H^s\times H^{s+1/2}$-norm specific to the
non-strictly hyperbolic nature of the water-waves equations 
(see \cite{Craig2} for a detailed comment on this point),
but the equivalence is not uniform with respect to the parameters,
and the $H^s\times H^{s+1/2}$ estimates of 
\cite{Wu1,Wu2,LannesJAMS,AmbroseMasmoudi} for 
instance, are useless for our purposes here. We thus have to work
with estimates in $\vert \cdot\vert_{\tX^s}$-type 
norms and the classical results
on Sobolev estimates of pseudodifferential operators cannot be used.
Consequently, we must rely on the structural properties of the water-waves
equations much more heavily than in the previous works quoted above.\\
Fundamental properties of the DN operators are given in \S \ref{sectfund},
while commutator estimates and further properties are investigated
in \S \ref{sectcomDN} and \S \ref{sectother}. 
We then give asymptotic expansions
of $\cG_{\mu,\gamma}[\eps\zeta,\beta b]\psi$ in terms of the parameters in 
\S \ref{sectas}.\\
Using the results of the previous sections, we study the 
Cauchy problem associated to the linearization
of (\ref{nondimww}) in Section \ref{sectLin}; the main energy
estimate is given in Proposition \ref{propmain}.\\
The full nonlinear equations are addressed in Section \ref{sectmain} and
our main result is stated in Theorem \ref{theomain}; it gives
a ``large-time'' (of order $O(\eps/\nu)$) 
existence result for the water-waves equations (\ref{nondimww}) and
a bound on its energy (defined in (\ref{introNRJ})). The most important
point is that this result is uniform with respect to \emph{all} the
parameters $(\eps,\mu,\gamma,\beta)$ satisfying (\ref{values}) and
such that the
steepness $\eps\sqrt{\mu}$ and the ratio $\beta/\eps$ remain bounded.
The theorem also requires a classical 
\emph{Taylor sign condition} on the initial data; we give in  
Proposition \ref{proptaylor} very simple sufficient conditions
(involving in particular the ``anisotropic Hessian'' of the bottom
parameterization $b$), which imply that the Taylor sign condition
is satisfied.\\
Both Theorem \ref{theomain} and Proposition \ref{proptaylor} can be
used for all the physical regimes given in the previous section, and the
solution they provide exists over a time scale relevant with respect to
the dynamics of the asymptotic models. We can
therefore study the asymptotic limits, which is done
in Section \ref{sectjustif}.
It is convenient to use
the classification introduced previously to present our results
(we also refer to \cite{LannesMoscou} for an overview of the
methods developed here):
\begin{itemize}
	\item Shallow-water, large amplitude ($\mu\ll 1$, $\eps\sim 1$). 
	We justify in \S \ref{sectjustifSW}
	the shallow-water equations without the restrictive
	assumptions of \cite{Iguchi2} and previous works. 
	For the Green-Naghdi model,
	we extend in \S \ref{sectjustifGNS} the
	result of \cite{LiCPAM} to non-flat bottoms, and to
	two dimensional surfaces. 
	\item Shallow water, medium amplitude ($\mu\ll 1$, 
$\eps\sim\sqrt{\mu}$). We rigorously justify the Serre approximation
	over the relevant $O(1/\sqrt{\mu})$ time scale 
	in \S \ref{sectjustifGNS}. 
	\item Shallow water, small amplitude ($\mu\ll 1$, 
$\eps\sim \mu$). 
\begin{itemize}
	\item Boussinesq systems: In 
	\S \ref{sectjustifLW}, we fully justify all the Boussinesq
	systems in the open case of two-dimensional surfaces 
	(flat or non-flat bottoms).
	\item Uncoupled models: We complete the full justification
	of the KP approximation in \S \ref{sectjustifKP}. 
\end{itemize}
\item Deep-water, small steepness ($\mu\geq 1$, $\eps\sqrt{\mu}\ll1$). 
We show in \S \ref{sectjustifFD} that the solutions of the full-dispersion model converge to
exact solutions of the water-waves equations as the steepness goes to zero
and give accurate error estimates.\\
We also give in \S \ref{sectjustifnum} 
an estimate on the precision of a model used for the
numerical computation of the water-waves equations (see \cite{CGHHS} for 
instance). 
\end{itemize}

\subsection{Notations}

\noindent
- We use the generic notation $C(\lambda_1,\lambda_2,\dots)$ to denote
a nondecreasing function of the parameters $\lambda_1,\lambda_2,\dots$.\\
- The notation $a\lesssim b$ means that $a \leq C b$, for some nonnegative
constant $C$ whose exact expression is of no importance ({\emph{in particular,
it is independent of the small parameters involved}).\\
- For all tempered distribution $u\in {\mathfrak S}'(\R^2)$, we denote
by $\widehat{u}$ its Fourier transform.\\
- Fourier multipliers: For all rapidly decaying $u\in {\mathfrak S}(\R^2)$
and all $f\in C(\R^2)$ with tempered growth, $f(D)$ 
is the distribution
defined by 
\begin{equation}
	\label{nota1}
	\forall \xi \in \R^2,\qquad
	\widehat{f(D)u}(\xi)=f(\xi)\widehat{u}(\xi);
\end{equation}
(this definition can be extended to wider spaces of functions).\\
- We write $\jap=(1+\vert \xi\vert^2)^{1/2}$ and $\Lambda=\langle D\rangle$.\\
- For all $1\leq p\leq \infty$, $\vert\cdot \vert_p$ denotes 
the classical norm of 
$L^p(\R^2)$ while $\Vert\cdot\Vert_p$ stands for the canonical norm of  
$L^p({\mathcal S})$, with ${\mathcal S}=\R^2\times (-1,0)$.\\
- For all $s\in\R$, $H^s(\R^2)$ is the classical Sobolev space defined as
$$
	H^s(\R^2)=\{u\in {\mathfrak S}'(\R^2), \vert u\vert_{H_s}:=\vert \Lambda^s u\vert_2<\infty\}.
$$
- For all $s\in\R$, $\Vert \cdot\Vert_{L^\infty H^s}$ denotes
the canonical norm of  \\
$L^\infty([-1,0];H^s(\R^2))$.\\
- If $B$ is a Banach space, then $\vert \cdot\vert_{B,T}$ stands for the
canonical norm of $L^{\infty}([0,T];B)$.\\
- For all $\gamma>0$, we write $\nag=(\partial_x,\gamma\partial_y)^T$,
so that $\nag$ coincides with the usual gradient when $\gamma=1$. We also
use the Fourier multiplier $\Dg$ defined as
$$
	\Dg=\sqrt{D_x^2+\gamma^2 D_y^2},
$$
as well as the anisotropic divergence operator
$$
	\divg=(\nag)^T.
$$
- We denote by $\Pig_{\mu,\gamma}$ (or simply $\Pig$ when no confusion is 
possible) the Fourier multiplier of order $1/2$
\begin{equation}\label{Pig}
	\Pig_{\mu,\gamma}(=\Pig):=\frac{\nu^{-1/2}\Dg}{(1+\sqrt{\mu}\Dg)^{1/2}}.
\end{equation}
- We write $X=(x,y)$ and $\nabla_{X,z}=(\dx,\dy,\dz)^T$; we also write
$$
	\nagm=(\sqrt{\mu}\dx,\gamma\sqrt{\mu}\dy,\dz)^T.
$$
- We use the condensed notation
\begin{equation}
	\label{nota3}
	A_s=B_s +\left\langle C_s\right\rangle_{s> \us}
\end{equation}
to say that $A_s=B_s$ if $s\leq \us$ and $A_s=B_s+C_s$ if $s> \us$.\\
- By convention, we take
\begin{equation}
	\label{conv1}
	\prod_{k=1}^0 p_k=1 \quad\mbox{ and }\quad
	\sum_{k=1}^0 p_k=0.
\end{equation}
- When the notation $\partial_n u_{\vert_{\partial\Omega}}$ is used for 
boundary conditions
of an elliptic equation of the form $\nabla_{X,z}\cdot P\nabla_{X,z}u=h$ in 
some open set $\Omega$,
it stands for the \emph{outward conormal derivative} associated to this operator, namely,
\begin{equation}
	\label{nota4}
	\partial_n u_{\vert_{\partial\Omega}}=
	{\bf n}\cdot P\nabla_{X,z}u_{\vert_{\partial\Omega}},
\end{equation}
${\bf n}$ standing for the \emph{outward} unit normal vector to 
$\partial\Omega$.

\section{Preliminary results}\label{sectprel}

\subsection{Commutator estimates and anisotropic Poisson regularization}
\label{sectcomm}

We recall first the tame product and Moser estimates in Sobolev spaces:
if $t_0>1$ and $s\geq 0$, then 
$\forall f\in H^{s}\cap H^{t_0}(\R^2), \forall g \in H^s(\R^2),$
\begin{equation}
	\label{tame}
	%\quad
	\vert fg\vert_{H^s}\lesssim
	\vert f\vert_{H^{t_0}}\vert g\vert_{H^s}
	+ \langle \vert f\vert_{H^{s}}
	\vert g\vert_{H^{t_0}}\rangle_{s> t_0}
\end{equation}
and, for all $F\in C^\infty(\R^n;\R^m)$ such that $F(0)=0$,
\begin{equation}\label{moser}
	\forall u\in H^s(\R^2)^n,\; \;
	F(u)\in H^s(\R^2)^m
	\; \; \mbox{ and }\; \;
	\vert F(u)\vert_{H^s}\leq C(\vert u\vert_\infty)\vert u\vert_{H^s}.
\end{equation}
In the next proposition, we give tame commutator estimates.
\begin{prop}[Ths. 3 and 6 of \cite{LannesJFA}]
	\label{propprel1}
	Let $t_0>1$ and $-t_0<r\leq t_0+1$.
	Then, for all $s\geq 0$, $f\in H^{t_0+1}\cap H^{s+r}(\R^2)$
	and $u\in H^{s+r-1}(\R^2)$,
	$$
	\big\vert
	[\Lambda^s,f]u\big\vert_{H^r}
	\lesssim \vert \nabla f\vert_{H^{t_0}}\vert u\vert_{H^{s+r-1}}
	+\left\langle  \vert \nabla f\vert_{H^{s+r-1}}\vert u\vert_{H^{t_0}}
	\right\rangle_{s> t_0+1-r},
	$$
	where we used the notation (\ref{nota3}).
\end{prop}
One can deduce from the above proposition some commutator estimates useful
in the present study.
\begin{coro}
	\label{propprel2}
	Let $t_0>1$, $s\geq 0$ and $\gamma\in (0,1]$. Then:\\ 
	{\bf i.} For all $\bbv\in H^{t_0+2}\cap H^{s+1}(\R^2)^2$
	and $u\in  H^{s}(\R^2)$, one has
	$$
	\big\vert\big[\Lambda^s,\divg(\bbv\cdot)\big]u\big\vert_2
	\leq \vert \bbv\vert_{H^{t_0+2}}\vert u\vert_{H^s}
	+\big\langle
	\vert u\vert_{H^{t_0+1}}\vert \bbv \vert_{H^{s+1}}
	\big\rangle_{s> t_0+1}.
	$$
	{\bf ii.} For all $0\leq r\leq t_0+1$,
	$f\in L^\infty((-1,0);H^{s+r}\cap H^{t_0+1}(\R^2))$
	and $u\in L^2((-1,0);H^{s+r-1}(\R^2))$,
\begin{eqnarray*}
%$$
	\Vert \Lambda^r [\Lambda^s,f] u\big\Vert_{2}
	&\lesssim&	\Vert f \Vert_{L^\infty H^{t_0+1}}
	\Vert \Lambda^{s+r-1} u\Vert_{2} \\
	&&+\left\langle  \Vert f\Vert_{L^\infty H^{s+r}}
	\Vert \Lambda^{t_0}u\Vert_{2}
	\right\rangle_{s> t_0+1-r}.
%$$
\end{eqnarray*}
\end{coro}
\begin{proof}
	For the first point, just remark that
	$$
	\big[\Lambda^s,\divg(\bbv\cdot)\big]u
	=\big[\Lambda^s,\divg(\bbv)\big]u+
	\big[\Lambda^s,\bbv\big]\cdot \nag u,
	$$
and use Proposition \ref{propprel1} to obtain the result (recall
that $\gamma\leq 1$).\\
For the second point of the corollary, remark that for all $z\in [-1,0]$,
$$
	\vert \Lambda [\Lambda^s,f] u(z)\big\vert_{2}
	\lesssim	\vert f(z)\vert_{H^{t_0+1}}
	\vert u(z)\vert_{H^{s}}
	+\left\langle  \vert f(z)\vert_{H^{s+1}}\vert u(z)\vert_{H^{t_0}}
	\right\rangle_{s> t_0},
$$
as a consequence of Proposition \ref{propprel1} (with $r=1$). The corollary
then follows easily. \qed
\end{proof}

Let us end this section with a result on anisotropic Poisson regularization
(when $\gamma=\mu=1$, the result below is just the standard gain of half
a derivative of the Poisson regularization).
\begin{prop}
	\label{propsharp}
	Let $\gamma\in (0,1]$, $\mu>0$ and 
	$\chi$ be a smooth, compactly supported function and 
	$u\in \mafrS'(\R^2)$. Define also
	$u^\dag:=\chi(\sqrt{\mu}z\Dg)u$.\\
	For all $s\in\R$, if $u\in H^{s-1/2}(\R^2)$,
	one has $\Lambda^su^\dag\in L^2(\cS)$ and
	$$
	c_1\big\vert \frac{1}{(1+\sqrt{\mu}\Dg)^{1/2}}u\big\vert_{H^s}\leq 
	\Vert \Lambda^s u^\dag\Vert_{2}\leq c_2
	\big\vert \frac{1}{(1+\sqrt{\mu}\Dg)^{1/2}}u\big\vert_{H^s}.
	$$
	Moreover, for all $s\in\R$, if $u\in H^{s+1/2}(\R^2)$,
	one has $\Lambda^s\nagm u^\dag\in L^2(\cS)^3$ and
	$$
	c_1'\big\vert \frac{\sqrt{\mu}\Dg}{(1+\sqrt{\mu}\Dg)^{1/2}}u\big\vert_{H^s}\leq 
	\Vert \Lambda^s\nagm u^\dag\Vert_{2}\leq c_2'
	\big\vert\frac{\sqrt{\mu}\Dg}{(1+\sqrt{\mu}\Dg)^{1/2}}u\big\vert_{H^s}.
	$$
	In the above estimates,  $c_1$, $c_2$, $c_1'$ and $c_2'$  
	are nonnegative constants which depend 
	only on $\chi$.
\end{prop}
\begin{proof}
Write classically (with $\vert \xi^\gamma\vert=\sqrt{\xi_1^2+\gamma^2\xi_2^2}$),
\begin{eqnarray*}
	\Vert \chi(\sqrt{\mu}z\Dg)u\Vert_{s,0}^2&=&
	\int_{\R^2}\int_{-1}^0 \jap^{2s}\chi(\sqrt{\mu}z\vert \xi^\gamma\vert)^2\vert \widehat{u}(\xi)\vert^2 dz d\xi\\
	&=&\int_{\R^2}\jap^{2s}\frac{F(0)-F(-\sqrt{\mu}\xig)}{\sqrt{\mu}\xig}\vert\widehat{u}(\xi)\vert^2d\xi,
\end{eqnarray*}
where $F$ denotes a primitive of $\chi^2$. The first estimate of the
proposition then follows from
the observation that
$$
	c_1^2 \frac{1}{1+\sqrt{\mu}\xig}\leq
	\frac{F(0)-F(-\sqrt{\mu}\xig)}{\sqrt{\mu}\xig}
	\leq c_2^2 \frac{1}{1+\sqrt{\mu}\xig},
$$
where the constants depend only on $\chi$.\\
For the second estimate of the proposition, remark that 
$$
	\big\Vert \nagm u^\dag\big\Vert_2\sim 
	\sqrt{\mu}\big\Vert \chi(\sqrt{\mu}z\Dg)\Dg u\big\Vert_2
	+\sqrt{\mu}\big\Vert \chi'(z\sqrt{\mu}\Dg)\Dg u\big\Vert,
$$ 
and use the first part of the proposition.  \qed
\end{proof}

\subsection{Elliptic estimates on a strip}\label{sectell}

We recall that the velocity potential $\Phi$ within the fluid domain solves
the boundary value elliptic problem
\begin{equation}
	\label{ell1}
	\left\lbrace
	\begin{array}{l}
	\dz^2\Phi+\mu \dx^2\Phi+\gamma^2\mu \dy^2\Phi=0,
	\qquad
	-1+\beta b <z <\eps \zeta,\\
	\Phi_{\vert_{z=\eps\zeta}}=\psi,\qquad
	\partial_n \Phi_{\vert_{z=-1+\beta b}}=0,
	\end{array}\right.
\end{equation}
with $(\eps,\mu,\gamma,\beta)\in (0,1]\times 
(0,\infty)\times (0,1]\times [0,1]$.\\
Denote by ${\mathcal S}$ the flat strip ${\mathcal S}=\R^2\times (-1,0)$, and
assume that the following assumption is satisfied:
\begin{equation}
	\label{ell2}
	\mbox{There exists } h_0>0 \mbox{ such that }
	1+\eps\zeta-\beta b \geq h_0.
\end{equation}
Under this assumption, one can define a diffeomorphism  $S$ mapping 
${\mathcal S}$ onto 
the fluid domain $\Omega$:
$$
	S:
	\begin{array}{ccc}
	{\mathcal S}&\to &\Omega\\
	(X,z)&\mapsto & S(X,z):=\big(X,z+\sigma(X,z)\big),
	\end{array}
$$
with
\begin{equation}
	\label{ell2bis}
	\sigma(X,z)=-\beta z b(X)+\eps (z+1)\zeta(X).
\end{equation}

\begin{rema}
	The mapping $\sigma$ used in (\ref{ell2bis}) to define the
	diffeomorphism $S$ is the most simple one can think of. If
	one wanted to have optimal estimates with respect to
	the fluid or bottom parameterization (but unfortunately
	not uniform with respect to the parameters), one should use
	instead \emph{regularizing diffeormorphisms} as in
	Prop. 2.13 of \cite{LannesJAMS}. 
\end{rema}

From Proposition 2.7 of \cite{LannesJAMS}, we know that the BVP (\ref{ell1})
is equivalent to the BVP (recall that we use the convention (\ref{nota4}) for
normal derivatives),
\begin{equation}
	\label{ell3}
	\left\lbrace
	\begin{array}{l}
	\nabla_{X,z}\cdot P[\sigma]\nabla_{X,z}\phi=0,
	\qquad
	\mbox{ in }{\mathcal S},\\
	\phi_{\vert_{z=0}}=\psi,\qquad
	\partial_n \phi_{\vert_{z=-1}}=0,
	\end{array}\right.
\end{equation}
with $\phi=\Phi\circ S$ and with the $(2+1)\times(2+1)$ matrix 
$P[\sigma]$ given by
$$
	P[\sigma]:=P_{\mu, \gamma}[\sigma]
	=\left(\begin{array}{ccc}
	\mu(1+  \dz \sigma) &0 & -\mu   \dx \sigma\\
	0 & \gamma^2\mu(1+  \dz \sigma) & -\gamma^2 \mu   \dy \sigma\\
	-\mu   \dx \sigma& -\gamma^2 \mu   \dy \sigma & 
	\frac{1+\mu (\dx \sigma)^2
	+\gamma^2\mu  (\dy \sigma)^2}{1+  \dz \sigma}
	   \end{array}\right).
$$
Remark also that it follows from the expression of $P[\sigma]$ that
$$
	\nabla_{X,z}\cdot P[\sigma]\nabla_{X,z}=\nagm \cdot 
	(1+  Q[\sigma]) \nagm,
$$
where
\begin{equation}
	\label{ell4}
	Q[\sigma] :=
	Q_{\mu,\gamma}[\sigma]=\left(\begin{array}{ccc}
	\dz \sigma &0 & -\sqrt{\mu} \dx \sigma\\
	0 & \dz \sigma & -\gamma \sqrt{\mu} \dy \sigma\\
	-\sqrt{\mu} \dx \sigma& -\gamma \sqrt{\mu} \dy \sigma & 
	\frac{-\dz\sigma+\mu(\dx \sigma)^2
	+\gamma^2\mu(\dy \sigma)^2}{1+\dz \sigma}
	   \end{array}\right).
\end{equation}
Below we provide two important properties satisfied by $Q[\sigma]$.
\begin{prop}
	\label{propQ}
	Let $t_0>1$, $s\geq 0$, and $\zeta,b \in H^{t_0+1}\cap H^{s+1}(\R^2)$ 
	be such that 
	(\ref{ell2}) is satisfied. Assume also that $\sigma$ is as
	defined in (\ref{ell2bis}). Then:\\
	{\bf i.} One has
	$$
	\Vert Q[\sigma]\big\Vert_{L^\infty H^s}
	\leq
	C\big(\frac{1}{h_0},
	\Vert \nagm \sigma\Vert_{L^\infty  H^{t_0}}\big)
	\Vert\nagm\sigma\Vert_{L^\infty H^{s}}
	$$
	and, when $\sigma$ is also time dependent,
	$$
	\Vert \dt Q[\sigma]\big\Vert_{\infty,T}
	\leq
	C\big(\frac{1}{h_0},
	\Vert \nagm \sigma\Vert_{\infty,T}\big)
	\Vert\nagm\dt \sigma\Vert_{\infty,T}.
	$$
	{\bf ii.} For all $j\geq 1$ and 
	${\bf h}\in H^{t_0+1}\cap H^{s+1}(\R^2)^j$, and denoting by 
	$Q^{(j)}[\sigma]\cdot{\bf h}$ the $j$-th derivative of
	$\zeta\mapsto Q[\sigma]$ in the direction ${\bf h}$, one has
	\begin{eqnarray*}
	%\lefteqn
        &&\!\!\!\! {\Vert Q^{(j)}[\sigma]\cdot{\bf h}\big\Vert_{L^\infty H^s}
	\leq \big(\frac{\eps}{\nu}\big)^j
	C\big(\frac{1}{h_0},\eps\sqrt{\mu},
	\Vert \nagm \sigma\Vert_{L^\infty  H^{t_0}}\big)}\\
	&&\!\!\!\! \times
	\Big(\!\sum_{k=1}^j \vert h_k\vert_{H^{s+1}}\! \prod_{l\neq k}\vert h_l\vert_{H^{t_0+1}}
	+\big\langle (1+\Vert\nagm\sigma\Vert_{L^\infty H^{s}})\! \prod_{k=1}^j
	\vert h_k\vert_{H^{t_0+1}}\big\rangle_{s>t_0}\! \Big).
	\end{eqnarray*}
	{\bf iii.} The matrix $1+Q[\sigma]$ is coercive in the sense that
	$$
	\forall \Theta\in \R^{2+1},\qquad
	\vert\Theta\vert^2\lesssim k[\sigma] 
	 (1+Q[\sigma])\Theta\cdot\Theta,
	$$
	with
	$$
	k[\sigma]:=k_{\mu,\gamma}[\sigma]=
	 1+\Vert \dz\sigma\Vert_\infty+\frac{1}{h_0} 
  	\Big(1+\sqrt{\mu}\Vert \nag\sigma\Vert_{\infty}\Big)^2.
	$$
\end{prop}
\begin{proof}
The first two points follow directly from the tame product and Moser's 
estimate (\ref{tame}) and (\ref{moser}),
and the explicit expression of $Q[\sigma]$.\\
It is not difficult to see that 
 $(1+  Q[\sigma])\Theta \cdot \Theta 
 = \frac{1}{1+  \dz \sigma} |B \Theta|^2$, where
$$
	B=\left(\begin{array}{ccc}
	1+  \dz \sigma &0 & -   \sqrt{\mu} \dx \sigma\\
	0 & 1+  \dz \sigma & -\gamma   \sqrt{\mu} \dy \sigma\\
	0& 0 & 1
	   \end{array}\right).
$$
The matrix $B$ is invertible and its inverse is given by
$$
	B^{-1}=\frac{1}{1+ \dz \sigma}\left(\begin{array}{ccc}
	1& 0 & 
	  \sqrt{\mu}  \dx \sigma\\
	0 & 1 & 
	\gamma   \sqrt{\mu} \dy \sigma\\
	0& 0 & 1+ \dz\sigma
	   \end{array}\right).
$$
Remark now that owing to (\ref{ell2}),  the mapping $\sigma$, 
as given by (\ref{ell2bis}), satisfies
$(1+ \dz\sigma)^{-1}\leq h_0^{-1}$, so that
$$
  \sqrt{1+\dz\sigma}|B^{-1}|_{\mathbb{R}^3 \mapsto \mathbb{R}^3} 
  \lesssim \sqrt{1+\Vert\dz\sigma\Vert_\infty}+\frac{1}{\sqrt{h_0}} 
  \Big(1+  \sqrt{\mu}\Vert \nag\sigma\Vert_{\infty}\Big).
$$
Since $|B \Theta| |B^{-1}|_{\mathbb{R}^3 \mapsto \mathbb{R}^3} 
\ge |\Theta|$, the third claim of the proposition follows.  \qed
\end{proof}

Since the Dirichlet condition in (\ref{ell3}) can be 'lifted' 
in order to take
homogeneous Dirichlet boundary condition, we are led to 
study the following class of elliptic BVPs:
\begin{equation}
	\label{ell5}
	\left\lbrace
	\begin{array}{l}
	\nagm \cdot (1+  Q[\sigma])\nagm u=
	\nagm \cdot{\bf g},
	\qquad
	\mbox{ in }{\mathcal S},\\
	u_{\vert_{z=0}}=0,\qquad
	\partial_n u_{\vert_{z=-1}}=-{\bf e_z}\cdot {\bf g}_{\vert_{z=-1}},
	\end{array}\right.
\end{equation}
where, according to the notation (\ref{nota4}), $\partial_n u_{\vert_{z=-1}}$
stands for
$$
	\partial_n u_{\vert_{z=-1}}=-{\bf e_z}\cdot (1+  Q[\sigma]) 
	\nagm u_{\vert_{z=-1}}.
$$
Before stating the main result of this section let us introduce
a notation:
\begin{nota}
	We generically write
	\begin{equation}
	\label{eqM}
	M[\sigma]:=C\big(\eps\sqrt{\mu},\frac{1}{h_0},
	\Vert \nagm \sigma\Vert_{L^\infty H^{t_0+1}} \big),
	\end{equation}
	where, as usual, $C(\cdot)$ is a nondecreasing function of its
	arguments.
\end{nota}
\begin{prop} \label{prop:exist}
	Let $t_0>1$, $s\geq 0$ and $\zeta,b\in H^{t_0+2}\cap H^{s+1}(\R^2)$
	be such that (\ref{ell2}) is satisfied, and let $\sigma$ be given
	by (\ref{ell2bis}).\\
	Then for all ${\bf g}\in C([-1,0];H^s(\R^2)^3)$, there 
	 exists a unique variational solution 
	$u \in H^{1} (\mathcal{S})$ to the BVP (\ref{ell5})
	and
\begin{eqnarray*}
	\Vert \Lambda^s \nagm u\Vert_{2}
	\leq
	M[\sigma]
	\big(\Vert \Lambda^s {\bf g}\Vert_2
	+
	\big\langle
	\Vert \nagm \sigma\Vert_{L^\infty H^s}
	\big\Vert \Lambda^{t_0}{\bf g}\Vert_2
	\big\rangle_{s> t_0+1}
	\big),
\end{eqnarray*}
	where $M[\sigma]$
	is defined in (\ref{eqM}).
\end{prop}
\begin{proof}
The existence of the solution can be obtained with very classical tools and
we therefore omit it. We thus focus our attention on the proof of the estimate.\\
Let $\chi(\cdot)$ be a smooth, compactly supported function such that
$\chi(\xi)=1$ in a neighborhood of $\xi=0$, and define 
$\Lambda_h:=\Lambda*\chi(h D)$. Using $\Lambda^{2s}_h u$ as test function
in the variational formulation of (\ref{ell5}), one gets 
$$
	\int_{\mathcal S} (1+  Q[\sigma]) \nagm u\cdot 
	\nagm \Lambda^{2s}_h u
	=\int_{\mathcal S} {\bf g}\cdot \nagm \Lambda^{2s}_hu,
$$
so that using 
the fact that $\Lambda^s_h$ is $L^2$-self-adjoint, one gets, with
$v_h=\Lambda^s_h u$,
$$
	\int_{\mathcal S} \Lambda^s_h (1+  Q[\sigma])\nagm u\cdot 
	\nagm v_h
	=\int_{\mathcal S} \Lambda^s_h {\bf g}\cdot 
	\nagm v_h,
$$
and thus 
\begin{equation*}
	\! \int_{\mathcal S} \! (1+  Q[\sigma]) \nagm v_h \cdot
	\nagm v_h \!
	=\!\!  \int_{\mathcal S} \!\! 
        \big(\Lambda^s_h {\bf g}\cdot \nagm v_h-\big[\Lambda^s_h,   Q[\sigma]\big]
	\nagm u\cdot \nagm v_h\big).
\end{equation*}
Thanks to the coercitivity property of Proposition \ref{propQ}, 
one gets 
\begin{equation}
	\label{ellter}
	k[\sigma]^{-1}\Vert \Lambda^s_h \nagm u\Vert_{2}
	\lesssim
	\Vert \big[\Lambda^s_h,  Q[\sigma]\big]\nagm u\Vert_2
	+\Vert \Lambda^s_h {\bf g}\Vert_2;
\end{equation}
since $[\Lambda^s,Q[\sigma]]$ is of order $s-1$, 
the above estimates allows one to conclude, after letting $h$ go to zero,
that $\Lambda^s\nagm u\in L^2(\cS)$; more precisely,
thanks to Corollary \ref{propprel2}, one deduces
\begin{eqnarray*}
	k[\sigma]^{-1}\Vert \Lambda^s \nagm u\Vert_{2}
	&\lesssim& 
	\Vert \Lambda^s {\bf g}\Vert_2
	+
	\big\Vert Q[\sigma]
	\big\Vert_{L^\infty H^{t_0+1}}
	\big\Vert \Lambda^{s-1} \nagm u\Vert_{2}\\
	& &+  \big\langle 
	\Vert Q[\sigma]\Vert_{L^\infty H^s}
	\big\Vert \Lambda^{t_0}\nagm u \big\Vert_{2}
	\big\rangle_{s> t_0+1},
\end{eqnarray*}
and thus
\begin{eqnarray}
	\label{eqell}
	\lefteqn{\Vert \Lambda^s \nagm u\Vert_{2}
	\leq
	C\big(k[\sigma],
	\big\Vert Q[\sigma]\big\Vert_{L^\infty H^{t_0+1}} \big)}\\
	\nonumber
	&\times&
	\big(\Vert \Lambda^s {\bf g}\Vert_2
	+
	\Vert \nagm u\Vert_{2}+
	\big\langle 
	\Vert Q[\sigma]\Vert_{L^\infty H^s}
	\big\Vert \Lambda^{t_0}\nagm u \big\Vert_{2}
	\big\rangle_{s> t_0+1}
	\big).
\end{eqnarray}
One also gets $\Vert \nagm u\Vert_2\leq k[\sigma] \Vert {\bf g} 
\Vert_2$ from (\ref{ellter}) after
remarking that the commutator in the r.h.s. vanishes when $s=h=0$; taking
$s=t_0$ in (\ref{eqell}) then gives
$\Vert\Lambda^{t_0} \nagm u\Vert_2\leq C\big(k[\sigma],
	\big\Vert Q[\sigma]\big\Vert_{L^\infty H^{t_0+1}} \big)\Vert \Lambda^{t_0}{\bf g}\Vert_2$, so that the r.h.s. of 
(\ref{eqell}) is bounded from above by
\begin{eqnarray*}
	C\big(k[\sigma],
	\Vert Q[\sigma]\Vert_{L^\infty H^{t_0+1}} \big)
	\big(\Vert \Lambda^s {\bf g}\Vert_2
	+
	\big\langle
	\Vert Q[\sigma]\Vert_{L^\infty H^s}
	\big\Vert \Lambda^{t_0}{\bf g}\Vert_2
	\big\rangle_{s> t_0+1}
	\big).
\end{eqnarray*}
The proposition follows therefore from Proposition \ref{propQ}. \qed
\end{proof}
Before stating a corollary to Proposition \ref{prop:exist}, let us
introduce a few notations:
\begin{nota}
    \label{notaDNoth1}
    {\bf i.} For all $u\in H^{3/2}(\R^2)$, we define $u^\flat$ as the
    solution to the BVP
    \begin{equation}
	\label{eqnota}
        \left\lbrace
        \begin{array}{l}
        \nagm \cdot (1+ Q[\sigma])\nagm u^\flat=0\\
        u^\flat _\surf=u,\qquad
        \partial_n u^\flat_\fond=0.
        \end{array}\right.
    \end{equation}
    {\bf ii.} For all $u\in \mafrS'(\R^2)$, one defines $u^\dag$ as
    $$
    \forall z\in [-1,0], u^\dag(\cdot,z)=\chi(\sqrt{\mu}z\Dg)u,
    $$
    where $\chi$ is a smooth, compactly supported function such that
    $\chi(0)=1$.
\end{nota}
The following corollary gives some control on the extension mapping
$u\mapsto u^\flat$.
\begin{coro}
	\label{coro1}
	Let $t_0>1$ and $s\geq 0$. Let also 
	$\zeta,b\in H^{t_0+2}\cap H^{s+1}(\R^2)$
	be such that (\ref{ell2}) is satisfied, and $\sigma$ be given
	by (\ref{ell2bis}).\\
	Then for all $u\in H^{s+1/2}(\R^2)$, there 
	 exists a unique solution 
	$u^\flat \in {H}^{1} (\mathcal{S})$ and
	$$
	\Vert \Lambda^s\nagm u^\flat\Vert_2
	\! \leq	\! \sqrt{\mu\nu}M[\sigma]
	\big(\vert \Pig u\vert_{H^s}
	+
	\big\langle 
	\Vert \nagm \sigma\Vert_{\! L^\infty H^s}
	\vert \Pig
	u\vert_{H^{t_0}} \!
	\big\rangle_{\! \! s> t_0+1} \!
	\big),
	$$
	with $\Pig$ as defined in (\ref{Pig}).
\end{coro}
\begin{proof}
Looking for $u^\flat$ under the form $u^\flat=v+u^\dag$, with 
$u^\dag$ given by Notation \ref{notaDNoth1}, one must solve
\begin{equation}
	\label{demcor1}
	\left\lbrace
	\begin{array}{l}
	\nagm\cdot \matr\nagm v=-\nagm\cdot \matr\nagm u^\dag,\\
	v\surf=0,\qquad \partial_n v\fond=
	{\bf e_z}\cdot \matr\nagm u^\dag\fond.
	\end{array}\right.
\end{equation}
Applying Proposition \ref{prop:exist} (with ${\bf g}=-\matr\nagm u^\dag$),
one gets
\begin{eqnarray*}
	\Vert \Lambda^s \nagm \! v\Vert_2
	\! \leq \!\! M[\sigma]
	\big(\Vert \Lambda^s \nagm \! u^\dag  \Vert_2
	\! + \! 
	\big\langle
	\Vert \nagm \! \sigma\Vert_{\! L^\infty H^s} \!
	\big\Vert \Lambda^{t_0} \! \nagm \! u^\dag\Vert_2
	\big\rangle_{\! s> t_0+1}
	\big),
\end{eqnarray*}
and since $u^\flat=u^\dag+v$, the corollary follows from Proposition \ref{propsharp}.
\qed
\end{proof}
\begin{rema}
	\label{remaprov}
	From the variational formulation of (\ref{demcor1}),
	one gets easily
$$
    \Vert (1+Q[\sigma])^{1/2}\nagm v\Vert_2
    \leq
    \Vert (1+Q[\sigma])^{1/2}\nagm u^\dag\Vert_2.
$$
\end{rema}

\section{The Dirichlet-Neumann operator}\label{sectDN}

As seen in the introduction, we define the Dirichlet-Neumann operator 
$\G\cdot$ as
$$
	\G \psi=\sqrt{1+\vert \eps \nabla\zeta\vert^2}
	\partial_n\Phi_{\vert_{z=\eps\zeta}},
$$
where $\Phi$ solves (\ref{ell1}).
Using Notation \ref{notaDNoth1}, one can give an alternate definition
of $\G\cdot$ (see Proposition 3.4 of \cite{LannesJAMS}), namely,
$$
	\G \psi =\partial_n \psi^\flat\surf
	\quad\big(={\bf e_z}\cdot P[\sigma]\nabla_{X,z}\psi^\flat\surf\big).
$$
More, precisely one has:
\begin{prop}
	Let $t_0>1$, $s\geq 0$ and $\zeta,b\in H^{t_0+2}\cap H^{s+1}(\R^2)$
	be such that (\ref{ell2}) is satisfied, and let $\sigma$ be given
	by (\ref{ell2bis}).\\
	Then one can define the mapping
	$\G\cdot$ (or simply $\GG\cdot$ when no confusion is possible) as
	$$
	\G\,(=\GG)\,:\begin{array}{ccc}
	H^{s+1/2}(\R^2)&\to& H^{s-1/2}(\R^2)\\
	u &\mapsto& \partial_n u^\flat\surf
	   \end{array}.
	$$
\end{prop}
\begin{proof}
The extension $u^\flat$ is well-defined owing to Corollary \ref{coro1}.
Moreover, we can use the definition of $P[\sigma]$ and $Q[\sigma]$ to see that
$$
	{\bf e_z}\cdot P[\sigma]\nabla_{X,z}u^\flat={\bf e_z}\cdot (1+Q[\sigma])
	\nagm u^\flat.
$$
We will now show that it makes sense to take the trace of the above
expression at $z=0$. This is trivially true for $Q[\sigma]$, so that
we are left with $u^\flat$. After a brief look at the proof of Corollary
\ref{coro1}, and using the same notations, one gets
$u^\flat=v+u^\dag$. Since one obviously has 
$\nagm u^\dag \in C([-1,0];H^{s-1/2}(\R^2)^3)$, 
the trace $\nagm u^\dag\surf$
makes sense. In order to prove that $\nagm v\surf$ is also defined, 
remark that $v$, which solves (\ref{demcor1}) satisfies
$\nagm v\in L^2((-1,0);H^{s}(\R^2)^3)$ and, using the equation, 
$\dz \nagm v\in L^2((-1,0);H^{s-1}(\R^2)^3)$. By the trace theorem,
these two properties show that $\nagm v\surf \in H^{s-1/2}(\R^2)^3$. \qed
\end{proof}

\subsection{Fundamental properties}\label{sectfund}

We begin this section with two basic properties of the Dirichlet-Neumann
operator which play a key role in the energy estimates.
\begin{prop}
	\label{propDNoth1}
	Let $t_0>1$ and $\zeta,b\in H^{t_0+2}(\R^2)$
	be such that (\ref{ell2}) is satisfied. Then \\
	{\bf i.} The Dirichlet-Neumann operator is self-adjoint:
	$$
	\forall u,v\in H^{1/2}(\R^2),\qquad
	(u,\GG v)=(v,\GG u).
	$$
	{\bf ii.} One has
	$$
	\forall u,v\in H^{1/2}(\R^2),\; \; 
	\big\vert (u,\GG v)\big\vert\leq
	(u,\GG u)^{1/2}
	(v,\GG v)^{1/2}.
	$$
\end{prop}
\begin{proof}
Using Notation \ref{notaDNoth1},  one gets by Green's identity
that
\begin{eqnarray}
	\label{DNoth0}
	(u,\GG v)&=&\int_\cS \! (1+Q[\sigma])\nagm u^\flat\cdot
	\nagm v^\flat \\
	\label{DNoth1}
	&=& \int_\cS  \! (1+Q[\sigma])^{1/2}\nagm u^\flat \! \cdot \!
	 (1+Q[\sigma])^{1/2}\nagm v^\flat,
\end{eqnarray}
where $(1+Q[\sigma])^{1/2}$ stands for the square root of the positive
definite matrix $(1+Q[\sigma])$ (note that the symmetry in $u$ and $v$ of 
the above expression proves the --very classical-- first point of the
proposition). It follows therefore from Cauchy-Schwartz inequality that
$$
	(u,\GG v)\leq \big\Vert (1+Q[\sigma])^{1/2}\nagm u^\flat\big\Vert_2
	\,
	 \big\Vert (1+Q[\sigma])^{1/2}\nagm v^\flat\big\Vert_2,
$$
which yields the second point of the proposition, since one has
\begin{equation}
	\label{DNoth2}
	(u,\GG u)=\big\Vert (1+Q[\sigma])^{1/2}\nagm u^\flat\big\Vert_2^2
\end{equation}
(just take $u=v$ in (\ref{DNoth1})). \qed
\end{proof}

The next proposition is related to the variational formula of Hadamard and
gives a uniform control of the operator norm of the
DN operator and its derivatives 
(recall that we use the convention (\ref{conv1})
and that $d_\zeta^j\cG[\eps\cdot]u\cdot{\bf h}=\cG[\eps\zeta]u$ when $j=0$).
\begin{prop}
	\label{unifcontrol}
	Let $t_0>1$, $s\geq 0$ and $\zeta,b\in H^{t_0+2}\cap H^{s+1}(\R^2)$
	be such that (\ref{ell2}) is satisfied, and let $\sigma$ be given
	by (\ref{ell2bis}).\\
	For all $u\in H^{s+1/2}(\R^2)$, $j\in\N$ and
	${\bf h}\in H^{t_0+2}\cap H^{s+1}(\R^2)^j$, one has
	\begin{eqnarray*}
        \lefteqn{ \!\!\!\!\!\!\!\!\!\!\!
        \big\vert 
	\frac{1}{\sqrt{\mu}}d_\zeta^j\cG[\eps\cdot]u\cdot{\bf h}
	\big\vert_{H^{s-1/2}}\leq
	(\frac{\eps}{\nu})^jM[\sigma]
	\Big(
	\big\vert \Pig u
	\big\vert_{H^{s}}  \prod_{k=1}^j\vert h_k\vert_{H^{t_0+1}}}\\
	& & \!\!\!\!\!\!\!\! 
        +\big\langle (1+\Vert \nagm\sigma\Vert_{L^\infty H^{s}}) 
	\big\vert \Pig u
	\big\vert_{H^{t_0}}  \prod_{k=1}^j\vert h_k\vert_{H^{t_0+1}}
	\big\rangle_{s>t_0}\\
	& & \!\!\!\!\!\!\!\! 
        +
        \big\langle
	\big\vert \Pig u
	\big\vert_{H^{t_0}}\sum_{k=1}^j\vert h_k\vert_{H^{s+1}}
	\prod_{l\neq k}\vert h_l\vert_{H^{t_0+1}}
	\big\rangle_{s>t_0}
	\Big),
	\end{eqnarray*}
	with $M[\sigma]$ as in (\ref{eqM}) while
	$\Pig$ is defined in (\ref{Pig}).
\end{prop}
\begin{rema}
	When $j=0$, the proposition gives 
	a much more precise estimate on 
	$\vert \cG[\eps \zeta] u \vert_{H^{s-1/2}}$ than 
	Theorem 3.6 of \cite{LannesJAMS}, but  requires
	$\zeta\in H^{s+1}(\R^2)$ while
	$\zeta\in H^{s+1/2}(\R^2)$ is enough, as shown in
	\cite{LannesJAMS} through the use of regularizing diffeomorphisms.
	This lack of optimality in the $\zeta$-dependence is the price
	to pay to obtain uniform estimates in terms of
	a $\cE^s(\cdot)$ rather than Sobolev-type norm. 
\end{rema}
\begin{rema}	\label{remopnorm}
	The r.h.s. of the estimate given in the proposition
	(when $j=0$)
	is itself bounded from above by
	$$
	M[\sigma]
	\big(
	\vert u \vert_{H^{s+1}}
	+\big\langle
	\Vert \nagm \sigma\Vert_{L^\infty H^s}
	\vert u \big\vert_{H^{t_0+1}}
	\big\rangle_{s> t_0}\big).
	$$
\end{rema}
\begin{proof}
First remark that
one has $\Lambda^{s-1/2} v^\dag\surf=\Lambda^{s-1/2} v$ 
(with $v^\dag$ as in
Notation \ref{notaDNoth1}), so that one gets 
by Green's identity,
\begin{eqnarray}
	\nonumber
	(\Lambda^{s-1/2}\cG[\eps\zeta] u,v)&=&
	(\cG[\eps\zeta] u,\Lambda^{s-1/2} v)\\
	\nonumber
	&=&\int_\cS (1+Q[\sigma])\nagm u^\flat\cdot
	\Lambda^{s-1/2}\nagm v^\dag\\
	\label{green1}
	&=&\int_\cS \! \Lambda^{s}
	(1+Q[\sigma])\nagm u^\flat\! \cdot \! 
	\Lambda^{-1/2}\nagm v^\dag.
\end{eqnarray}
A Cauchy-Schwartz inequality then yields,
\begin{equation}\label{bof}
	(\Lambda^{s-1/2}\GG u,v)\leq
	\Vert \Lambda^{s}
	(1+Q[\sigma])\nagm u^\flat\Vert_2
	\Vert \Lambda^{-1/2}\nagm v^\dag\Vert_2,
\end{equation}
and since it follows from the product estimate (\ref{tame})
that $\Vert \Lambda^{s}
	(1+Q[\sigma])\nagm u^\flat\Vert_2$ is bounded from above by
$$
	(1+\Vert Q[\sigma]\Vert_{L^\infty H^{t_0}})
	\Vert \Lambda^{s}\nagm u^\flat\Vert_2
	+\big\langle\Vert Q[\sigma]\Vert_{L^\infty H^{s}}
	\Vert \Lambda^{t_0}\nagm u^\flat\Vert_2\big\rangle_{s> t_0},
$$
one can deduce from Propositions \ref{propsharp} and \ref{propQ} that
(recall that $\nu=\frac{1}{1+\sqrt{\mu}}$),
\begin{eqnarray*}
	\lefteqn{(\Lambda^{s-1/2}\GG u,v)\leq \nu^{-1/2}
	M[\sigma]
	\vert v\vert_2}\\
	&\times&\Big(\Vert \Lambda^{s}\nagm u^\flat\Vert_2
	+\big\langle\Vert \nagm \sigma \Vert_{L^\infty H^{s}}
	\Vert\Lambda^{t_0}\nagm u^\flat\Vert_2\big\rangle_{s> t_0}\Big),
\end{eqnarray*}
and  the 
proposition thus follows directly from Corollary \ref{coro1} and
a duality argument in the case $j=0$.\\
In the case $j\neq 0$, after differentiating (\ref{green1}),
and using the same notation as in Proposition \ref{propQ}, one gets
\begin{eqnarray}	
	\nonumber
	(\Lambda^{s-1/2}d_\zeta^j\cG[\eps\cdot] u\cdot{\bf h},v)
	=
	\int_\cS \Lambda^{s}
	(Q^{(j)}[\sigma]\cdot{\bf h})\nagm u^\flat\cdot
	\Lambda^{-1/2}\nagm v^\dag\\
	\label{eqsupp}
	+
	\sum_{k=1}^j\,\sum_{{\bf h}_k,{\bf h}_{j-k}}
	\int_\cS\Lambda^sB({\bf h}_k,{\bf h}_{j-k})\cdot \Lambda^{-1/2}\nagm v^\dag,\qquad
\end{eqnarray}
where the second summation is over all the $k$-uplets ${\bf h}_k$ and
$(j-k)$-uplets ${\bf h}_{j-k}$ such that $({\bf h}_k,{\bf h}_{j-k})$
is a permutation of ${\bf h}$, and where $B({\bf h}_k,{\bf h}_{j-k})$
is given by
$$
	B({\bf h}_k,{\bf h}_{j-k})=(Q^{(j-k)}[\sigma]\cdot {\bf h}_{j-k})
	\nagm (u^{\flat,k}\cdot {\bf h}_k)
$$
($u^{\flat,k}\cdot {\bf h}_k$ standing for the $k$-th order derivative
of $\zeta\mapsto u^\flat$ at $\uz$ and in the direction ${\bf h}_k$).\\
Proceeding as for the case $j=0$ and using the estimates on 
$\Vert Q^{(j)}[\sigma]\cdot{\bf h}\Vert_{L^\infty H^s}$
provided by Proposition \ref{propQ}, one arrives at the desired estimate
for the first term of the r.h.s. of (\ref{eqsupp}). For the other terms,
one has to remark first that $u^{\flat,k}\cdot {\bf h}_k$ 
solves a bvp like (\ref{ell5}) with
$$
	{\bf g}=-\sum_{l=0}^{k-1} \,\sum_{{\bf h}_{k,l},{\bf h}_{k,k-l}}
	B({\bf h}_{k,l},{\bf h}_{k,k-l}),
$$
where the second summation is taken over all the $l$ and $k-l$-uplets
such that $({\bf h}_{k,l},{\bf h}_{k,k-l})$ is a permutation of ${\bf h}_k$.
A control of 
$\Vert \nagm u^{\flat,k}\cdot{\bf h}_k\Vert_{L^\infty H^s}$ 
 in terms of 
$\Vert Q^{(k-l)}[\sigma]\cdot{\bf h}_{k,k-l}\Vert_{L^\infty H^s}$
 and $\Vert \nagm u^{\flat,l}\cdot{\bf h}_{k,l}\Vert_{L^\infty H^s}$
is therefore provided by Proposition \ref{prop:exist}. It is then easy
to conclude by a simple induction.
\qed
\end{proof}
\begin{rema}\label{remlemfin}
	Instead of (\ref{bof}), one can easily get
	$$
	(\Lambda^{s-1/2}\GG u,v)\leq
	\Vert \Lambda^{s+1/2}
	(1+Q[\sigma])\nagm u^\flat\Vert_2
	\Vert \Lambda^{-1}\nagm v^\dag\Vert_2,
	$$
	and since 
	$\Vert \Lambda^{-1}\nagm v^\dag\Vert_2\lesssim 
	\sqrt{\mu}\vert v\vert_2$ one also has the estimate
	(with $\zeta=0$ for the sake of simplicity):
	$$
	\big\vert \frac{1}{\mu}\cG[0]\psi\big\vert_{H^{s-1/2}}
	\leq C(\frac{1}{h_0},\beta\sqrt{\mu},\vert b\vert_{H^{s+3/2}})
	\big\vert \frac{\Dg}{(1+\sqrt{\mu}\Dg)^{1/2}}\psi\big\vert_{H^{s+1/2}},
	$$
	showing that $\frac{1}{\mu}\cG[0]\cdot$ 
	can be uniformly controlled		when $\mu$ goes to zero.
\end{rema}

The proposition below show that controls in terms of $\vert\Pig u\vert_2$
or $(u,\frac{1}{\mu\nu}\GG u)^{1/2}$  are equivalent. This result can 
be seen as a version of the Garding inequality for the DN operator.
\begin{prop}
	\label{propDNoth2}
	Let $t_0>1$ and $\zeta,b\in H^{t_0+2}(\R^2)$
	be such that (\ref{ell2}) is satisfied, and let $\sigma$ be given
	by (\ref{ell2bis}), $k[\sigma]$ be as defined in Proposition
	\ref{propQ} and $\Pig$ be given by (\ref{Pig}).
	For all $u\in H^{1/2}(\R^2)$, one has
	$$
	(u,\frac{1}{\mu\nu}\GG u)\leq
	M[\sigma]
	\vert \Pig u \vert_2^2
	\quad\mbox{ and }\quad
	k[\sigma]^{-1}
	\vert \Pig u \vert_2^2
	\lesssim
	(u,\frac{1}{\mu\nu}\GG u).
	$$
\end{prop}
\begin{proof}
The first estimate of the proposition  follows directly 
from (\ref{DNoth2}) and Corollary  \ref{coro1}. 

The second estimate is more delicate. Let $\varphi$ be a smooth function,
with compact support in $(-1,0]$ and such that $\varphi(0)=1$; define also
$v(X,z)=\varphi(z)u^\flat$
(with $u^\flat$ defined as in Notation \ref{notaDNoth1}).
Since $v_\fond=0$, one can get, after taking the Fourier transform
with respect to the horizontal variables,
$$
	\frac{\xig^2}{1+\sqrt{\mu}\xig}\vert\widehat{u}(\xi)\vert^2
	\leq 2
	\int_{-1}^0
	\frac{\xig^2}{1+\sqrt{\mu}\xig}\vert\widehat{v}(\xi,z)\vert
	\,\vert\dz\widehat{v}(\xi,z)\vert dz.
$$
Remarking that 
$$
	\vert\widehat{v}\vert\leq \vert\varphi\vert_\infty 
	\vert\widehat{u^\flat}\vert
	\quad\mbox{ and }\quad
	\vert\dz \widehat{v}\vert\leq \vert\dz\varphi\vert_\infty\vert\widehat{u^\flat}\vert
	+\vert\varphi\vert_\infty 
	\vert\dz\widehat{u^\flat}\vert,
$$
one gets
\begin{eqnarray*}
	\lefteqn{\frac{\xig^2}{1+\sqrt{\mu}\xig}\vert\widehat{u}(\xi)\vert^2
	\leq 2\vert\varphi\vert_\infty\vert\dz\varphi\vert_\infty
	\int_{-1}^0
	\frac{\xig^2}{1+\sqrt{\mu}\xig}\vert\widehat{u^\flat}(\xi,z)\vert^2
	dz}\\
	& &+2\vert\varphi\vert^2_\infty \int_{-1}^0
	\frac{\xig^2}{1+\sqrt{\mu}\xig}\vert\widehat{u^\flat}(\xi,z)\vert
	\,\vert\dz\widehat{u^\flat}(\xi,z)\vert dz,\\
	&\leq& 2\vert\varphi\vert_\infty\vert\dz\varphi\vert_\infty
	\int_{-1}^0
	\xig^2\vert\widehat{u^\flat}(\xi,z)\vert^2
	dz\\
	& &+\vert\varphi\vert^2_\infty \int_{-1}^0
	\frac{\mu\xig^4}{(1+\sqrt{\mu}\xig)^2}
	\vert\widehat{u^\flat}(\xi,z)\vert^2dz
	+ \vert\varphi\vert^2_\infty \int_{-1}^0
	\frac{1}{\mu}\vert\dz\widehat{u^\flat}(\xi,z)\vert^2 dz,
\end{eqnarray*}
where Young's inequality has been used to obtain the last line. Remarking
now that
$\frac{\mu\xig^4}{(1+\sqrt{\mu}\xig)^2}\leq \xig^2$, one has
$$
	\frac{\xig^2}{1+\sqrt{\mu}\xig}\vert\widehat{u}(\xi)\vert^2
	\lesssim
	\int_{-1}^0
	\xig^2\vert\widehat{u^\flat}(\xi,z)\vert^2
	dz+
	\int_{-1}^0
	\frac{1}{\mu}\vert\dz\widehat{u^\flat}(\xi,z)\vert^2 dz,
$$	
so that, integrating with respect to $\xi$, one gets
$$
	\Big\vert \frac{\Dg}{(1+\sqrt{\mu}\Dg)^{1/2}}u\Big\vert_2^2
	\lesssim
	\frac{1}{\mu}\Vert \nagm u^\flat\Vert_2^2.
$$
Owing to Proposition \ref{propQ} and (\ref{DNoth0}) (with $v=u$), one has 
$$
	\Vert \nagm u^\flat\Vert_2^2\lesssim k[\sigma] (u,\GG u),
$$
 and the proposition follows. \qed
\end{proof}

\subsection{Commutator estimates}\label{sectcomDN}

In the following proposition, we show how to control in the energy estimates, 
the terms involving commutators between the Dirichlet-Neumann operator
and spatial or time derivatives and in terms of $\cE^s(\cdot)$ rather 
than Sobolev-type norms.
\begin{prop}
	\label{propestcom}
	Let $t_0>1$, $s\geq 0$ and $\zeta,b\in H^{t_0+2}\cap H^{s+2}(\R^2)$
	be such that (\ref{ell2}) is satisfied, and let $\sigma$ be given
	by (\ref{ell2bis}). Then, for all $v\in H^{s+1/2}(\R^2)$,
	\begin{eqnarray*}
	\big\vert
	[\Lambda^s,\frac{1}{\mu\nu}\GG] v
	\big\vert_2
	&\leq &
	M[\sigma]
		\Big(\Vert\nagm\sigma\Vert_{L^\infty H^{t_0+1}}
	\vert \Pig v \vert_{H^s}\\
	&+&\big\langle 
	\Vert\nagm\sigma\Vert_{L^\infty H^{s+1}}
	\vert \Pig v \vert_{H^{t_0}}
	\big\rangle_{s> t_0}\Big),
\end{eqnarray*}
where $M[\sigma]$ is as in (\ref{eqM}) while $\Pig$ is defined in (\ref{Pig}).
\end{prop}
\begin{proof}
First remark that for all $u\in \mafrS(\R^2)$, 
$$
	\big(u,[\GG,\Lambda^s]v\big)
	=\big(u,\GG \Lambda^s v\big)
	-\big(\Lambda^s u,\GG v\big).
$$
Since $u^\dag\surf=u$ and 
$\Lambda^s u^\dag\surf=\Lambda^s u$ 
(we use here Notation \ref{notaDNoth1}), 
it follows from Green's identity that
\begin{eqnarray}
	\nonumber
	\lefteqn{\big(u,[\GG,\Lambda^s]v\big)
	=\int_\cS \nagm u^\dag \cdot 	\matr
	\nagm (\Lambda^s v)^\flat}\\
	\nonumber
	& &	-\int_\cS \matr \nagm v^\flat\cdot\nagm \Lambda^s 
	u^\dag\\
	\label{eqcom1}
	&=&\int_\cS \nagm u^\dag\cdot
	\big(\matr \nagm \big((\Lambda^s v)^\flat-\Lambda^s v^\flat\big)
	-[\Lambda^s,Q[\sigma]]\nagm v^\flat\big). \nonumber \\
\end{eqnarray}
Let us now prove the following lemma:
\begin{lemm}
	For all $f\in L^2(\R^2)$ and ${\bf g}\in H^1(\cS)^3$, one has
		$$
		\int_{\cS} \nagm f^\dag \cdot {\bf g}\lesssim \sqrt{\mu }
		\sqrt{\nu}\vert f \vert_2\big\Vert \Lambda{\bf g}\big\Vert_2.
		$$
\end{lemm}
\begin{proof}
By definition of $f^\dag$, one has
$$
	\nagm f^\dag \! = \!
	\sqrt{\mu}\left(\begin{array}{c}
	\chi(z\sqrt{\mu}\Dg)\dx f\\
	\gamma\chi(z\sqrt{\mu}\Dg)\dy f\\
	\chi'(z\sqrt{\mu}\Dg)\Dg f
			\end{array}\right).
$$
	Replacing $\nagm f^\dag$ in the integral to control by
	 this expression, and using the self-adjointness of $\Lambda$,
	one gets easily
	from Proposition \ref{propsharp} that
$$
	\int_{\cS} \nagm f^\dag \cdot {\bf g}\lesssim \sqrt{\mu }\sqrt{\nu}
	\vert \Pig \Lambda^{-1}f\vert_2
	\big\Vert \Lambda {\bf g}\big\Vert_2;
$$
recalling that $\nu=\frac{1}{1+\sqrt{\mu}}$ and $\gamma\leq 1$, one can check
that $\vert \Pig \Lambda^{-1}f\vert_2
\lesssim \vert f\vert_2$, \emph{uniformly with respect to $\mu$ and $\gamma$},
and the lemma follows. \qed
\end{proof}
It is then a simple consequence of the lemma, (\ref{eqcom1}) and 
Proposition \ref{propQ} that
\begin{eqnarray}
	\nonumber
	\lefteqn{\big(u,[\GG,\Lambda^s]v\big)
	\lesssim \sqrt{\mu}\sqrt{\nu}
	\vert u\vert_2
	\Big(	\Vert\Lambda [\Lambda^s,Q[\sigma]]\nagm v^\flat\Vert_{2}}\\
	\label{eqcom2}
	&+ &(1+\Vert \nagm\sigma \Vert_{L^\infty H^{t_0+1}})
	\Vert \Lambda \nagm\big((\Lambda^s v)^\flat
	-\Lambda^s v^\flat\big)\Vert_{2}
	\Big),
\end{eqnarray}
which motivates the following lemma:
\begin{lemm}
	\label{lemmg}
	One has
	$$
	\Vert \Lambda\nagm \big((\Lambda^s v)^\flat-\Lambda^s v^\flat\big)
	\Vert_2
	\leq
	M[\sigma]
	\Vert \Lambda[\Lambda^s,Q[\sigma]]\nagm v^\flat\Vert_2.
	$$
\end{lemm}
\begin{proof}
Just remark that $w:=(\Lambda^s v)^\flat-\Lambda^s v^\flat$ solves
$$
	\left\lbrace\begin{array}{l}
	\nagm\cdot \matr \nagm w=\nagm\cdot {\bf g},\\
	w\surf=0,\qquad \partial_n w\fond=-{\bf e_z}\cdot {\bf g}\fond,
	\end{array}\right.
$$
with ${\bf g}=[\Lambda^s,Q[\sigma]]\nagm v^\flat$, and use 
Proposition \ref{prop:exist}. \qed
\end{proof}
With the help of the lemma, one deduces from (\ref{eqcom2}) 
that 
$$
	\big(u,[\GG,\Lambda^s]v\big)
	\leq
	\sqrt{\mu}\sqrt{\nu}M[\sigma]
	\big\Vert\Lambda [\Lambda^s,Q[\sigma]]\nagm v^\flat\big\Vert_{2}
	\vert u\vert_2,
$$
and thus, owing to Corollary \ref{propprel2} and Proposition \ref{propQ}
\begin{eqnarray*}
	\!\! \lefteqn{\big(u,[\GG,\Lambda^s]v\big)
	\leq \sqrt{\mu}\sqrt{\nu} M[\sigma]
	\vert u\vert_2}\\
	\!\! &\times& \!\! 
	\big(\Vert \nagm\! \sigma\Vert_{L^\infty H^{t_0+1}}
	\Vert \Lambda^s \nagm \! v^\flat\Vert_{2}
	\! + \! \big\langle 
	\Vert\nagm \! \sigma\Vert_{L^\infty H^{s+1}}
	\Vert \Lambda^{t_0}\nagm \! v^\flat\Vert_2
	\big\rangle_{\! s> t_0}\big),
\end{eqnarray*}	
and the result follows therefore from Corollary \ref{coro1} and a
duality argument. \qed
\end{proof}

The next proposition gives control of the commutator between the
Dirichlet-Neumann operator and a time derivative.
\begin{prop}
	\label{time}
	Let $t_0>1$, $T>0$ and $\zeta,b\in C^1([0,T];H^{t_0+2}(\R^2))$
	be such that (\ref{ell2}) is satisfied (uniformly with
	respect to $t$), and let $\sigma$ be given
	by (\ref{ell2bis}). Then, for all 
	$u\in C^1([0,T];H^{1/2}(\R^2))$ and $t\in [0,T]$,
	$$
	\big\vert \big([\dt,\frac{1}{\mu\nu}\GG] u(t),u(t)\big)
	\big\vert
	\leq
	M[\sigma(t)]\, \Vert \nagm \dt \sigma\Vert_{\infty,T}
	\vert \Pig u(t) \vert_{2}^2,
	$$
	where $M[\sigma(t)]$ is as in (\ref{eqM}) while 
	$\Pig$ is defined in (\ref{Pig}).
\end{prop}

\begin{proof}
First remark that $$(u,[\dt,\GG] u)=\dt (u,\GG u)-2(u,\GG \dt u),$$ so that
using Green's identity, one gets
\begin{eqnarray*}
	\lefteqn{(u,[\dt,\GG]u)=
	\dt \int_\cS \matr \nagm u^\flat\cdot \nagm u^\flat}\\
	& &-2 \int_\cS \matr \nagm (\dt u)^\flat\cdot\nagm u^\flat\\
	&=&\int_\cS (\dt Q[\sigma])\nagm u^\flat\cdot\nagm u^\flat \\
	& & -2\int_\cS \matr\nagm\big((\dt u)^\flat-\dt u^\flat\big)
	\cdot \nagm u^\flat.
\end{eqnarray*}
It follows directly that
\begin{eqnarray*}
	(u,[\dt,\GG]u)&\lesssim&
	\Vert \dt Q[\sigma]\Vert_\infty\Vert \nagm u^\flat\Vert_2^2\\
	&+&(1+\Vert Q[\sigma]\Vert_\infty) 
	\Vert \nagm\big((\dt u)^\flat-\dt u^\flat\big)
	\Vert_2\big\Vert \nagm u^\flat\Vert_2.
\end{eqnarray*}
Proceeding exactly as in the proof of Lemma \ref{lemmg}, one gets
$$\Vert \nagm \big((\dt u)^\flat-\dt u^\flat\big)\Vert_2
\lesssim \Vert \dt Q[\sigma]\Vert_\infty \Vert \nagm u^\flat\Vert_2,$$ and
the result follows therefore from Corollary \ref{coro1} and 
Proposition \ref{propQ}. \qed
\end{proof}

\subsection{Other properties}\label{sectother}

Propositions \ref{propDNoth1} and \ref{propDNoth2} allow one to control
$(u, \GG v)$ in general. However, it is sometimes necessary to have
more precise estimates, when  $u$ and $v$ have some special structure that can 
be exploited.
\begin{prop}
	\label{propDNoth3}
	Let $t_0>1$, $s\geq 0$ and $\zeta,b\in H^{t_0+2}\cap H^{s+2}(\R^2)$
	be such that (\ref{ell2}) is satisfied, and let $\sigma$ be given
	by (\ref{ell2bis}). \\
	{\bf i.} For all $\bv\in H^{s+1}\cap H^{t_0+2}(\R^2)^2$ 
	and $u\in H^{s+1/2}(\R^2)$,
	one has
	\begin{eqnarray*}
	\lefteqn{\big([\Lambda^s,\bv]\cdot \nag u,\frac{1}{\mu\nu}\GG 
	([\Lambda^s,\bv]\cdot \nag u)\big)^{1/2}
	\leq
	M[\sigma]}\\
	&\times &\big(
		\vert \bv\vert_{H^{t_0+2}}\big\vert \Pig u\big\vert_{H^{s}}
	+\big\langle
	\vert \bv\vert_{H^{s+1}}\big\vert \Pig u\big\vert_{H^{t_0+1}}
	\big\rangle_{s> t_0+1}\big).
	\end{eqnarray*}
	{\bf ii.}
	For all $\bv\in H^{t_0+1}(\R^2)^2$ and $u\in H^{1/2}(\R^2)$, one has
$$
	\big((\bv\cdot\nag u),\frac{1}{\mu\nu}\GG u\big)
	\leq M[\sigma] \, \vert \bv\vert_{W^{1,\infty}}
	\vert \Pig u \vert_2^2.
$$
\end{prop}
\begin{proof}
In order to prove the first point of the proposition, define
$U=[\Lambda^s,\bv]\cdot \nag u$ and recall that we saw in
(\ref{DNoth2}) that
$$
	(U,\GG U)=\big\Vert (1+Q[\sigma])^{1/2}\nagm U^\flat\big\Vert_2^2.
$$
Using Notation \ref{notaDNoth1}, we define $U^\natural
=[\Lambda^s,\bv]\cdot\nag u^\dag$;
as in Remark \ref{remaprov} (with $U^\natural$ instead of 
$U^\dag$), one deduces
$$
	(U,\GG U)\leq \big\Vert (1+Q[\sigma])^{1/2}\nagm U^\natural
	\big\Vert_2^2.
$$
Since $\gamma\leq 1$, one has $\Vert \nag U^\natural\Vert_2\lesssim \Vert \Lambda U^\natural\Vert_2$ and one gets with 
Proposition \ref{propprel1},
$$
	\Vert \nag U^\natural\Vert_2 \lesssim
	\vert \bv\vert_{H^{t_0+1}}\Vert \Lambda^s\nag u^\dag\Vert_{2}
	+\big\langle
	\vert \bv\vert_{H^{s+1}}\Vert \Lambda^{t_0}\nag u^\dag\Vert_{2}
	\big\rangle_{s> t_0};
$$
similarly, since $\dz U^\natural=[\Lambda^s,\bv]\cdot\nag \dz u^\dag$,
one gets
$$
	\Vert \dz U^\natural\Vert_2 \lesssim
	\vert \bv\vert_{H^{t_0+1}}\Vert \Lambda^s\dz u^\dag\Vert_{2}
	+\big\langle
	\vert \bv\vert_{H^{s}}\Vert \Lambda^{t_0+1}\dz u^\dag\Vert_{2}
	\big\rangle_{s> t_0+1},
$$
so that, finally,
$$
	\Vert \nagm U^\natural\Vert_2 \lesssim
	\vert \bv\vert_{H^{t_0+2}}\Vert \Lambda^s\nagm u^\dag\Vert_{2}
	+\big\langle
	\vert \bv\vert_{H^{s+1}}\Vert \Lambda^{t_0+1}\nagm u^\dag\Vert_{2}
	\big\rangle_{s> t_0+1}.
$$
It follows therefore from Proposition \ref{propsharp} that
\begin{eqnarray}
	\nonumber
	\Vert \nagm U^\natural\Vert_2 &\lesssim&
	\sqrt{\mu}\vert \bv\vert_{H^{t_0+2}}\big\vert \frac{\Dg}{(1+\sqrt{\mu}\Dg)^{1/2}}u\big\vert_{H^{s}}\\
	\label{DNoth3}
	& &
	+\big\langle\sqrt{\mu}
	\vert \bv\vert_{H^{s+1}}\big\vert \frac{\Dg}{(1+\sqrt{\mu}\Dg)^{1/2}}u\big\vert_{H^{t_0+1}}
	\big\rangle_{s> t_0+1},\quad
\end{eqnarray}
and the result follows.\\
To establish the second point of the proposition, first remark that, 
owing to Green's identity,
\begin{equation}
	\label{DNoth4}
	\big((\bv\cdot\nag u),\GG u \big)
	=\int_{\cS}\matr \nagm u^\flat\cdot\nagm(\bv\cdot\nag u^\flat),
\end{equation}
so that, 
\begin{eqnarray}
	\label{DNoth5}
	& &\big((\bv\cdot\nag u),\GG u \big) =
	\int_{\cS} \matr\nagm u^\flat\cdot [\nagm,\bv\cdot\nag]u^\flat 
	\nonumber \\
	& &+\int_{\cS}\nagm u^\flat\cdot [Q[\sigma],(\bv\cdot\nag)]\nagm 
	u^\flat 
        \nonumber  \\
        & &+\int_\cS \nagm u^\flat\cdot (\bv\cdot\nag)\matr\nagm u^\flat.
\end{eqnarray}
Integrating by parts, one finds
\begin{eqnarray}
	\nonumber
	\lefteqn{\int_\cS \nagm u^\flat\cdot (\bv\cdot\nag)\matr\nagm u^\flat}\\
	&=&-\int_\cS \big((\divg \bv)+\bv\cdot\nag\big)\nagm u^\flat\cdot \matr\nagm u^\flat 
	 \label{DNoth5a}\\
	\nonumber
	&=&-\int_\cS (\divg \bv)\nagm u^\flat\cdot \matr \nagm u^\flat \\
        & & \nonumber  
	-\int_\cS [\bv\cdot\nag,\nagm] u^\flat\cdot \matr \nagm u^\flat\\
	\label{DNoth6}
	& &
	-\int_\cS \nagm (\bv\cdot\nag u^\flat)\cdot \matr \nagm u^\flat.
\end{eqnarray}
From (\ref{DNoth5}), (\ref{DNoth5a}) and (\ref{DNoth6}), one gets therefore
\begin{eqnarray*}
	\big((\bv\cdot\nag u),\GG u\big)&=&
	\int_{\cS} \matr\nagm u^\flat\cdot [\nagm,\bv\cdot\nag]u^\flat\\
	&+&\frac{1}{2}\int_{\cS}\nagm u^\flat\cdot [Q[\sigma],(\bv\cdot\nag)]\nagm u^\flat\\
	&-&\frac{1}{2}\int_\cS (\divg \bv) \nagm u^\flat\cdot \matr\nagm u^\flat.
\end{eqnarray*}
Remarking that 
$
	[\nagm,\bv\cdot\nag]=\left(\begin{array}{c}
	\nag \bv_1\sqrt{\mu}\dx+	\nag \bv_2\gamma\sqrt{\mu}\dy\\
	0			   \end{array}\right),
$
one deduces easily that
$$
	\big((\bv\cdot\nag u),\GG u\big)
	\lesssim \vert \bv\vert_{W^{1,\infty}}(1+\Vert Q[\sigma]\Vert_{W^{1,\infty}})
	\Vert \nagm u^\flat\Vert_2^2
$$
and the results follows from Corollary \ref{coro1}. \qed
\end{proof}

We finally state the following theorem, which gives an explicit
formula for the shape derivative of the Dirichlet-Neumann operator. This
theorem is a particular case of Theorem 3.20 of \cite{LannesJAMS}.
\begin{theo}
	\label{theoshape}
	Let $t_0>1$, $s\geq t_0$ and $\uz,b\in H^{s+3/2}(\R^2)$
	be such that (\ref{ell2}) is satisfied.  For all
	$\up\in H^{s+3/2}(\R^2)$,  the mapping
	$$
	\zeta\mapsto \cG[\eps\zeta]\up \in H^{s+1/2}(\R^2)
	$$
	is well defined and differentiable in a neighborhood of
	$\uz$ in $H^{s+3/2}(\R^2)$, and
	$$
	\forall h\in H^{s+3/2}(\R^2),\qquad
	d_{\uz}\cG[\eps\cdot]\up\cdot h=
	-\eps \cG[\eps\uz](h\uZ)
	-\eps\mu\nag\cdot(h\bv),
	$$
	with $\uZ:=\cZ[\eps\uz]\up$ and $\bv:=\nag\up-\eps\uZ\nag\uz$, and
	where
	$$
	\cZ[\eps\uz]:=\frac{1}{1+\eps^2\mu\vert \nag\uz\vert^2}
	(\cG[\eps\uz]+\eps\mu\nag\uz\cdot\nag).
	$$
\end{theo}
\begin{rema}
We take this opportunity to correct a harmless misprint in the statement
of Theorem 3.20 of \cite{LannesJAMS}. It should read
$$
	d_{\underline{a}}G(\cdot,b)f\cdot h=-G(\underline{a},b)(h\uZ)
	-\left(\begin{array}{c}\nabla_X\\0\end{array}\right)
	\cdot\Big[hP \left(\begin{array}{c}\bv\\ \uZ\end{array}\right)\Big],
$$
and $\widetilde{P}_{\underline{a}}$ should be replaced by $P$ on the 
right hand side of the equation in the statement of Lemma 3.24.
\end{rema}
\subsection{Asymptotic expansions}\label{sectas}

This subsection is devoted to the asymptotic expansion of the DN
operator $\GG\psi(=\G\psi)$ in terms of one or several of the parameters
$\eps$, $\mu$, $\gamma$  and $\beta$. We consider two cases which cover
all the physical regimes described in the introduction. 

\subsubsection{Expansions in shallow-water ($\mu\ll 1$)}

In shallow water, that is when $\mu\ll 1$, the Laplace equation (\ref{ell1})
--or its straightened version (\ref{ell3})-- reduces at first order to
the ODE $\dz^2\Phi=0$. 
This fact can be exploited to find an approximate solution
$\Phi_{app}$ of the Laplace equation by a standard BKW expansion. This method 
has been used in the long-waves regime in \cite{BCL,LannesSaut,Chazel} 
(see also \cite{NichollsReitich}) where 
the corresponding expansions of the DN operator can be found. We prove
here that it can be used uniformly with respect to $\eps$ and $\beta$,
which allows one to consider at once  
the shallow-water/Green-Naghdi and Serre scalings.
The difference
between both regimes is that $\eps=\beta=1$ in the former (large amplitude for 
the surface and bottom variations), while $\eps=\beta=\sqrt{\mu}$ in the
latter (medium amplitude variations for the surface and bottom variations).\\
Let us first define the first order linear
operator ${\mathcal T}[h,b]$ as
\begin{equation}
	\label{Tb}
	{\mathcal T}[h,b] V\!\!:=\!\!
	-\frac{1}{3}\nabla(h^3\nabla\cdot V)
	\!+\!\frac{1}{2}\big[
	\nabla(h^2\nabla b \cdot V)\!-\!h^2\nabla b \nabla\cdot V\big]
	\!+\!h\nabla b\nabla b\cdot V.
\end{equation}
\begin{prop}[Shallow-water and Serre scalings]\label{propshallow}
	Let $\gamma=1$, 
	$s\geq t_0>1$, $\nabla \psi\in H^{s+11/2}(\R^2)$, 
	$b\in H^{s+11/2}(\R^2)$ 
	and
	$\zeta\in H^{s+9/2}(\R^2)$ 
	and assume that (\ref{ell2}) is satisfied.\\
	With $h:=1+\eps\zeta-\beta b$, one then has
	\begin{eqnarray*}
	\big\vert \GG\psi-\nabla\cdot\big(-\mu h\nabla\psi\big)
	\big\vert_{H^s}
	&\leq&
	\mu^2 C_0 \\
	\big\vert \GG\psi-\nabla\cdot\big(-\mu h\nabla\psi
	+\mu^2
	{\mathcal T}[h,\beta b] \nabla\psi\big)
	\big\vert_{H^s}
	&\leq&
	\mu^3 C_1,
	\end{eqnarray*}
	with $C_j=C(\frac{1}{h_0},\vert \zeta\vert_{H^{s+5/2+2j}},
	\vert b\vert_{H^{s+7/2+2j}},
	\vert \nag\psi\vert_{H^{s+7/2+2j}})$ ($j=0,1$),
	and uniformly with respect to $\eps,\beta\in [0,1]$.
\end{prop}
\begin{proof}
We 
look for an approximate solution
$\phi_{app}$ to the exact solution $\phi$ of the potential equation
(\ref{ell3}) under the form
$$
	\phi_{app}(X,z)=\psi(X)+\mu\phi_1(X,z).
$$
Plugging this ansatz into (\ref{ell3}), and expanding the result into
powers of $\mu$, one can cancel the leading term by a good choice of
$\phi_1$, namely, 
$$
	\phi_1(X,z)=-h\big( h(\frac{z^2}{2}+z)\Delta\psi
	-z\beta\nabla b\cdot\nabla\psi\big).
$$
One can then check that
$$
	\left\lbrace
	\begin{array}{l}
	\nabla_{X,z}\cdot P[\sigma]\nabla_{X,z}\phi_{app}=\mu^2 R_\mu,
	\qquad
	\mbox{ in }{\mathcal S},\\
	\phi_{app}\,_{\vert_{z=0}}=\psi,\qquad
	\partial_n \phi_{app}\,_{\vert_{z=-1}}=
	\mu^2 r_\mu,
	\end{array}\right.
$$
with $(R_\mu,r_\mu)$ satisfying, uniformly with respect to $\mu\in (0,1)$,
\begin{equation}\label{eqreste}
	\Vert \Lambda^{s+1/2} R_\mu\Vert_2
	+\vert r_\mu\vert_{H^{s+1/2}}\leq 
	C(\vert \zeta\vert_{H^{s+5/2}},\vert b\vert_{H^{s+7/2}},
	\vert \nag\psi\vert_{H^{s+7/2}}).
\end{equation}
Since $\GG \psi-
%\sqrt{1+\vert\nabla\zeta\vert^2}
\partial_n\phi_{app}\,_{\vert_{z=0}} = 
%\sqrt{1+\vert\nabla\zeta\vert^2}
\partial_n(\phi-\phi_{app})_{\vert_{z=0}}$,
the truncation error can be estimated using the trace theorem and
an elliptic estimate on the BVP solved by $\phi-\phi_{app}$; this
is exactly what is done in Theorem 1.6 of \cite{Chazel} for instance, 
which gives here:
$$
 \vert \GG \psi-
 %\sqrt{1  + \eps^2\vert\nabla\zeta\vert^2}
 %\partial_n\phi_{app}\,_\fond\vert_{\!H^s}
 \partial_n\phi_{app}\,_{\vert_{z=0}}\vert_{\!H^s}
  \leq  \mu^2C_s (\Vert \Lambda^{s+1/2} R_\mu\Vert_2
 +\vert r_\mu\vert_{H^{s+1/2}}) ,
$$
with $C_s=C(\vert \zeta\vert_{H^{s+5/2}},\vert b\vert_{H^{s+5/2}})$.
Together with (\ref{eqreste}), this gives the result.\\
In order to prove the second estimate of the proposition, one must
look for a higher order approximate solution of (\ref{ell3}), namely
$\phi_{app}=\psi+\mu \phi_1+\mu^2\phi_2$. The computations can
be performed by any software of symbolic calculus and
the estimates are exactly the same as above; we thus omit this
technical step \qed
\end{proof}

\subsubsection{The case of small amplitude waves ($\eps\ll 1$)}

Expansions of the Dirichlet-Neumann operator for small
amplitude waves has been developed in \cite{CSS1,CSS2}.
This method is very efficient to compute the formal
expansion, but instead of adapting it in the present case
to give uniform estimates 
on the truncation error, we
rather propose a very simple method based on Theorem \ref{theoshape}.
\begin{prop}\label{proptrunc}
	Let $s\geq t_0>1$, $\Pig\psi\in H^{s+1/2}(\R^2)$ and
	$\zeta\in H^{s+3/2}(\R^2)$ 
	be such that (\ref{ell2}) is satisfied for some $h_0>0$.
	Then one has
	\begin{eqnarray*}
	& &\big\vert \GG\psi-\big[\cG[0]\psi-\eps \cG
        [0]\big(\zeta (\cG[0]\psi)\big)
	-\eps\mu\nag\cdot(\zeta \nag\psi)\big]\big\vert_{H^s} \\
	& & \leq (\frac{\eps}{\nu})^2\sqrt{\mu}
	C\big(\frac{1}{h_0},\eps\sqrt{\mu},\vert \zeta\vert_{H^{s+3/2}},
	\vert \Pig
	\psi\vert_{H^{s+1/2}}\big).
	\end{eqnarray*}
\end{prop}
\begin{proof}
A second order Taylor expansion of $\GG\psi$ gives
$$
	\GG\psi=\cG[0]\psi
	+d_0\cG[\eps\cdot]\psi\cdot \zeta
	+\int_0^1 (1-z)d^2_{z\zeta} \cG[\eps\cdot]\psi 
	\cdot(\zeta,\zeta)dz.
$$
Using Theorem \ref{theoshape}, one computes 
$$
	d_0\cG[\eps\cdot]\psi\cdot \zeta=
	-\eps \cG[0]\big(\zeta (\cG[0]\psi)\big)
	-\eps\mu \nag\cdot(\zeta\nag\psi),
$$ 
while for all $z\in [-1,0]$, Proposition \ref{unifcontrol}
controls $d^2_{z\zeta} \cG[\eps\cdot]\psi 
	\cdot(\zeta,\zeta)$ in $H^s$ by the
r.h.s. of the estimate given in the statement. \qed
\end{proof}

\section{Linear analysis}\label{sectLin}

The water-waves equations (\ref{nondimww}) can be written in condensed form as
$$
	\dt U+{\mathcal L}U+
	\frac{\eps}{\nu}{\mathcal A}[U]=0,
$$
with $U=(\zeta,\psi)^T$, $\cA[U]=(\cA_1[U],\cA_2[U])^T$ and where
\begin{equation}
	\label{L}
	\cL:=\left(\begin{array}{cc}
	0 & -\frac{1}{\mu\nu}\cG[0]\cdot\\
	1 & 0
	     \end{array}\right)
\end{equation}
and
\begin{equation}
	\label{A}
	\begin{array}{rcl}
	\dsp\cA_1[U]&=&-\frac{1}{\eps\mu}
	(\cG[\eps\zeta]\psi-\cG[0]\psi),
	\vspace{1mm}\\
	\dsp\cA_2[U]&=&\dsp\frac{1}{2}\vert\nag\psi\vert^2-\frac{(\frac{1}{\sqrt{\mu}}\cG[\eps\zeta]\psi+\eps\sqrt{\mu}\nag\zeta\cdot\nag\psi)^2}{2(1+\eps^2\mu\vert\nag\zeta\vert^2)}.
	\end{array}
\end{equation}

By definition, the linearized operator $\cLu$ around some reference
state $\uU=(\uz,\up)^T$ is given by
$$
	\cLu= \dt + \cL+\frac{\eps}{\nu}d_{\uU}\cA;
$$
assuming that $\uU$ is such that the assumptions of Theorem \ref{theoshape}
are satisfied, one computes that $\cLu$ is equal to
\begin{equation}
	\label{lin1}
	\dt+
	\left(
	\begin{array}{cc}
	\frac{\eps}{\mu\nu} \cG[\eps\uz](\uZ\cdot)
	+\frac{\eps}{\nu}\nag\cdot(\cdot \bv) &
	-\frac{1}{\mu\nu}\cG[\eps\uz]\cdot\\
	\frac{\eps^2}{\mu\nu}\uZ \cG[\eps\uz](\uZ\cdot)+(1+\frac{\eps^2}{\nu}\uZ \nag\cdot\bv) & \quad\frac{\eps}{\nu}\bv\cdot\nag\cdot-\frac{\eps}{\mu\nu}\uZ \cG[\eps\uz]\cdot
	\end{array}\right),
\end{equation}
where $\bv$ and $\uZ$ are as in the statement of Theorem \ref{theoshape}.

\bigbreak

This section is devoted to the proof of energy estimates for the 
associated initial value problem,
\begin{equation}
	\label{lin2}
	\left\lbrace
	\begin{array}{l}
	\cLu U=\frac{\eps}{\nu} G\\
	U_\init=U^0.
	\end{array}\right.
\end{equation}
Defining
\begin{equation}
	\label{lin3}
	\mfa=1+\frac{\eps}{\nu} \mfb,
	\quad
	\mbox{ and }
	\quad
	\mfb=\eps \bv\cdot\nag \uZ +\nu\dt\uZ,
\end{equation}
we first introduce the notion of \emph{admissible} reference state:
\begin{defi} \label{defiadmissible}
	Let $t_0>1$, $T>0$ and $b\in H^{t_0+2}(\R^2)$.
	We say that $\uU=(\uz,\up)$ is \emph{admissible on} 
	$[0,\frac{\nu T}{\eps}]$ if
	\begin{itemize}
	\item The surface and bottom parameterizations $\uz$ and $b$
	satisfy (\ref{ell2}) for some $h_0>0$, uniformly on 
	$[0,\frac{\nu T}{\eps}]$;
	\item There exists $c_0>0$ such that $\mfa\geq c_0$, uniformly
	on $[0,\frac{\nu T}{\eps}]$.
	\end{itemize}
\end{defi}
We also need to define some functional spaces and notations
linked to the energy (\ref{introNRJ}) mentioned in the introduction.
\begin{defi} \label{defispaces}For all $s\in\R$ and $T>0$,\\
	{\bf i.} We denote by $X^s$ the vector 
	space $H^s(\R^2)\times H^{s+1/2}(\R^2)$ endowed with the
	norm
	$$
	\forall U=(\zeta,\psi)^T\in X^s,
	\qquad
	\vert U\vert_{X^s}:=\vert \zeta\vert_{H^s}+
	\frac{\eps}{\nu}\vert \psi\vert_{H^s}+
	\vert \Pig \psi\vert_{H^s},
	$$
	while $X^s_T$ stands for $C([0,\frac{\nu T}{\eps}];X^s)$
	endowed with its canonical norm.
	{\bf ii.} We define the space $\widetilde{X}^s$ as
	$$
	\widetilde{X}^s:=\{U=(\zeta,\psi)^T,\zeta\in H^s(\R^2),
	\nabla \psi\in H^{s-1/2}(\R^2)^2\},
	$$
	and endow it with the semi-norm
	$\vert U \vert_{\widetilde{X}^s}
	:=\vert \zeta\vert_{H^s}+\vert\Pig\psi\vert_{H^s}$.\\
	{\bf iii.} We define the semi-normed space 
	$(Y^s_T,\vert \cdot\vert_{Y^s_T})$ as
	$$
	Y^s_T\!:=\! \bigcap_{k=0}^2 C^k([0,\frac{\nu T}{\eps}];\widetilde{X}^{s-\frac{3}{2}k})
	\quad \mbox{and}\quad
	\vert U\vert_{Y^s_T}\!=\!\sum_{k=0}^2 
	\sup_{[0,\frac{\nu T}{\eps}]}
	\vert \dt^k U\vert_{\widetilde{X}^{s-\frac{3}{2}k}}.
	$$
	{\bf iv.} For all  $(G,U^0)\in X^s_T\times X^s$,
	we define
	$$
	{\mathcal I}^s(t,U^0,G):=
	\vert U^0\vert_{X^s}+
	\frac{\eps}{\nu}\int_0^t \sup_{0\leq t''\leq t'}\vert G(t'')\vert_{X^s}		dt'.
	$$
\end{defi}
We can now state the energy estimate associated to (\ref{lin2}), and
whose proof is given in the next two subsections. 
\begin{prop}
	\label{propmain}
	Let $s\geq t_0>1$, $T>0$,  $b\in H^{s+9/2}(\R^2)$,
	and
	$\uU=(\uz,\up)\in Y^{s+9/2}_{T}$
	be admissible on 
	$[0,\frac{\nu T}{\eps}]$
	for some $h_0>0$ and $c_0>0$. Let also
	$(G,U^0)\in X^{s+2}_{T}\times X^{s+2}$.\\
	There exists
	a unique solution $U\in X^s_{T}$ 
	to (\ref{lin2}); moreover, for all $0\leq t\leq \frac{\nu T}{\eps}$,
	one has
	$$
	\vert U(t)\vert_{X^s}\leq \uC \big({\mathcal I}^{s+2}(t,U^0,G)+
	\vert \uU\vert_{Y^{s+9/2}_{T}}
	{\mathcal I}^{t_0+2}(t,U^0,G)  
	\big),
	$$
	where 
	$
	\uC=
	C\big(T,\frac{1}{h_0},\frac{1}{c_0},
	\frac{\eps}{\nu},\frac{\beta}{\eps},
	\vert b\vert_{H^{s+9/2}},\vert\uU\vert_{Y^{t_0+9/2}_{T}}\big)
	$.
\end{prop}

\subsection{Energy estimates for the trigonalized linearized operator}

As shown in \cite{LannesJAMS}, the operator $\cLu$ is non-strictly hyperbolic,
in the sense that its principal symbol has a double purely imaginary 
eigenvalue, with a nontrivial Jordan block. It was shown in Prop. 4.2 of 
\cite{LannesJAMS} that a simple change of basis can be used to put
the principal symbol of $\cLu$ under a canonical trigonal form. This
result is generalized to the present case. More precisely, with $\mfa$ as defined in (\ref{lin3}) and
 defining
the operator $\cMu=\dt+\Mu$ with
\begin{equation}
	\label{lin4}
	\Mu=
	\left(\begin{array}{cc}
	\frac{\eps}{\nu}\nag\cdot(\cdot \bv) &
	-\frac{1}{\mu\nu}\GGu\cdot\\
	\mfa  & \frac{\eps}{\nu}\bv\cdot\nag\cdot
	\end{array}\right),
\end{equation}
one reduces the study of (\ref{lin2}) 
to the study of the initial value problem
\begin{equation}
	\label{lin5}
	\left\lbrace
	\begin{array}{l}
	\cMu V=\frac{\eps}{\nu} H\\
	V_\init=V^0,
	\end{array}\right.
\end{equation}
as shown in the following proposition (whose proof relies on simple
computations and is omitted).
\begin{prop}
	The following two assertions are equivalent:
	\begin{itemize}
	\item The pair $U=(\zeta,\psi)^T$ solves (\ref{lin2});
	\item The pair $V=(\zeta,\psi-\eps\uZ\zeta)^T$ solves
	(\ref{lin5}), with $H=(G_1,G_2-\eps\uZ G_1)^T$ and
	$V^0=(\zeta^0,\psi^0-\eps\uZ_\init\zeta^0)^T$.
	\end{itemize}
\end{prop}

In view of this proposition, it is a key step to understand 
(\ref{lin5}), and the rest of this subsection is thus 
devoted to the proof of energy estimates for this initial value problem.

First remark that a symmetrizer for $\cMu$ is given by
\begin{equation}
	\label{lin6}
	S=\left(
	\begin{array}{cc}
	\mfa & 0\\
	0 & \frac{\eps^2}{\nu^2}+\frac{1}{\mu\nu}\GGu\cdot
	\end{array}\right),
\end{equation}
so that (provided that $\mfa$ is nonnegative), a natural energy for the IVP (\ref{lin5}) is given by
\begin{eqnarray}
	\nonumber
	E^s(V)^2&=&(\Lambda^s V,S\Lambda^s V)\\
	\label{lin7}
	&=& \vert \sqrt{\mfa}\Lambda^sV_1\vert_{2}^2
	+\frac{\eps^2}{\nu^2}\vert  V_2\vert_{H^s}^2
	+(\Lambda^s V_2,\frac{1}{\mu\nu}\GGu \Lambda^s V_2).
\end{eqnarray}
\begin{rema}
The introduction of the term $\eps^2/\nu^2$ in (\ref{lin6}) --and thus of $\eps^2/\nu^2\vert V_2\vert_{H^s}^2$ in (\ref{lin7})-- is not necessary to the energy estimate below. But this constant term plays a crucial role in the
iterative scheme used to solve the nonlinear problem because it controls
the low frequencies. It also turns out that the order $O(\eps^2/\nu^2)$
of this constant term is the only one which allows uniform estimates.
\end{rema}

We can now give the energy estimate associated to (\ref{lin5});
in the statement below, we use the notation 
$$
	{I}^s(t,V^0,H):=
	E^s(V^0)+
	\frac{\eps}{\nu}\int_0^t \sup_{0\leq t''\leq t'} E^s(H(t''))dt',
	$$
while $s\vee t_0:=\max\{s,t_0\}$ and  $\uC$ is as defined
in Proposition \ref{propmain}.
\begin{prop}
	\label{proptrig}
	Let $s\geq 0$, $t_0>1$, $T>0$,  $b\in H^{s\vee t_0+9/2}(\R^2)$,
	and
	$\uU=(\uz,\up)\in Y^{s\vee t_0+9/2}_{T}$
	be admissible on 
	$[0,\frac{\nu T}{\eps}]$
	for some $h_0>0$ and $c_0>0$.\\
	Then, for all $(H,V^0)\in X^{s}_{T}\times X^s$,
	there exists
	a unique solution $V\in X^{s}_{T}$ 
	to (\ref{lin5}) and
	for all $0\leq t\leq \frac{\nu T}{\eps}$,
	$$
	E^s(V(t))\leq
	\underline{C} \big(I^s(t,V^0,H)
	+\big\langle \vert \uU\vert_{Y^{s+7/2}_{T}}
	I^{t_0+1}(t,V^0,H)
	\big\rangle_{s> t_0+1}\big).
	$$
\end{prop}
\begin{proof}
\emph{Throughout this proof, $\uC_0$ denotes a nondecreasing function of 
$\frac{1}{c_0}$, $\frac{\eps}{\nu}$, $M[\usi]$, $\vert \bv\vert_{H^{t_0+2}}$, $\vert \mfb\vert_{H^{t_0+2}}$, and $\vert \dt \mfb\vert_{\infty}$ 
which may vary from one line to another, and $\usi$ is given by 
(\ref{ell2bis}) with $\zeta=\uz$}.\\
Existence of a solution to the IVP (\ref{lin5}) is achieved by classical
means, and we thus focus our attention on the proof of the energy estimate. 
For any given $\kappa\in\R$, we compute
\begin{eqnarray}
	\nonumber
	& & e^{\frac{\eps\kappa}{\nu} t}\frac{d}{dt}(e^{-\frac{\eps\kappa}{\nu} t}E^s(V)^2)
        =-\frac{\eps\kappa}{\nu} E^s(V)^2
	+2\frac{\eps}{\nu}(\Lambda^s H,S\Lambda^s V)\\
	\label{lin8}
	& &-2(\Lambda^s \Mu V,S\Lambda^s V)
	+(\Lambda^s V, [\dt,S]\Lambda^s V).
\end{eqnarray}
We now turn to bound from above the different components of the r.h.s. of
(\ref{lin8}).

\noindent
$\bullet$ Estimate of $(\Lambda^s H,S\Lambda^s V)$. We can rewrite this 
term as 
$$
	(\sqrt{\mfa}\Lambda^s H_1,\sqrt{\mfa}\Lambda^s V_1)
	+(\frac{\eps}{\nu} \Lambda^s H_2,\frac{\eps}{\nu}\Lambda^s V_2)
	+(\Lambda^sH_2,\frac{1}{\mu\nu}\GGu \Lambda^sV_2),
$$
so that Cauchy-Schwartz inequality and Proposition \ref{propDNoth1} yield
\begin{equation}
	\label{lin9}
	(\Lambda^s H,S\Lambda^s V)\leq E^s(H)E^s(V).
\end{equation}

\noindent
$\bullet$ Estimate of $(\Lambda^s \Mu V,S\Lambda^s V)$. One computes
\begin{eqnarray*}
	& & (\Lambda^s \Mu V,S\Lambda^s V)=
	\big(\Lambda^s (\frac{\eps}{\nu} \divg(\bv V_1)-\frac{1}{\mu\nu}\GGu V_2),\mfa \Lambda^s V_1\big)\\
	& &+\big(\Lambda^s(\mfa V_1+\frac{\eps}{\nu}\bv\cdot \nag V_2), 
           (\frac{\eps^2}{\nu^2}+\frac{1}{\mu\nu}\GGu)\Lambda^s V_2\big),
\end{eqnarray*}
so that one can write
$$
	(\Lambda^s \Mu V,S\Lambda^s V)=A_1+A_2+A_3+A_4+A_5,
$$
with
\begin{eqnarray*}
	A_1&=&\frac{\eps}{\nu}\big(\Lambda^s \divg(\bv V_1),\mfa \Lambda^s V_1\big),\\
	A_2&=&\frac{\eps}{\nu}\big(\Lambda^s (\bv\cdot\nag V_2),\frac{\eps^2}{\nu^2} \Lambda^s V_2\big),\\
	A_3&=&\frac{\eps}{\nu}\big(\Lambda^s(\bv\cdot \nag V_2),\frac{1}{\mu\nu}\GGu\Lambda^s V_2 \big),\\
	A_4&=&\big(\Lambda^s (\mfa V_1),\frac{1}{\mu\nu}\GGu\Lambda^s V_2\big)
	-\big(\mfa\Lambda^s V_1,\Lambda^s (\frac{1}{\mu\nu}\GGu V_2)\big),\\
	A_5&=&\big(\Lambda^s (\mfa V_1),\frac{\eps^2}{\nu^2}\Lambda^s V_2).
\end{eqnarray*}
We now turn to prove the following estimates:
\begin{eqnarray}
	\nonumber
	\!\!\!\!\!\!\!\!\!\!\!\!\!\!
        & & A_j\leq \frac{\eps}{\nu} \uC_0E^s(V)
	\Big(
	\big(1+\frac{\nu}{\eps}\Vert\nagm\underline{\sigma}\Vert_{L^\infty H^{t_0+2}}\big)
	E^s(V)\\
	\label{estgen}
	\!\!\!\!\!\!\!\!\!\!\!\!\!\!
        & &+ \big\langle \big(\vert \bv\vert_{H^{s+1}}
	+\vert \mfb\vert_{H^{s+1}}
	+\frac{\nu}{\eps}\Vert\nagm\underline{\sigma}\Vert_{L^\infty H^{s+1}}\big)
	E^{t_0+1}(V)\big\rangle_{s> t_0+1}
	\Big),
\end{eqnarray}
for $j=1,\dots,5$.
\begin{itemize}
\item Control of $A_1$ and $A_2$. Integrating by parts, one obtains
\begin{eqnarray*}
	\frac{\nu}{\eps}A_1&=&\big(
	[\Lambda^s,\divg(\bv \cdot)]V_1,\mfa \Lambda^s V_1
	\big)
	-\frac{1}{2}\big(\Lambda^s V_1,(\bv\cdot\nag\mfa)\Lambda^s V_1\big)\\
	& &+\frac{1}{2}\big(
	\mfa\Lambda^s V_1,(\divg \bv)\Lambda^s V_1\big)
\end{eqnarray*}
and
$$
	\frac{\nu}{\eps}A_2=\big(
	[\Lambda^s,\bv]\nag V_2,\frac{\eps^2}{\nu^2}\Lambda^s V_2\big)-\frac{\eps^2}{2\nu^2}
	\big(\Lambda^s V_2,(\divg \bv)\Lambda^s V_2\big).
$$
Recalling that $\mfa=1+\frac{\eps}{\nu}\mfb$, one can then deduce easily (with the help of Proposition \ref{propprel1} and Corollary \ref{propprel2} 
to control the commutators in the above expressions) that (\ref{estgen})
holds for $j=1,2$.

\item Control of $A_3$. First write $A_3=A_{31}+A_{32}$ with
\begin{eqnarray*}
	A_{31}
	&=&
	\frac{\eps}{\nu}\big(\frac{1}{\mu\nu}\GGu \Lambda^s V_2,[\Lambda^s,\bv]\cdot\nag V_2\big)\\
	A_{32}&=&
	\frac{\eps}{\nu}
	\big(\bv\cdot\nag \Lambda^s V_2,\frac{1}{\mu\nu}\GGu\Lambda^sV_2\big).
\end{eqnarray*}
Thanks to Proposition \ref{propDNoth1}, one gets
$$
	A_{31}\leq \frac{\eps}{\nu}
	\big(\frac{1}{\mu\nu}\GGu [\Lambda^s,\bv]\cdot\nag V_2
	,[\Lambda^s,\bv]\cdot\nag V_2\big)^{1/2}
	E^s(V),
$$
and Propositions \ref{propDNoth3}{\bf.i} and \ref{propDNoth2} can then 
be used to
show that $A_{31}$ is bounded from above by the r.h.s. of (\ref{estgen}).
This is also the case of $A_{32}$, as a direct consequence of 
 Propositions  \ref{propDNoth3}{\bf .ii} and \ref{propDNoth2}. It follows
that  (\ref{estgen}) holds for $j=3$.

\item Control of $A_4$. One computes, remarking that $[\Lambda^s,\mfa]=\frac{\eps}{\nu}[\Lambda^s,\mfb]$, 
\begin{eqnarray*}
	A_4&=&\big(
	\mfa \Lambda^s V_1,[\frac{1}{\mu\nu}\GGu,\Lambda^s]V_2
	\big)
	+\frac{\eps}{\nu}
	\big(
	[\Lambda^s,\mfb]V_1,\frac{1}{\mu\nu}\GGu \Lambda^s V_2	
	\big)\\
	&:=&A_{41}+A_{42}.
\end{eqnarray*}
Using successively Cauchy-Schwartz inequality, Proposition \ref{propestcom}, 
and Proposition \ref{propDNoth2}, one obtains directly that $A_{41}$
is bounded from above by the r.h.s. of (\ref{estgen}).
In order to control $A_{42}$, first remark that using Propositions
\ref{propDNoth1} and \ref{propDNoth2}, one gets
\begin{eqnarray*}
	\frac{\nu}{\eps}A_{42}&\leq&\big([\Lambda^s,\mfb]V_1,\frac{1}{\mu\nu}\GGu[\Lambda^s,\mfb]V_1\big)^{1/2}E^s(V)\\
	&\leq&M[\underline{\sigma}]
	\big\vert 
	\Pig
	[\Lambda^s,\mfb]V_1\big\vert_2 E^s(V).
\end{eqnarray*}
Recalling that $\nu=\frac{1}{1+\sqrt{\mu}}$ one can check that
for all $\xi\in \R^2$, $\frac{\nu^{-1/2}\xig}{(1+\sqrt{\mu}\xig)^{1/2}}
\lesssim \jap$, \emph{uniformly with respect to $\mu$ and $\gamma$}, so that
one deduces
$$
	\frac{\nu}{\eps}A_{42}
	\leq M[\underline{\sigma}]
	 \big\vert [\Lambda^s,\mfb]V_1\big\vert_{H^{1}} E^s(V).
$$
Remarking that owing to Proposition \ref{propprel1}, one has
\begin{eqnarray*}
	\big\vert [\Lambda^s,\mfb]V_1\vert_{H^{1}}&\leq&
	 \vert \mfb\vert_{H^{t_0+1}}\vert V_1\vert_{H^{s}}
	+\big\langle \vert \mfb\vert_{H^{s+1}}
	\vert V_1\vert_{H^{t_0}}\big\rangle_{s> t_0}\\
	&\leq&\frac{1}{\sqrt{c_0}}\big(
	 \vert \mfb\vert_{H^{t_0+2}}E^s(V)
	+\big\langle \vert \mfb\vert_{H^{s+1}}
	E^{t_0+1}(V)\big\rangle_{s> t_0+1}\big),
\end{eqnarray*}
and $A_{42}$ is thus bounded from above by the r.h.s. of (\ref{estgen}). This
shows that
(\ref{estgen}) also holds for $j=4$.
\item Control of $A_5$. First remark that 
$$
	A_5=(\Lambda^s V_1,\frac{\eps^2}{\nu^2}\Lambda^s V_2)
	+\frac{\eps}{\nu}(\Lambda^s (\mfb V_1),
	\frac{\eps^2}{\nu^2}\Lambda^s V_2),
$$
so that Cauchy-Schwartz inequality and the tame
product estimate (\ref{tame}) yield
\begin{eqnarray*}
	\!\!\! & & A_5 \leq  
        \frac{\eps}{\nu} \big((1+\vert \frac{\eps}{\nu}\mfb\vert_{H^{t_0}})\vert V_1\vert_{H^s}
	+ \big\langle \vert \frac{\eps}{\nu}\mfb \vert_{H^s}\vert V_1\vert_{H^{t_0}}
        \big\rangle_{s> t_0}\big)\frac{\eps}{\nu} \vert V_2\vert_{H^s}\\
	\!\!\! & & \leq \frac{\eps}{\nu}\frac{1}{\sqrt{c_0}} 
	\big((1+\vert \frac{\eps}{\nu}\mfb\vert_{H^{t_0+1}})
	E^s(V)
	+ \big\langle \vert \frac{\eps}{\nu}\mfb \vert_{H^s}E^{t_0+1}(V)
	\big\rangle_{s> t_0+1}\big)E^s(V),
\end{eqnarray*}
and (\ref{estgen}) thus holds for $j=5$.
\end{itemize}
From (\ref{estgen}), we obtain directly
\begin{eqnarray}
	\nonumber
	\lefteqn{(\Lambda^s \Mu V,S\Lambda^s V)\leq
	\frac{\eps}{\nu} \uC_0E^s(V)
	\big(
	\big(1+\frac{\nu}{\eps}\Vert\nagm\underline{\sigma}
	\Vert_{L^\infty H^{t_0+2}}\big)
	E^s(V)}\\
	\label{lin10}
	&+&\big\langle \big(\vert \bv\vert_{H^{s+1}}
	+\vert \mfb\vert_{H^{s+1}}
	+\frac{\nu}{\eps}\Vert\nagm\underline{\sigma}\Vert_{L^\infty H^{s+1}}\big)
	E^{t_0+1}(V)\big\rangle_{s> t_0+1}
	\big).\nonumber  \\
\end{eqnarray}

\noindent
$\bullet$ Estimate of $(\Lambda^s V, [\dt,S]\Lambda^s V)$. One has
$$
	(\Lambda^s V, [\dt,S]\Lambda^s V)=
	\frac{\eps}{\nu} (\Lambda^s V_1,\dt \mfb \Lambda^s V_1)+(\Lambda^s V_2,[\dt,\frac{1}{\mu\nu}\GGu]\Lambda^s V_2),
$$
so that, using Proposition \ref{time} to control the second component of 
the r.h.s., one gets easily
\begin{equation}
	\label{lin11}
	(\Lambda^s V, [\dt,S]\Lambda^s V)\leq 	
	\frac{\eps}{\nu} \uC_0(1+\frac{\nu}{\eps}\Vert \nagm\dt
	\underline{\sigma}\Vert_{\infty})
	E^s(V)^2.
\end{equation}

According to (\ref{lin8}), (\ref{lin9}), (\ref{lin10}) and (\ref{lin11}), 
we have
\begin{equation}
	\label{lin12}
	e^{\frac{\eps\kappa}{\nu} t}\frac{d}{dt}(e^{-\frac{\eps \kappa}{\nu} t}E^s(V)^2)\leq 
	\frac{\eps}{\nu} E^s(V)
	\big(2E^s(H)+
	\uC_0\big\langle D_s
	E^{t_0+1}(V)\big\rangle_{s> t_0+1}
	\big),
\end{equation}
with $D_s:=\big(\vert \bv\vert_{H^{s+1}}
	+\frac{\nu}{\eps}\Vert\nagm\usi\Vert_{L^\infty H^{s+1}}
	+\vert \mfb\vert_{H^{s+1}}\big)$,
provided that $\kappa$ is large enough, how large depending only on 
$$
	\sup_{t\in[0,\frac{\nu T}{\eps}]}\Big[\uC_0(t)
	\big(1+\frac{\nu}{\eps}\Vert \nagm\usi(t)\Vert_{L^\infty H^{t_0+2}}
	+\frac{\nu}{\eps}\Vert \nagm \dt\usi(t)\Vert_{\infty}
	\big)\Big].
$$
It follows from (\ref{lin12}) that,
\begin{eqnarray}
	\nonumber
	E^s(V(t))&\leq& e^{\frac{\eps \kappa}{\nu} t}E^s(V^0)+
	\frac{\eps}{\nu} \int_0^t e^{\frac{\eps \kappa}{\nu} (t-t')}
	E^s(H(t'))dt'\\
	\label{lin13}
	&+&
	\big\langle 
	\frac{\eps}{\nu} \uC_0(\sup_{[0,\nu T/\eps]}D_s)\int_0^t e^{\frac{\eps \kappa}{\nu} (t-t')}
	E^{t_0+1}(V(t'))dt'
	\big\rangle_{s> t_0+1}; \nonumber \\
\end{eqnarray}
using (\ref{lin13}) with $s=t_0+1$ gives
$$
	E^{t_0+1}(V(t))\leq e^{\frac{\eps\kappa}{\nu} t}E^{t_0+1}(V^0)
	+\frac{\eps}{\nu} te^{\frac{\eps\kappa}{\nu} t }\sup_{0\leq t'\leq t}E^{t_0+1}(H(t')),
$$
and plugging this expression back into (\ref{lin13}) gives therefore
\begin{eqnarray*}
	E^s(V(t))&\leq& 
	\uC_1\big(I^s(t,V^0,H)
	+\big\langle (\sup_{t\in [0,\nu T/\eps]}D_s)
	I^{t_0+1}(t,V^0,H)
	\big\rangle_{s> t_0+1}\big),
	\end{eqnarray*}
	where $\uC_1$ is a nondecreasing function of 
	$T, \frac{1}{c_0}, \frac{1}{h_0},\frac{\eps}{\nu}$
	and of the supremum on the time interval	
	$[0,\frac{\nu T}{\eps}]$ of 
	$\frac{\nu}{\eps}\Vert\nagm \underline{\sigma}
	\Vert_{L^\infty H^{t_0+2}}$,
	$\frac{\nu}{\eps}\Vert \nagm \dt \underline{\sigma}\Vert_{\infty}$,
	$\vert \bv\vert_{H^{t_0+2}}$, $\vert \mfb\vert_{H^{t_0+2}}$ and
	$\vert \dt \mfb\vert_{L^\infty}$.
The proposition follows therefore from the following lemma:
\begin{lemm}
	\label{lemmcontrol}
	With $\uC$ and $\vert\cdot\vert_{Y^s_{T}}$
	as defined in the statement
	of Proposition \ref{propmain} and Definition \ref{defispaces}, 
	one has, 
	$$
		\forall s\geq t_0+1, \qquad
		\sup_{t\in[0,\nu T/\eps]}D_s(t)\leq \uC \vert \uU\vert_{Y^{s+7/2}_{T}}
	\quad\mbox{ and }\quad
	\uC_1\leq \underline{C}.
	$$
\end{lemm}
\begin{proof}
Remark first that, as a consequence of Proposition \ref{unifcontrol}, one has
 for all $r\geq t_0+1$,
\begin{equation}
	\label{marre1}
	\big\vert \frac{1}{\sqrt{\mu}}\GGu \up\big\vert_{H^r}
	\leq C\big(\frac{1}{h_0},\frac{1}{c_0},
	\frac{\eps}{\nu},\frac{\beta}{\eps},
	\vert\uU\vert_{\tX^{t_0+2}},
	\vert b\vert_{H^{r+3/2}}
	\big)
	\vert\uU\vert_{\tX^{r+3/2}};
\end{equation}
since $\xig\leq\frac{\nu^{-1/2}\xig(1+\vert\xi\vert)^{1/2}}{(1+\sqrt{\mu}\xig)^{1/2}}$, uniformly with respect to $\gamma\in (0,1]$ and $\mu>0$, one also has
\begin{equation}
	\label{marre2}
	\vert \nag\up\vert_{H^r}\leq \vert \Pig \up \vert_{H^{r+1/2}}\leq 
	\vert\uU\vert_{\tX^{r+1/2}}.
\end{equation}
It follows from the explicit expression of $\uZ$ given in 
Theorem \ref{theoshape} that $\eps\uZ$ is a smooth function of
$\eps\sqrt{\mu}\leq \eps/\nu$, $\nag\uz$, $\nag\up$ and $\frac{1}{\sqrt{\mu}}\GGu\up$. Moser's type estimates then imply
that for all $r\geq t_0+1$, $\vert \eps\uZ\vert_{H^r}$
-- and hence $\vert \bv\vert_{H^{s+1}}$ -- is bounded from
above by $\uC\vert \uU\vert_{Y^{s+7/2}_{T}}$.
 This is also the
case of the second component of $D_s$, as a direct consequence of 
(\ref{ell2bis}), and because $\nu\sqrt{\mu}\leq 1$.\\
To control the third component of $D_s$, namely, 
$\sup_{[0,\frac{\nu T}{\eps}]}\vert \mfb\vert_{H^{s+1}}$ (with $\mfb$ given by (\ref{lin3})), 
we need to bound
$\vert \nu\dt\uZ\vert_{H^{s+1}}$ from above. Using Theorem
\ref{theoshape} to compute explicitly $\dt\uZ$, one finds
\begin{eqnarray*}
	& &\nu\dt\uZ=
	\frac{\sqrt{\mu}\nu}{1+\eps^2\mu\vert\nag\uz\vert^2}
	\Big(
	\frac{1}{\sqrt{\mu}}\GGu \dt\up+\eps\sqrt{\mu}\nag\uz\cdot\nag\dt\up\\
	& &-\eps\sqrt{\mu}\nag\dt\uz\cdot (\eps\uZ)\nag\uz
	-\eps\sqrt{\mu}\dt\uz \divg\bv
	-\frac{1}{\sqrt{\mu}}\GGu \big(\dt\uz (\eps\uZ)\big)
	\Big),
\end{eqnarray*}
which is a smooth function of $\eps\sqrt{\mu}\leq \eps/\nu$, 
$\nag\uz$, $\dt\uz$, $\nag \dt \uz$, $\nag\dt\up$, $\eps\uZ$, $\frac{1}{\sqrt{\mu}}\GGu\dt \up$ and $\frac{1}{\sqrt{\mu}}\GGu \big(\dt\uz (\eps\uZ)\big)$.
The sought after estimate on $\vert \mfb\vert_{H^{s+1}}$ thus follows
from Moser's type estimates (note that Remark \ref{remopnorm} is used
to control $\frac{1}{\sqrt{\mu}}\GGu \big(\dt\uz (\eps\uZ)\big)$ in terms
of Sobolev norms of $\dt\uz$ and $\eps\uZ$).\\
The estimate on $\uC_1$ is obtained
 exactly in the same way and we omit the proof. \qed
\end{proof}
\qed
\end{proof}

\subsection{Proof of Proposition \ref{propmain}}

Deducing Proposition \ref{propmain} from Proposition \ref{proptrig}
is only a technical step, essentially based on the equivalence
of the norms $E^s$ and $\vert \cdot\vert_{X^s}$ stemming from
Proposition \ref{propDNoth2}.  We only give the main steps of the proof.\\
{\bf Step 1.} Since $U=(V_1,V_2+\eps\uZ V_1)$, one can expand 
$\vert U\vert_{X^s}$
in terms of $V_1$ and $V_2$ and control the different components
using the norm $E^s$ to obtain:
\begin{equation}
	\label{eqmain1}
	\vert U\vert_{X^s}
	\leq \uC \times\big(E^{s+1}(V)+\langle \vert\uU\vert_{\tX^{s+5/2}}
	E^{t_0+1}(V)\rangle_{s> t_0}\big).
\end{equation}

\noindent
{\bf Step 2.} Using Proposition \ref{proptrig} to control $E^{s+1}(V)$
and $E^{t_0+1}(V)$ in terms of $V^0=(U_1^0,U_2^0-\eps\uZ\init U_1^0)$ and
$H=(G_1,G_2-\eps\uZ G_1)$ in (\ref{eqmain1}), one gets
$$
	\vert U(t)\vert_{X^s}\leq \uC \big({\mathcal I}^{s+1}(t,V^0,H)+
	\langle \vert \uU\vert_{Y^{s+9/2}_{T}} 
	({\mathcal I}^{t_0+1}(t,V^0,H))
	\rangle_{s> t_0}\big).
$$

\noindent
{\bf Step 3.} Replacing  $H$ by $(G_1,G_2-\eps\uZ G_1)$ and $V^0$ by
$(U_0^1,U^0_2-\eps\uz\init U^0_1)$ one obtains the following control
on ${\mathcal I}^{r+1}(t,V^0,H)$ ($r=s,t_0$):
$$
	{\mathcal I}^{r+1}(t,V^0,H)\leq \uC ({\mathcal I}^{r+2}(t,U^0,G)
	+\langle \vert\uU\vert_{\tX^{r+7/2}}
	{\mathcal I}^{t_0+2}(t,U^0,G)\rangle_{r> t_0}).
$$

\noindent
{\bf Step 4.} The proposition follows from Steps 2 and 3.

\section{Main results}\label{sectmain}

\subsection{Large time existence for the water-waves equations}

In this section we prove the main result of this paper, which proves the
well-posedness of the water-waves equations over large times and provides
a uniform energy control which will allow us to justify all the
asymptotic regimes evoked in the introduction. Recall first that
the semi-normed spaces $(\tX^s,\vert\cdot\vert_{\tX^s})$ have been 
defined in Definition \ref{defispaces} as
$$
	\widetilde{X}^s:=\{(\zeta,\psi),\zeta\in H^s(\R^2),
	\nabla \psi\in H^{s-1/2}(\R^2)^2\},
$$
and $\vert(\zeta,\psi) \vert_{\widetilde{X}^s}
	:=\vert \zeta\vert_{H^s}+\vert\Pig\psi\vert_{H^s}$,
and define also the mapping ${\mathfrak a}$ by
\begin{eqnarray*}
	{\mathfrak a}(\zeta,\psi)&:=&
	\frac{\eps^2}{\nu}(\nag\psi-\eps\cZ[\eps\zeta]\nag\zeta)
	\cdot\nag\cZ[\eps\zeta]\psi\\
	&-&\eps\cZ[\eps\zeta](\zeta+\cA_2[(\zeta,\psi)])
	+\eps d_{\zeta}\cZ[\eps\cdot]\psi\cdot 
	\cG[\eps\zeta]\psi +1,
\end{eqnarray*}
where $\cZ[\eps\zeta]$ is as defined in Theorem \ref{theoshape} and
$\cA_2$ is defined in (\ref{A}), and where
	$\nu=(1+\sqrt{\mu})^{-1}$. The only condition we
impose on the parameters is that the steepness $\eps\sqrt{\mu}$
and the ratio $\beta/\eps$ remain bounded. More precisely, 
$(\eps,\mu,\gamma,\beta)\in {\mathcal P}_M$ ($M>0$) with
$$
	{\mathcal P}_M\!\!=\!\!\{(\eps,\mu,\gamma,\beta) \!\in \! (0,1]\times
	(0,\infty)\times (0,1]\times [0,1],\eps\sqrt{\mu} \! \leq \! M
	\mbox{ and } \frac{\beta}{\eps} \! \leq \! M\}.
$$
We can now state the theorem:
\begin{theo}\label{theomain}
	Let $t_0>1$, $M>0$  and ${\mathcal P}\subset {\mathcal P}_M$.\\
	There exists $P>D>0$ such that for all
	$s\geq s_0$, $b\in H^{s+P}(\R^2)$,  and all family 
	$(\zeta^0_p,\psi^0_p)_{p\in {\mathcal P}}$
	bounded in $\tX^{s+P}$ satisfying
	$$
	\inf_{\R^s} 1+\eps\zeta^0_p-\beta b>0
	\quad\mbox{ and }\quad
	\inf_{\R^2} {\mathfrak a}(\zeta^0_p,\psi^0_p)>0
	$$
	(uniformly with respect to 
	$p=(\eps,\mu,\gamma,\beta)\in {\mathcal P}$), there
	exist $T>0$ and a unique family 
	$(\zeta_p,\psi_p)_{p\in {\mathcal P}}$
	bounded in $C([0,\frac{\nu T}{\eps}];\tX^{s+D})$
	solving (\ref{nondimww}) with initial condition 
	$(\zeta^0_p,\psi^0_p)_{p\in {\mathcal P}}$.
\end{theo}
\begin{rema}
	The time interval of the solution varies 
	with $p\in {\mathcal P}$ (through $\eps$ and $\nu$); when
	we say that $(\zeta_p,\psi_p)_{p\in {\mathcal P}}$
	is bounded in  \\ 
        $C([0,\frac{\nu T}{\eps}];\tX^{s+D})$, we mean
	that there exists $C$ such that
	$$
	\forall p\in {\mathcal P},\quad
	\forall t\in [0,\frac{\nu T}{\eps}],\qquad
	\vert \zeta_p(t)\vert_{H^{s+D}}
	+\vert\Pig
	\psi_p(t)\vert_{H^{s+D}}\leq C.
	$$
\end{rema}
\begin{rema}
	For the shallow water
	regime for instance, one has $\eps=\beta=\gamma=1$ and
	$\mu$ is small (say, $\mu<1$); thus, we can take 
	${\mathcal P}=\{1\}\times (0,1)\times \{1\}\times \{1\}$;
	for the KP regime (with flat bottom), one takes
	${\mathcal P}=\{(\eps,\eps,\sqrt{\eps}),\eps\in (0,1)\}\times \{0\}$,
	etc.
\end{rema}
\begin{rema}
The condition $\inf_{\R^2}{\mathfrak a}(\zeta^0_p,\psi^0_p)>0$
is the classical Taylor sign condition proper to the water-wave equations
(\cite{Wu1,Wu2,LannesJAMS,Lindblad1,Lindblad2,AmbroseMasmoudi,CoutandShkoller,ShatahZeng}, among others). 
It is obviously true for small data and we give in 
Proposition \ref{proptaylor} some simple sufficient conditions on the initial
data and the bottom parameterization $b$ which ensure that it is satisfied.
\end{rema}
\begin{rema}\label{remstab}
	One also has the following stability property (see Corollary 1 in
	\cite{AlvarezLannes}): let $\underline{T}>0$ and
	$(U^{app}_p)_{p\in {\mathcal P}}
	=(\zeta^{app}_p,\psi^{app}_p)_{p\in{\mathcal P}}$,
	bounded in $Y^{s+P}_{\underline{T}}$ 
	(see Definition \ref{defispaces}), be an approximate solution
	of (\ref{nondimww}) in the sense that
	$$
	\dt U^{app}_p+\cL U^{app}_p+\frac{\eps}{\nu}\cA[U^{app}_p]=
	\frac{\eps}{\nu}\delta_p R_p,
	\qquad
	U^{app}_p\,_\init=(\zeta^0_p,\psi^0_p)+\delta_p r_p,
	$$
	with $(R_p,r_p)_{p}$ bounded in
	$C([0,\frac{\nu \underline{T}}{\eps}];X^{s+P})
	\cap C^1([0,\frac{\nu \underline{T}}{\eps}];X^{s+P-5/2})
	\times X^{s+P}$ (and $\delta_p\geq 0$). If moreover
	the $U^{app}_p$ are admissible, then one has
	$$
	\forall t\in [0,\frac{\nu}{\eps}\inf\{T,\underline{T}\}],\qquad
	\vert U_p(t)-U_p^{app}(t)\vert_{\tX^{s+D}}\leq \cst \delta_p,
	$$
where $U_p\in C([0,\frac{\nu T}{\eps}];\tX^{s+D})$ is the solution furnished
by the theorem. For $\delta_p$ small enough, one can take $T=\underline{T}$.
\end{rema}
\begin{rema}
	The numbers $P$ and $D$ could be explicited in the
	above theorem (as in Theorem 1 of \cite{AlvarezLannes} for instance),
	but since the focus here is not on the regularity of the solutions,
	we chose to alleviate the proof as much as possible. 
	For the same reason, we use a Nash-Moser iterative
	scheme which allows us to deal with all the different
	regimes at once, though it is possible in some cases
	to push further the analysis of the linearized operator
	and use a standard
	iterative scheme (as shown in \cite{Iguchi2} 
	for the shallow-water regime).
\end{rema}
\begin{proof}
Let us denote in this proof $\epsi=\eps/\nu$ and omit
the index $p$ for the sake of clarity. Rescaling the time by
$t\leadsto t/\epsi$ and using the same notations
as in (\ref{L}) and (\ref{A}), the theorem reduces to proving 
the well-posedness
of the IVP
$$
	\left\lbrace
	\begin{array}{l}
	\dt U+\frac{1}{\epsi}\cL U+\cA[U]=0,\\
	U\init=U^0,
	\end{array}\right.
$$
on a time interval $[0,T]$, with $T>0$ \emph{independent of all the
parameters}.\\
Define first the evolution operator $S^\eps(\cdot)$ associated to
the linear part of the above IVP. The following lemma shows that
the definition
\begin{equation}
	\label{defS}
	S^\epsi(t)U^0:=U(t), \quad\mbox{ with }\quad
	\dt U+\frac{1}{\epsi}\cL U=0
	\quad\mbox{ and }\quad U_\init=U^0
\end{equation}
makes sense for all data $U^0\in \tX^s$. 
\begin{lemm}
	For all $U^{0}\in \tX^{s}$, $S^\epsi(\cdot)U^0$ is
	well defined in $C([0,T]; \tX^s)$ by (\ref{defS}).
	Moreover, for all $0\leq t\leq T$,
	$$
 	 \vert S^{\epsilon} (t)U^0\vert_{\tX^s}
 	 \le C(T,\frac{1}{h_0}, |b|_{H^{s+7/2}}, \frac{\beta}{\varepsilon}, 
 	 \frac{\varepsilon}{\nu}) 
 	 \vert U^0\vert_{\tX^s}.
	$$
\end{lemm}
\begin{proof}
Proceeding as in the proof of Proposition \ref{propmain} (in the very simple case 
$\underline{U} = (0,0)^T$), one checks that $S^{\epsilon}(t)U^0$ makes
sense and that the estimate of the lemma holds if $U^0 \in X^s$.\\
Now, let us extend this result to data $U^{0}\in \tX^s$.
Let $\iota$ be a smooth function vanishing in a neighborhood of the 
origin and being constant equal to one outside the unit disc, and define, 
for all $\delta>0$, $\iota^{\delta}=\iota(|D|/\delta)$.  The couple 
$U^{0,\delta}:= (\zeta^0, \iota^{\delta} \psi^{0})^T
=(\zeta^0, \psi^{0,\delta})^T$ then belongs to 
$X^s$ and 
$U^{\delta}(t):= S^{\epsilon} (t)U^{0,\delta}
=(\zeta^{\delta}(t), \psi^{\delta}(t))^T$ is well defined in 
$X^s$.  Since 
$$
 \vert U^{\delta}(t)-U^{\delta'}(t)\vert_{\tX^s}
  \le C(T,\frac{1}{h_0}, |b|_{H^{s+7/2}}, \frac{\beta}{\varepsilon}, 
  \frac{\varepsilon}{\nu}) 
  \vert U^{0,\delta}-U^{0,\delta'}\vert_{\tX^s},
$$
it follows by dominated convergence that 
$(\zeta^{\delta})_{\delta\to 0}$ and $(\Pig \psi^{\delta})_{\delta\to 0}$ 
are Cauchy 
sequences in $C([0,T];H^s(\mathbb{R}^2))$. Therefore,
$(\zeta^{\delta}) \rightarrow \zeta$ and 
$(\Pig \psi^{\delta}) \rightarrow \omega$ 
in $C([0,T];H^s(\mathbb{R}^2))$, as $\delta$ goes to $0$. 
Defining $\psi(t)=\psi^{0}-\frac{1}{\epsilon} 
\int_0^t\zeta(t')dt'$ and using  
$\psi^{\delta}(t)=\psi^{0,\delta}-\frac{1}{\epsilon} 
\int_0^t\zeta^{\delta}(t')dt'$, one deduces $\omega=\Pig\psi$, from
which one infers that 
$\nabla \psi \in C([0,T];H^{s-1/2}(\mathbb{R}^2)^2)$. \\
From the convergence $\Pig\psi^\delta\to \omega=\Pig \psi$ in $H^s(\mathbb{R}^2)$ and 
Proposition \ref{unifcontrol} one deduces also that
$\GGuo\psi^\delta\to \GGuo\psi$ in $H^{s-1/2}(\mathbb{R}^2)$. One can thus take the
limit as $\delta\to 0$ in the relation
$\partial_t \zeta^{\delta}(t)-\frac{1}{\epsilon} \frac{1}{\mu\nu}
\GGuo \psi^{\delta}(t)=0$, thus proving that $(\zeta,\psi)\in \tX^s$ solves the IVP (\ref{defS}).
Since the solution to this IVP is obviously unique, this shows that
$S^\eps(\cdot )U^0$ makes sense in $\tX^s$ when $U^0\in \tX^s$.\\
The last assertion of the lemma follows by taking the 
limit when $\delta \rightarrow 0$ in the following expression 
$$
\vert S^{\epsilon} (t)U^{0,\delta}\vert_{\tX^s}
  \le C(T,\frac{1}{h_0}, |b|_{H^{s+7/2}}, \frac{\beta}{\varepsilon}, 
  \frac{\varepsilon}{\nu}) 
  \vert U^{0,\delta}\vert_{\tX^s}.  
$$ 
\qed
\end{proof}
We now look for the exact solution under the form
$U=S^\eps(t)U^0+V$, which is equivalent to
solving
\begin{equation}\label{IVPW}
	\left\lbrace
	\begin{array}{l}
	\dt V +\frac{1}{\epsi}\cL V+\cF[t,V]=h\\
	V\init=(0,0)^T,
	\end{array}\right.
\end{equation}
with $\cF[t,V]:=\cA[S^\epsi(t)U^0+V]-\cA[S^\eps(t)U^0]$ and 
$h :=-\cA[S^\eps(t)U^0]$.\\
We can now state two important properties satisfied by $\cL$
and $\cF$ (in the statement below, the notation $\cF^{(i)}_{(j)}$ means that
$\cF$ has been differentiated $i$ times with respect to time and $j$ with 
respect to its second argument).
\begin{lemm}\label{lemmcheck}
	Let $T>0$, $p=1$ and $m=5/2$. Then:\\
	{\bf i.} For all $s\geq t_0$, 
	the mapping $\cL: X^{s+m}\to X^s$ is well defined
	and continuous; moreover, the family of evolution 
	operators $(S^\eps(\cdot))_{0<\epsi<\epsi_0}$ is uniformly
	bounded in $C([-T,T];Lin(X^{s+m},X^s))$.\\
	{\bf ii.} For all $0\leq i\leq p$ and $0\leq i+j\leq p+2$, and for all
	$s\geq t_0+im$, one has 
	\begin{eqnarray*}
	\lefteqn{\sup_{t\in[0,T]}\vert 
	\epsi^i {\mathcal F}^{(i)}_{(j)}[t,U]
	(V_1,\dots,V_{j})
	\vert_{s-im}\leq
	C(s,T,\vert U\vert_{t_0+(i+1)m })}\\
	& &\times\big(\sum_{k=1}^{j} \vert V_k\vert_{s+m}
	\prod_{l\neq k }\vert V_l\vert_{t_0+(i+1)m}
	+\vert U\vert_{s+m}
	\prod_{k=1}^{j} \vert V_k\vert_{t_0+(i+1)m}\big),
	\end{eqnarray*}
	for all $U\in H^{s+m}(\R^2)$ and $(V_1,\dots, V_j)\in H^{s+m}(\R^2)^j$.
\end{lemm}
\begin{proof}
{\bf i.} The property on $S^\epsi(\cdot)$ follows from Proposition
\ref{proptrig} with $\uU=(0,0)$ (recall that we rescaled the time
variable). In order to prove the continuity of $\cL$, let us write,
for all $W=(\zeta,\psi)^T$,
\begin{eqnarray*}
	\vert\cL W\vert_{X^s} &\leq&
	\vert \frac{1}{\mu\nu}\cG[0]\psi\vert_{H^s}
	+\frac{\eps}{\nu}\vert \zeta \vert_{H^s}
	+\vert\Pig \zeta
	\vert_{H^s} \\
	&\leq&  \vert \frac{1}{\mu\nu}\cG[0]\psi\vert_{H^s}+
	C(\epsi)\vert \zeta\vert_{H^{s+1}}.
\end{eqnarray*}
One therefore deduces the continuity property on $\cL$  from the
following inequality:
$$
	\big\vert\frac{1}{\mu\nu}\cG[0]\psi
	\big\vert_{H^s}\leq
	C(\frac{1}{h_0},\eps\sqrt{\mu},\frac{\beta}{\eps},
	\vert b\vert_{H^{s+2}})
	\vert\Pig \psi
	\vert_{H^{s+1}};
$$
for $\mu\geq 1$, one has the uniform bound $\frac{1}{\mu\nu}\lesssim 1/\sqrt{\mu}$, and the inequality is a direct consequence of 
Proposition \ref{unifcontrol}; for $\mu\leq 1$, one has
$\nu\sim 1$ and we rather use 
Remark \ref{remlemfin}.\\
{\bf ii.} Since by definition 
${\mathcal F}[t,U]=\cA[S^\eps(t)U^0+U]-\cA(S^\eps(t)U^0)$,
it follows from the first point that it suffices to prove the estimates
in the case $i=0$ and with $\cF$ replaced by $\cA$. Recall that $\cA$ is
explicitly given by (\ref{A}) and remark that
\begin{eqnarray*}
	\cA_1[U]&=&-\frac{1}{\eps\mu}\int_0^1
	d_{ z\zeta}\cG[\eps\cdot]\psi\cdot \zeta  dz,\\
	&=&\int_0^1 \frac{1}{\sqrt{\mu}}\cG[\eps z\zeta]
	(\zeta\frac{1}{\sqrt{\mu}}\uZ)
	+z\nag\cdot(\zeta\bv)dz,
\end{eqnarray*}
where $\uZ$  and $\bv$ are as in Theorem \ref{theoshape} (with $\uz=\zeta$
and $\up=\psi$).\\
The estimates on $\cA$ are therefore a direct consequence of Proposition
\ref{unifcontrol}. \qed
\end{proof}
The well-posedness of (\ref{IVPW}) is deduced from the general
Nash-Moser theorem for singular evolution equations of \cite{AlvarezLannes}
(Theorem 1'),
provided that the three assumptions (Assumptions 1',2' and 3' in 
\cite{AlvarezLannes}) it requires are satisfied. The first
two, which concern the linear operator $\cL$ and the nonlinearity
$\cF[t,\cdot]$, are exactly the results stated in Lemma \ref{lemmcheck}.
The third assumption concerns the linearized operator around $\uV$ associated 
to (\ref{IVPW}); after remarking that
$$
	\dt+\frac{1}{\eps}\cL+d_{\uV}\cF[t,\cdot]
	=
	\cLu,
$$
with $\uU=(\uz,\up)^T=S^\epsi(t)U^0+\uV$, one can check that 
this last assumption
is exactly the result stated in Proposition \ref{propmain}, provided
that the following 
quantity (which is the first iterate of the Nash-Moser scheme,
see Remark 3.2.2 of \cite{AlvarezLannes})
\begin{equation}\label{U0}
	U_0:=t\mapsto S^\epsi(t)U^0+\int_0^t S^\epsi(t-t')\cF[t',U^0]dt'
\end{equation}
is an
admissible reference state in the sense of Definition \ref{defiadmissible}
on the time interval $[0,T]$ (recall that we rescaled the time variable).
Taking a smaller $T$ if necessary, it is sufficient to check the
admissibility at $t=0$, which is equivalent to the two assumptions
made in the statement of the theorem (after remarking that 
${\mathfrak a}(\zeta^0,\psi^0)=\mfa(t=0)$, with 
$\mfa$ as defined in (\ref{lin3}) and $\uU=U_0$). The proof is thus complete. 
\qed
\end{proof}

We end this section with a proposition showing that the Taylor sign condition
\begin{equation}\label{Taylorsign}
	\inf_{\R^2}{\mathfrak a}(\zeta^0_p,\psi^0_p)_{p\in {\mathcal P}}
	>0,
	\quad\mbox{ uniformly with respect to} \quad
	p\in {\mathcal P}
\end{equation}
can be replaced in Theorem \ref{theomain} 
by a much simpler condition. We need to introduce first 
the ``anisotropic Hessian''
${\mathcal H}_b^\gamma$ associated to the bottom parameterization $b$,
$$
	{\mathcal H}^\gamma_b:=
	\left(\begin{array}{cc} 
	\dx^2 b & \gamma^2\partial_{xy}^2 b\\
	\gamma^2\partial_{xy}^2b &\gamma^4 \dy^2 b
 	\end{array}\right)
$$ 
and the initial velocity potential $\Phi^0_p$ given by the BVP
\begin{equation}\label{bvpinit}
	\left\lbrace
	\begin{array}{l}
	\dsp \mu\partial_{x }^2\Phi^0_p
	+\gamma^2\mu\partial_{y }^2\Phi^0_p
	+\partial_{z }^2\Phi^0_p =0,
	\qquad -1+\beta b \leq z \leq 
	\eps\zeta^0_p,
	\vspace{1mm}\\
	\dsp \Phi^0_p\,_{\vert_{z=\eps\zeta^0_p}}=\psi^0_p,\qquad
	\partial_n\Phi^0_p\,_{\vert_{z=-1+\beta b}}=0.	
	\end{array}\right.
\end{equation}
\begin{prop}\label{proptaylor}
	Let $t_0>1$, $M>0$ and ${\mathcal P}\subset{\mathcal P}_M$; 
	let also $b\in H^{t_0+2}(\R^2)$,
	$(\zeta^0_p,\psi^0_p)_{p\in{\mathcal P}}$
	be  bounded in $\tX^{t_0+1}$ and
	$(\Phi^0_p)_{p\in{\mathcal P}}$ solve the BVPs (\ref{bvpinit}).
	Then, \\
	{\bf i.} There exists $\epsi_0>0$ such that 
	(\ref{Taylorsign})
	is satisfied if one replaces ${\mathcal P}$ by
	${\mathcal P}_{\eps_0}:=
	\{p=(\eps,\mu,\gamma,\beta)\in {\mathcal P},
	\eps\nu^{-1}\leq \epsi_0\}$;\\
	{\bf ii.} If there exist $\mu_1>0$ and 
	$\underline{\gamma}\in C((0,1]\times (0,\mu_1])$
	such that for all  $p=(\eps,\mu,\gamma,\beta)\in {\mathcal P}$ 
	one has $\mu\leq \mu_1$ and $\gamma=\underline{\gamma}(\eps,\mu)$,
	 and if
	$$
	-\eps^2\beta\mu{\mathcal H}^\gamma_b
	(\nabla\Phi^0_p\,_{\vert_{z=-1+\beta b}})\leq 1,
	\quad\mbox{ for all} \quad
	p\in {\mathcal P},
	$$
	then the Taylor sign condition (\ref{Taylorsign}) is satisfied.
\end{prop}
\begin{rema}
The first point of the proposition is used to check the Taylor condition
in deep water regime; in this latter case, one has indeed $\eps/\nu\sim \eps\sqrt{\mu}$ which is the \emph{steepness} of the wave, the small
parameter with respect to which asymptotic models are derived.\\
The second point of the proposition is essential in the shallow-water regime
(since $\eps/\nu$ does not go to zero as $\mu\to 0$). It is important to notice
that it implies that the Taylor condition \emph{is automatically satisfied for
flat bottoms}.
\end{rema}
\begin{rema}
S. Wu proved in  \cite{Wu1,Wu2} that the Taylor sign condition 
(\ref{Taylorsign}) is automatically satisfied in infinite depth; this
result was extended in \cite{LannesJAMS} to finite depth with flat bottoms.
The result needed here is stronger, since we want (\ref{Taylorsign}) to
be satisfied uniformly with respect to the parameters. In the $1DH$-case,
for flat bottoms, and in the particular case of the shallow-water regime, 
such a result was established in \cite{LiCPAM}.
\end{rema}
\begin{proof}
\emph{As in the proof of Theorem \ref{theomain}, we omit the
index $p$ to alleviate the notations}.\\
{\bf i.} As seen in the proof of Theorem \ref{theomain}, one has 
${\mathfrak a}(\zeta^0,\psi^0)=\mfa(t=0)$, where $\mfa$  is 
as defined in (\ref{lin3}) (with $\uU=U_0$ and $U_0$ given 
by (\ref{U0})). Thus, ${\mathfrak a}(\zeta^0,\psi^0)=1+\frac{\eps}{\nu}\mfb$
and $\vert \mfa\vert_{L^\infty}\geq 1-\epsi_0\vert \mfb\vert_{L^\infty}$.
It follows from Lemma \ref{lemmcontrol} that for the range of
parameters considered here, $\vert \mfb\vert_{L^\infty}$ is uniformly 
bounded on $[0,T]$, so that the result follows
 when $\epsi_0$ is small enough.\\
{\bf ii.} Step 1: \emph{ There exists $\mu_0>0$ such that (\ref{Taylorsign}) is
satisfied for all $p=(\eps,\mu,\gamma,\beta)\in {\mathcal P}$ 
such that $\mu\leq \mu_0$}. It is indeed a
consequence of Remark \ref{remlemfin} that $\vert \mfb\vert_{L^\infty}=O(\sqrt{\mu})$ as $\mu\to 0$; since moreover $\epsi=\frac{\eps}{\nu}$ remains
bounded, one can conclude as in the first step.\\
{Step 2.} \emph{The case $\mu\geq \mu_0$}.
For all time $t$, let $\Phi(t)$ denote the solution of the BVP (\ref{bvpinit}),
with the Dirichlet condition at the surface replaced by
$\Phi^0\,_{\vert_{z=\eps\zeta^0}}=\psi_0(t)$, 
where $U_0(t)=(\zeta_0(t),\psi_0(t))$ is given by (\ref{U0}).
Let us also define the ``pressure'' $P$ as
$$
	-\frac{1}{\eps}P:=
	\dt \Phi+\frac{1}{2}\big(\frac{\eps}{\nu}\vert \nag\Phi\vert^2
+\frac{\eps}{\mu\nu}(\dz \Phi)^2\big)+\frac{1}{\eps}z.
$$
Since $U_0=(\zeta_0,\psi_0)$ solves (\ref{nondimww}) at $t=0$,
one can check as in Proposition 4.4 of \cite{LannesJAMS} that
$P(t=0,X,\eps\zeta^0(X))=0$. Differentiating this relation
with respect to $X$ shows that 
$-\nag\zeta^0\cdot\nag P=\eps\vert\nag\zeta^0\vert^2\dz P$ on the 
surface, from which one deduces the identity,
$$
	(1+\eps^2\vert\nabla\zeta^0\vert^2)^{1/2}
	\partial_n P_{\vert_{z=\eps\zeta^0}}
	=
	(1+\eps^2\mu\vert\nag\zeta^0\vert^2)\dz 
	P_{\vert_{z=\eps\zeta^0}},
$$
where $\partial_n P_{\vert_{z=\eps\zeta^0}}$ stands for
the outwards conormal derivative associated to the elliptic
operator $\dsp \mu\partial_{x }^2
	+\gamma^2\mu\partial_{y }^2 
	+\partial_{z }^2$. Expressing $\Phi$ and its derivatives evaluated
at the surface in terms of $\Psi$, one can then
check that
\begin{equation}\label{mfadnp}
	\mfa(t=0)=-\frac{1}{1+\eps^2\mu^2\vert\nag\zeta^0\vert^2}
	(1+\eps^2\vert\nabla\zeta^0\vert^2)^{1/2}
	\partial_n P_{\vert_{z=\eps\zeta^0}}.
\end{equation}
Let us now remark that $P$ solves the BVP
$$
	\left\lbrace
	\begin{array}{l}
	\dsp (\mu\partial_{x }^2 
	+\gamma^2\mu\partial_{y }^2 
	+\partial_{z }^2)P=h ,
	\qquad -1+\beta b \leq z \leq 
	\eps \zeta^0,
	\vspace{1mm}\\
	\dsp P_{\vert_{z=\eps\zeta^0}}=0,\qquad
	\partial_n P_{\vert_{z=-1+\beta b}}=
	g,	
	\end{array}\right.
$$
with 
\begin{eqnarray*}
	h&:=&-\frac{1}{2}(\mu\partial_{x }^2 
	+\gamma^2\mu\partial_{y }^2 
	+\partial_{z }^2)(\frac{\eps^2}{\nu}\vert \nag\Phi\vert^2+\frac{\eps^2}{\mu\nu}(\dz\Phi)^2)\\
	g&:=&-\frac{1}{2}\partial_n(\frac{\eps^2}{\nu}\vert \nag\Phi\vert^2+\frac{\eps^2}{\mu\nu}(\dz\Phi)^2)_{\vert_{z=-1+\beta b}}
	-\partial_n(z)_{\vert_{z=-1+\beta b}}.
\end{eqnarray*}
Exactly as in 
the proof of  Proposition 4.15 of \cite{LannesJAMS},
one can check that $h\leq 0$ and use a maximum principle
(using the fact that (\ref{mfadnp}) links $\mfa$ to the normal
derivative of $P$ at the surface)
to show that if $g\leq 0$ then there exists a constant $c(\eps,\mu,\beta)>0$ 
such that
$\mfa(t_0)\geq c(\eps,\mu,\beta)$. 
We thus turn to prove that $g\leq 0$.\\
Recall that by construction of $\Phi$,
$$
	(1+\beta^2\vert\nabla b\vert^2)^{1/2}
	\partial_n\Phi_{\vert_{z=-1+\beta b }}
	\big(=
	\beta\mu \nag b\cdot\nag 
	\Phi_{\vert_{z=-1+\beta b }}
	-\dz \Phi_{\vert_{z=-1+\beta b }}
	\big)=0.
$$
Differentiating this relation with respect to $j$ ($j=x,y$), one gets
\begin{eqnarray*}
	& &(1+\beta^2\vert\nabla b\vert^2)^{1/2}
	\partial_n(\partial_j\Phi)_{\vert_{z=-1+\beta b }} =
	-\beta\mu \nag\partial_j b\cdot \nag\Phi_{\vert_{z=-1+\beta b }}\\
	& & +\beta\partial_j b 
	(1+\beta^2\vert\nabla b\vert^2)^{1/2}
	\partial_n(\dz \Phi)_{\vert_{z=-1+\beta b }},
\end{eqnarray*}
and using this formula one computes
\begin{eqnarray*}
	\lefteqn{\frac{1}{2}(1+\beta^2\vert\nabla b\vert^2)^{1/2}
	\partial_n(\frac{\eps^2}{\nu}\vert \nag\Phi\vert^2
	+\frac{\eps^2}{\mu\nu}
	(\dz\Phi)^2)_{\vert_{z=-1+\beta b}}}\\
	&=&-\frac{\eps^2\beta\mu}{\nu}
	\big((\dx\Phi)^2\dx^2 b+2\gamma^2\dx\Phi\dy\Phi\partial^2_{xy}b
	+\gamma^4(\dy\Phi)^2\dy^2 b\big)\\
	&=&- \frac{\eps^2\beta\mu}{\nu}{\mathcal H}^\gamma_b
	(\nabla\Phi_{\vert_{z=-1+\beta b}});
\end{eqnarray*}
since moreover $(1+\beta^2\vert\nabla b\vert^2)^{1/2}\partial_n(z)_{\vert_{z=-1+\beta b}}=1$, one gets $g\geq 0$ if the condition given in the
statement of the proposition is fulfilled.\\
As detailed above, we therefore have $\mfa(t=0)\geq c(\eps,\mu,\beta)>0$; 
moreover, there
exists by assumption $\mu_1$ such that for all 
$p=(\eps,\mu,\gamma,\beta)\in {\mathcal P}$,
one has $\mu\leq \mu_1$;  due to the first point of the proposition,
Step 1 and the 
fact the $\gamma=\underline \gamma(\eps,\mu)$, 
it is sufficient to
prove the proposition for all the parameters $p\in {\mathcal P}_1$ with
$$
	{\mathcal P}_1:=[\eps_0,1]\times [\mu_0,\mu_1]\times 
	\underline \gamma([\eps_0,1]\times [\mu_0,\mu_1])\times [0,1]
	\;\; (\eps_0:=(1+\sqrt{\mu_1})^{-1}\epsi_0).
$$
The dependence of $c(\eps,\mu,\beta)>0$ on $\eps$, $\mu$ and $\beta$
is continuous and therefore, 
$\inf_{[\eps_0,1]\times[\mu_0,\mu_1]\times [0,1]}c(\eps,\mu,\beta)>0$. \qed
\end{proof}

\section{Asymptotics for $3D$ water-waves}\label{sectjustif}

We will now provide a rigorous justification of the  main 
asymptotic models used in coastal oceanography.
\begin{rema}
	Throughout this section, we assume the following:\\
	- $P$ and $D$ are as in the statement of
	Theorem \ref{theomain};\\
	- $\Phi^0$ stands for the initial
	velocity potential as in Proposition \ref{proptaylor};\\
	- The bottom parameterization
	satisfies $b\in H^{s+P}(\R^2)$;\\
	- Except for the KP equations, we always consider fully transverse
	waves ($\gamma=1$), but one could easily use the methods set in
	this paper to derive and justify
	weakly transverse models in the other regimes. 
\end{rema}
\subsection{Shallow-water and Serre regimes}

We recall that the so-called  ``shallow-water'' regime corresponds
to the conditions $\mu\ll 1$ (so that $\nu\sim 1$) and 
$\eps=\gamma=1$; we also consider bottom
variations which can be of large amplitude ($\beta=1$).  
Without restriction, we can assume that $\nu=1$
(which corresponds to the nondimensionalization (\ref{nondimww})). 
The shallow-water model -- which goes back to Airy \cite{Airy}
and Friedrichs \cite{Friedrichs}-- consists of neglecting the $O(\mu)$ terms
in the water-waves equations, while the Green-Naghdi equations 
\cite{GreenLawsNaghdi,GreenNaghdi,SuGardner} is a more
precise approximation, which neglects only the $O(\mu^2)$ quantities. The
Serre equations \cite{Serre,SuGardner}
are quite similar to the Green-Naghdi equations, but
assume that the bottom and surface variations are
of medium amplitude: $\eps=\beta=\sqrt{\mu}$.

\subsubsection{The shallow-water equations}\label{sectjustifSW}
The shallow water equations are
\begin{equation}
	\label{shall1}
	\left\lbrace
	\begin{array}{l}
	\dt V+\nabla \zeta +\frac{1}{2}\nabla \vert V\vert^2=
	0,\\
	\dt \zeta+\nabla\cdot \big((1+\zeta-b) V\big)=0,
	\end{array}\right.
\end{equation}
and the following theorem shows that they provide a good approximation
to the exact solution of the water-waves equations.
\begin{theo}[Shallow-water equations]
\label{theoSW}
	Let $s\geq t_0>1$ and
	$(\zeta^{0}_\mu,\psi^{0}_\mu)_{0<\mu<1}$ be bounded
	in $\tX^{s+P}$. Assume moreover that
	there exist $h_0>0$ and $\mu_0>0$ such that for all 
	$\mu\in (0,\mu_0)$,
	$$
	\inf_{\R^2}(1+\zeta^0_\mu-b)\geq h_0
	\quad\mbox{ and }\quad
	-\mu{\mathcal H}^\gamma_b
	(\nabla\Phi^0_\mu\,_{\vert_{z=-1+b}})\leq 1.
	$$
	Then there exists $T>0$ and:
	\begin{enumerate}
	\item a unique family 
	$(\zeta_\mu,\psi_\mu)_{0<\mu<\mu_0}$ bounded in 
	$C([0,T];\tX^{s+D})$
	 and solving (\ref{nondimww}) with initial conditions
	$(\zeta^{0}_\mu,\psi_{\mu}^0)_{0<\mu<\mu_0}$;
	\item a unique family
	$(V^{SW}_\mu,\zeta^{SW}_\mu)_{0<\mu<\mu_0}$ bounded in \\ 
        $C([0,T];H^{s+P-1/2}(\R^2)^3)$ and solving (\ref{shall1})
	with initial conditions 
	$(\zeta^{0}_\mu,\nabla\psi_{\mu}^0)_{0<\mu<\mu_0}$.
	\end{enumerate}
	Moreover, one has, for some
	$C>0$,
	$$
	\forall 0<\mu<\mu_0, \;\;
	\vert \zeta_\mu-\zeta^{SW}_\mu\vert_{L^\infty([0,T]\times\R^2)}
	+\vert \nabla\psi_\mu-V^{SW}_\mu\vert_{L^\infty([0,T]\times\R^2)}
	\leq C \mu.
	$$
\end{theo}
\begin{rema}
The existence time provided by Theorem \ref{theomain} is $O(1)$, but is
large in the sense that it does not shrink to zero when $\mu\to 0$.
\end{rema}
\begin{rema}
	Instead of assuming that the initial data
	$(\zeta_\mu^0,\psi_\mu^0)_{0<\mu<1}$ are bounded
	in $\tX^{s+P}$, we could assume that
	$(\zeta_\mu^0,\nabla\psi_\mu^0)_{0<\mu<1}$ is
	bounded in $H^{s+P}(\R^2)^3$
	(because 
	$\vert\Pig\psi\vert_{H^{s+P}}
	\lesssim \vert \nabla\psi\vert_{H^{s+P}}$,
	uniformly with respect to $\mu\in (0,1)$). 
\end{rema}
\begin{rema}
	The $2DH$ shallow-water model has been justified rigorously 
	by Iguchi in a recent work \cite{Iguchi2}, but under
	two restrictions: a) The velocity potential is assumed to have
	Sobolev regularity which implies that the velocity
	must satisfy some restrictive zero mass assumptions and b) the 
	theorem holds only for very small values of $\mu$.
	 These assumptions are removed
	in the above result.
\end{rema}
\begin{proof}
The assumptions allow us to use Theorem \ref{theomain} and 
Proposition \ref{proptaylor} with ${\mathcal P}=\{1\}\times (0,\mu_0)\times
\{1\}\times\{1\}$, which proves the first part of the theorem.\\
The second point of the theorem is straightforward since  \\
$\vert \nabla\psi_\mu^0\vert_{H^{s+P-1/2}}\leq 
\vert \Pig\psi^0_\mu\vert_{H^{s+P}}$ (recall that $\mu<1$ here), and because
(\ref{shall1}) is a 
quasilinear hyperbolic system (since $\inf_{\R^2}(1+\zeta-b)>0$).
In order to prove the error estimate, 
plug the expansion furnished by Proposition \ref{propshallow} into
(\ref{nondimww}) and take the gradient of the second equation in
order to obtain a system of equations on $\zeta_\mu$ and 
$V_\mu=\nabla\psi_\mu$.
One gets
\begin{equation}
	\label{shall2}
	\left\lbrace
	\begin{array}{l}
	\dt V_\mu+\nabla \zeta_\mu +\frac{1}{2}\nabla \vert V_\mu\vert^2=
	\mu R^1_\mu,\\
	\dt \zeta_\mu+\nabla\cdot \big((1+\zeta_\mu-b) V_\mu\big)=\mu R^2_\mu,
	\end{array}\right.
\end{equation}
with $(R^1_\mu,R^2_\mu)$ uniformly bounded in 
$L^\infty([0,T];H^{t_0}(\R^2)^{2+1})$. An energy estimate on (\ref{shall2})
thus gives a Sobolev  error estimate from which one
deduces the $L^\infty$ estimate of the theorem using
 the classical continuous embedding $H^{t_0}\subset L^\infty$. \qed
\end{proof}
\subsubsection{The Green-Naghdi and Serre equations}\label{sectjustifGNS}

Though corresponding to two different physical regimes, 
the Green-Naghdi and Serre equations can both be written at the same time
if one assumes that $\eps=1$ for the Green-Naghdi equations and 
$\eps=\sqrt{\mu}$ for the Serre equations in the formulation below:
\begin{equation}\label{GN0}
	\left\lbrace
	\begin{array}{l}
	(h+\mu{\mathcal T}[h,\eps b])\dt V+h\nabla\zeta
	+\eps h(V\cdot\nabla)V\\
	\qquad \qquad+\mu\eps  \Big[
	\frac{1}{3}\nabla\big(h^3{\mathcal D}_{V}\mbox{div}(V)\big)+
	\cQ[h,   \eps b](V)\Big]=0	\\
	\dt\zeta+\nabla\cdot(hV)=0,
	\end{array}\right.
\end{equation}
where $h:=1+ \eps(\zeta-b)$ while the linear operators
${\mathcal T}[h,b]$ and ${\mathcal D}_{V}$ and the quadratic form 
$\cQ[h,b](\cdot)$ are 
defined as
\begin{eqnarray*}
	{\mathcal T}[h,b] V&:=&
	-\frac{1}{3}\nabla(h^3\nabla\cdot V)
	+\frac{1}{2}\big[
	\nabla(h^2\nabla b \cdot V)-h^2\nabla b \nabla\cdot V\big] \\
        & & +h\nabla b\nabla b\cdot V,\\
	{\mathcal D}_V&:=&-(V\cdot \nabla)+\mbox{div}(V), \\
	\cQ[h,b](V)&:=&	\frac{1}{2}\nabla\big(h^2(V\cdot\nabla)^2b\big)+
	h\big(\frac{h}{2}{\mathcal D}_{V}\mbox{div}(V)+(V\cdot\nabla)^2b\big)
	\nabla b.
\end{eqnarray*}
Both the Green-Naghdi and Serre models are rigorously justified in the
theorem below:
\begin{theo}[Green-Naghdi and Serre equations]
	\label{theoGN}
	Let $s\geq t_0>1$ and
	$(\zeta^{0}_\mu,\psi^{0}_\mu)_{0<\mu<1}$ be bounded
	in $\tX^{s+P}$. Let $\eps=1$ (Green-Naghdi)
	or $\eps=\sqrt{\mu}$ (Serre) and assume that for
	some $h_0>0$, $\mu_0>0$ and for all $\mu\in (0,\mu_0)$,
	$$
	\inf_{\R^2}(1+\eps(\zeta^0_\mu-b))\geq h_0
	\quad\mbox{ and }\quad
	-\mu\eps^3{\mathcal H}^\gamma_b
	(\nabla\Phi^0_\mu\,_{\vert_{z=-1+\eps b}})\leq 1.
	$$
	Then there exists $T>0$ and:
	\begin{enumerate}
	\item a unique family 
	$(\zeta_\mu,\psi_\mu)_{0<\mu<\mu_0}$ bounded in 
	$C([0,\frac{T}{\eps}];\tX^{s+D})$
	 and solving (\ref{nondimww}) with initial conditions
	$(\zeta^{0}_\mu,\psi_{\mu}^0)_{0<\mu<\mu_0}$;
	\item a unique 
	family $(V^{GN}_\mu,\zeta^{GN}_\mu)_{0<\mu<\mu_0}$ bounded
	in $C([0,\frac{T}{\eps}];H^{s}(\R^2)^3)$ and solving (\ref{GN0})
	with initial conditions $(\zeta^{0}_\mu,(1-\frac{\mu}{h^0}
	{\mathcal T}[h^0,\eps b])
	\nabla\psi_{\mu}^0)$ (with $h^0=1+\eps(\zeta^0-b)$).
	\end{enumerate}
	Moreover, one has for some
	$C>0$ independent of $\mu\in (0,\mu_0)$,
	$$
	\vert \zeta_\mu-\zeta^{GN}_\mu\vert_{L^\infty([0,\frac{T}{\eps}]
	\times\R^2)}
	+\vert \nabla\psi_\mu-(1+\frac{\mu}{h}{\mathcal T}[h,\eps b])
	V^{GN}_\mu\vert_{L^\infty([0,\frac{T}{\eps}]\times\R^2)}
	\leq C \frac{\mu^2}{\eps}.
	$$
\end{theo}
\begin{rema}
	The precision of the GN approximation ($\eps=\beta=1$) 
	is therefore one order
	better than the shallow-water equations. This model had been 
	justified in $1DH$ and for flat bottoms by Y. A. Li \cite{LiCPAM}.
	The theorem above is stated in $2DH$ but one can cover the open case
	of $1DH$ non-flat bottoms with a straightforward adaptation.
\end{rema}
\begin{rema}
	In the Serre scaling, one has $\eps=\sqrt{\mu}$, and the precision
	of the theorem is therefore $O(\mu^{3/2})$, which is worse than the
	$O(\mu^2)$ precision of the GN model, but the approximation
	remains valid over a larger time scale (namely, $O(\mu^{-1/2})$ versus
	$O(1)$ for GN). Notice also that at first order in $\mu$,
	the Serre equations reduce to a simple wave equation (speed $\pm 1$)
	on $\zeta$ and $V$, which is not the case for GN where
	the shallow-water equations (\ref{shall1}) are found at first order.
\end{rema}
\begin{proof}
The first assertion of the theorem is exactly the same as in 
Theorem \ref{theoSW} in the GN case. For the Serre equations,
it is also a direct consequence of Theorem \ref{theomain} and Proposition
\ref{proptaylor}, with 
${\mathcal P}=\{(\sqrt{\mu},\mu,1,\sqrt{\mu}),\mu\in (0,\mu_0)\}$.\\
For the second assertion, we replace $\G$ 
in (\ref{nondimww}) by the expansion given in  Proposition \ref{propshallow}
and take the gradient of the equation on $\psi$ 
to obtain
\begin{equation}
	\label{GN1}
	\left\lbrace
	\begin{array}{l}
	\dt V_\mu+\nabla \zeta_\mu +\frac{\eps}{2}\nabla \vert V_\mu\vert^2
	-\frac{\eps\mu}{2}\nabla(h\nabla\cdot V_\mu-\eps\nabla b\cdot V_\mu)^2=
	\mu^2 R^1_{\mu},\\
	\dt \zeta_\mu+\nabla\cdot \big(h V\big)=\mu^2 R^2_\mu,
	\end{array}\right.
\end{equation}
with $(R_\mu^1,R^2_\mu)_\mu$ bounded in 
$L^\infty([0,\frac{T}{\eps}];H^{t_0}(\R^2)^{2+1})$ while $V$ is defined as
$
	V:=V_\mu-\frac{\mu}{h}{\mathcal T}[h,b]V_\mu,
$
so that $V_\mu=V+\frac{\mu}{h}{\mathcal T}[h,b]V+O(\mu^2)$. 
Replacing $V_\mu$ by this expression in (\ref{GN1}) and neglecting the
$O(\mu^2)$ terms then gives (\ref{GN0}). The theorem is then a direct
consequence of the well-posedness theorem for the Green-Naghdi and Serre equations proved
in \cite{AlvarezLannes} and of the error estimates given in Theorem 3
of that reference. \qed
\end{proof}

\subsection{Long-waves regime: the Boussinesq approximation}
\label{sectjustifLW}

The long-wave regime is characterized by the scaling $\gamma=1$, 
$\mu=\eps\ll 1$, so
that one has $\nu\sim 1$. As for the shallow-water equations, we take $\nu=1$
for notational convenience. 
When the bottom is non-flat, it is assumed
that its variations are of the order of the size of the waves, that is,
$\beta=\eps$. Since the pioneer work of Boussinesq
\cite{Boussinesq}, many formally equivalent systems, generically
called Boussinesq systems, have been derived to model the dynamics of the
waves under this scaling. Following \cite{BonaSmith}, these systems where 
derived in a systematic way in \cite{BCS,BCL,Chen,Chazel}. 
In \cite{BCL,Chazel} some interesting
\emph{symmetric} systems where introduced:
$$
	S'_{\theta,p_1,p_2}
	\left\lbrace
	\begin{array}{l}
	\displaystyle
	(1-\eps a_2\Delta)\dt V+\nabla \zeta+
	\eps\big(
	\frac{1}{4}\nabla\vert V\vert^2+\frac{1}{2}(V\cdot \nabla)V
	+\frac{1}{2}V \nabla\cdot V\vspace{1mm}\\
	\displaystyle
	\qquad\qquad\qquad\qquad\qquad+\frac{1}{4}\nabla\vert\zeta\vert^2
	-\frac{1}{2}b\nabla\zeta+a_1\Delta\nabla\zeta\big)=0,\vspace{1mm}\\
	\dsp (1-\eps a_4\Delta) \dt \zeta+\nabla \! \cdot \! V+\frac{\eps}{2}\big(
	\nabla \! \cdot \! \big((\zeta-b) V\big)+a_3\Delta\nabla\cdot V)=0,
	\end{array}\right.
$$
where the coefficients $a_j$ ($j=1,\dots,4$) 
depend on $p_1,p_2\in \R$ and $\theta\in [0,1]$ through the relations
$a_1=(\frac{\theta^2}{2}-\frac{1}{6})p_1$,
$a_2=(\frac{\theta^2}{2}-\frac{1}{6})(1-p_1)$,
$a_3=\frac{1-\theta^2}{2}p_2$,
and $a_4=\frac{1-\theta^2}{2}(1-p_2)$;
some choices of parameters yield $a_1=a_3$ and $a_2\geq 0$, $a_4\geq 0$,
and the 
corresponding systems $S'_{\theta,p_1,p_2}$ are the completely
symmetric systems mentioned above.
The so-called Boussinesq approximation associated to
a family of initial data 
$(\zeta^0_\eps,\psi^0_\eps)_{0<\eps<1}$ is given by
\begin{equation}
	\label{appB}
	\zeta_\eps^{app}=\zeta_\eps^{B}
	\quad\mbox{ and }\quad
	V_\eps^{app}=(1-\frac{\eps}{2}(1-\theta^2)\Delta)
	\big(1-\frac{\eps}{2}(\zeta_\eps^{B}-b)V^{B}_\eps\big),
\end{equation}
where
$(V^{B}_\eps,\zeta_\eps^{B})_{0<\eps<1}$ 
solves $S'_{\theta,p_1,p_2}$ with initial data
\begin{equation}\label{initBouss}
	V^{B,0}_{\eps}=\big(1+\frac{\eps}{2}(\zeta_{\eps}^0-b)\big)
	\big(1-\frac{\eps}{2}(1-\theta^2)\Delta\big)^{-1}\nabla\psi_{\eps}^0
	\quad\mbox{ and }\quad	
	\zeta^{B,0}_\eps=\zeta_\eps^0.
\end{equation}
The following theorem fully justifies this approximation.
\begin{theo}[Boussinesq systems]
\label{theoBouss}
	Let $s\geq t_0>1$ and \\
	$(\zeta^{0}_\eps,\psi^{0}_\eps)_{0<\eps<1}$ be bounded
	in $\tX^{s+P}$ and assume that there exist $h_0>0$ and
	$\eps_0>0$ such that for all $\eps\in (0,\eps_0)$,
	$$
	\inf_{\R^2}(1+\eps(\zeta^0_\eps-b))\geq h_0
	\quad\mbox{ and }\quad
	-\eps^4{\mathcal H}^\gamma_b
	(\nabla\Phi^0_\mu\,_{\vert_{z=-1+\eps b}})\leq 1.
	$$
	Then there exists $T>0$ and:
	\begin{enumerate}
	\item a unique family 
	$(\zeta_\eps,\psi_\eps)_{0<\eps<\eps_0}$ bounded in 
	$C([0,\frac{T}{\eps}];\tX^{s+D})$
	 and solving (\ref{nondimww}) with initial conditions
	$(\zeta^{0}_\eps,\psi^{0}_{\eps})_{0<\eps<\eps_0}$;
	\item a unique family
	$(V^{B}_\eps,\zeta^{B}_\eps)_{0<\eps<\eps_0}$ bounded
	in $C([0,\! \frac{T}{\eps}];H^{\!s+P-\frac12}(\R^2)^3)$ and solving
	$S'_{\theta,p_1,p_2}$ with initial conditions (\ref{initBouss}).
	\end{enumerate}
	Moreover, for some $C>0$ independent of $\eps\in (0,\eps_0)$, one has
	$$
	\forall 0\leq t\leq \frac{T}{\eps},
	\qquad
	\vert \zeta_\eps(t)-\zeta^{app}_\eps(t)\vert_{\infty}
	+\vert \nabla\psi_\eps(t)-V_\eps^{app}(t)\vert_\infty	
	\leq C\eps^2 t,
	$$
	where $(V^{app}_\eps,\zeta^{app}_\eps)$ is given by (\ref{appB}).
\end{theo}
\begin{rema}
	The above theorem justifies \emph{all} the Boussinesq systems
	and not only  the completely
	symmetric Boussinesq systems considered here: 
	it is proved in \cite{BCL,Chazel} 
	that the justification of \emph{all} the Boussinesq systems
	follows directly from the justification of \emph{one} of them, 
	in the sense that their
	solutions (if they exist!) provide an approximation of 
	order $O(\eps^2t)$ to the water-waves equations.
\end{rema}
\begin{proof}
The ``if-theorems'' of \cite{BCL,Chazel} prove the result assuming that
the first statement of the theorem holds, which is a direct consequence
of Theorem \ref{theomain} and Proposition \ref{proptaylor} with
${\mathcal P}=\{(\eps,\eps,1,\eps),\eps\in (0,\eps_0)\}$. \qed
\end{proof}

\subsection{Weakly transverse long-waves: the KP  approximation}
\label{sectjustifKP}

We recall that the KP regime is the same as the long-waves regime,
but with $\gamma=\sqrt{\eps}$. Moreover, we assume here that the
bottom is flat, for the sake of simplicity. The KP
approximation \cite{KP}
consists in replacing the exact water elevation $\zeta_\eps$
by the sum of two counter propagating waves, slowly modulated
by a KP equation; more precisely, one defines $\zeta^{KP}_\eps$
as
\begin{equation}\label{appKP}
	\zeta^{KP}_\eps(t,x)=\frac{1}{2}\big(
	\zeta_+(\eps t, \sqrt{\eps}y,x-t)
	+\zeta_-(\eps t, \sqrt{\eps}y,x+t)\big)
\end{equation}
where $\zeta_\pm(\tau,Y,X)$ solve the KP equation
$$
	\partial_\tau \zeta_\pm \pm 
	\frac{1}{2}\partial_X^{-1}\partial_Y^2\zeta_\pm 
	\pm \frac{1}{6}\partial_X^3\zeta_\pm 
	\pm \frac{3}{2}\zeta_\pm \partial_X\zeta_\pm =0.
	\qquad (KP)_\pm
$$
This approximation is rigorously justified in the theorem below:
\begin{theo}[KP equation]
	Let $s\geq t_0>1$ and
	$(\zeta^{0},\psi^{0})\in \tX^{s+P}$ and assume that 
	there exist $h_0>0$ and $\eps_0>0$ such that
	for all $\eps\in(0,\eps_0)$,
	$$
	\inf_{\R^2}(1+\eps\zeta^0)\geq h_0,
	$$
	and assume also that
	$(\dy^2\dx\psi^{0},\dy^2\zeta^{0})\in \partial_x^2 H^{s+P}(\R^2)^2$.\\
	Then there exists $T>0$ and:
	\begin{enumerate}
	\item a unique family 
	$(\zeta_\eps,\psi_\eps)_{0<\eps<\eps_0}$
	solving (\ref{nondimww}) with initial conditions
	$(\zeta^{0},\psi^{0})$ and such that
	$(\zeta_\eps)_{0<\eps<\eps_0}$, $(\dx\psi_\eps)_{0<\eps<\eps_0}$ and 
	$(\sqrt{\eps}\dy \psi_\eps)_{0<\eps<\eps_0}$ are bounded in 
	$C([0,\frac{T}{\eps}];H^{s+D-1/2})$;
	\item a unique solution 
	$\zeta_\pm\in C([0,T];H^{s+P-1/2}(\R^2))$ to
	(KP)$_\pm$ with initial condition $(\zeta^0\pm\dx\psi^0)/2$.
	\end{enumerate}
	Moreover, one has the following error estimate for the 
	approximation (\ref{appKP}):
	$$
	\lim_{\eps\to 0}
	\vert \zeta_\eps -\zeta^{KP}_\eps\vert_{L^\infty([0,\frac{T}{\eps}]
	\times\R^2)}	
	=0.
	$$
\end{theo}
\begin{rema}
	The very restrictive ``zero mass'' assumptions that \\
	$\dy^2\dx\psi^0$ and $\dy^2\zeta^0$ are twice the derivative
	of a Sobolev function comes from the singular component
	$\partial_X^{-1}\partial_Y^2$ of the KP equations (KP)$_\pm$.  
	Furthermore, the error estimate is much worse than for the Boussinesq
	approximations. These two drawbacks are removed if one replaces
	the KP approximation by the approximation furnished by the
	\emph{weakly transverse Boussinesq systems} introduced 
	in \cite{LannesSaut}. As shown in \cite{LannesSaut}, the first
	assertion of the above theorem rigorously justify these systems: 
	they provide an approximation
	of order $O(\eps^2 t)$ on the time interval $[0,T/\eps]$, and
	do not require the ``zero mass'' assumptions.
\end{rema}
\begin{proof}
As for the Boussinesq systems, we only have to prove the first assertion
of the theorem, and the whole result then follows from the ``if-theorem''
of \cite{LannesSaut}. Taking 
${\mathcal P}=\{(\eps,\eps,\sqrt{\eps},0),\eps\in(0,\eps_0)\}$, 
Theorem \ref{theomain} and Proposition \ref{proptaylor} 
give a family of solutions $(\zeta_\eps,\psi_\eps)_{0<\eps<\eps_0}$ 
bounded in 
$C([0,\frac{T}{\eps}];\tX^{s+D})$. In particular, 
$(\vert \Pig\psi_\eps\vert_{H^{s+D}})_\eps$ is bounded, and
thus $(\vert \nag\psi_\eps\vert_{H^{s+D-1/2}})_\eps$ is also bounded. 
Since $\gamma=\sqrt{\eps}$, one has $\vert \dx\psi\vert_{H^{s+D-1/2}}+\sqrt{\eps}\vert \dy\psi_\eps\vert_{H^{s+D-1/2}}
\lesssim \vert \nag\psi_\eps\vert_{H^{s+D-1/2}}$ and the claim follows.  \qed
\end{proof}

\subsection{Deep water}

\subsubsection{Full dispersion model}\label{sectjustifFD}

We present here the so-called full dispersion (or Matsuno) model for
deep water-waves. Contrary to all the asymptotic models seen above, 
the shallowness parameter $\mu$ is allowed to take large values
(deep water) provided that the \emph{steepness} of the waves $\eps\sqrt{\mu}$
remains small; without restriction, we can therefore take 
$\nu=\mu^{-1/2}$ here (i.e., we use the nondimensionalization 
(\ref{nondimdeep})). Introducing $\epsi=\eps\sqrt{\mu}$ 
the full dispersion model derived 
in \cite{Matsuno1,Matsuno2,Choi} can be written in the case of flat 
bottoms ($\beta=0$):
\begin{equation}\label{matsuno}
	\left\lbrace
	\begin{array}{l}
	\dt \zeta -\cT_\mu V+\epsi\big(\cT_\mu(\zeta\nabla\cT_\mu V)
	+\nabla\cdot(\zeta V)\big)=0,\vspace{1mm}\\
	\dt V+\nabla\zeta+\epsi\big(\frac{1}{2}\nabla \vert V\vert^2
	-\nabla\zeta\cT_\mu\nabla\zeta\big)=0,
	\end{array}\right.
\end{equation}
where $\cT_\mu$ is a Fourier multiplier defined as
$$
	\forall V\in {\mathfrak S}(\R^2)^2,\qquad
	\widehat{\cT_\mu V}(\xi)=-
	\frac{\tanh(\sqrt{\mu}\vert\xi\vert)}{\vert\xi\vert}(i\xi)\cdot \widehat{V}(\xi).
$$
Since $\beta=0$ (flat bottom) and $\gamma=1$ (fully transverse), the full
dispersion model depends on two parameters $(\eps,\mu)$ which are linked by a
small steepness assumption:
	$$
	\exists \epsi_0>0,\qquad
	(\eps,\mu)\in {\mathcal P}_{\epsi_0}\subset\{(\eps,\mu)\in (0,1]\times [1,\infty),
	\epsi:=\eps\sqrt{\mu}\leq \epsi_0\}.
	$$
The well-posedness of the full-dispersion model has not been investigated yet,
but we can prove that if a solution exists on 
$[0,\frac{\underline{T}}{\epsi}]$ ($\epsi>0$ small enough), 
then the solution of the water-waves equations exists over the same time
interval and is well approximated by the solution of the full-dispersion model:
\begin{theo}[Full-dispersion model]\label{theodeep}
	Let $\epsi_0>0$, $Q\geq P$ large enough, $s\geq t_0>1$ and
	$(\zeta^{0},\psi^{0})\in \tX^{s+P}$, and assume that 
	$$
	\forall \eps \in (0,1],\qquad
	\inf_{\R^2}(1+\eps(\zeta^0-b))\geq h_0>0.
	$$
	Let also $\underline{T}>0$ and let $(\zeta^{FD}_{\eps,{\mu}},
	V^{FD}_{\eps,\mu})_{(\eps,\mu)\in {\mathcal P}_{\epsi_0}}$
	be bounded in \\
        $C([0,\frac{\underline{T}}{\epsi}], H^{s+Q}(\R^2)^3)$
	and solving (\ref{matsuno}) with initial condition 
	$(\zeta^0,\nabla\psi^0-\epsi(\cT_\mu \nabla\psi^0)\nabla\zeta^0)$.\\
	Then, if $\epsi_0$ is small enough, there is 
	a unique family 
	$(\zeta_{\eps,\mu},
	\psi_{\eps,\mu})_{(\eps,\mu)\in {\mathcal P}_{\epsi_0}}$
	bounded in 
	$C([0,\frac{\underline{T}}{\epsi}];\tX^{s+D})$ 
	 and solving (\ref{nondimww}) with initial conditions
	$(\zeta^{0},\psi^{0})$. In addition,  
	for some $C>0$ independent of $(\eps,\mu)\in {\mathcal P}_{\epsi_0}$, 
	one has
	$$
	\vert \zeta_{\eps,\mu}-\zeta^{FD}_{\eps,\mu}\vert_{L^\infty([0,\frac{\underline{T}}{\epsi}]\times \R^2)}
	+\vert \nabla\psi_\eps-V_\eps^{FD}
	\vert_{L^\infty([0,\frac{\underline{T}}{\epsi}]\times \R^2)}	
	\leq C\epsi \qquad (\epsi:=\eps\sqrt{\mu}).
	$$
\end{theo}
\begin{proof}
Since  $\cT_\mu:H^r(\R^2)^2\mapsto H^r(\R^2)$ is continuous with operator norm
bounded from above by $1$, the mapping 
$V\in H^r(\R^2)^2\mapsto V-\epsi (\cT_\mu V)\nabla\zeta\in H^r(\R^2)^2$ 
is continuous for all $r\geq t_0$ and 
$\zeta\in H^{r+1}(\R^2)$. 
Moreover, this mapping is invertible for $\epsi$ small
enough, and one can accordingly define $\tV:=(1-\epsi \nabla\zeta \cT_\mu)^{-1}V $, 
so that $V=\tV-\epsi(\cT_\mu \tV)\nabla\zeta$. 
Replacing $V$ by this expression in (\ref{matsuno})
gives
\begin{equation}\label{matsuno1}
	\left\lbrace
	\begin{array}{l}
	\dt \zeta-{\mathcal T}_\mu \tV+\epsi
	\big(\nabla\cdot(\zeta \tV)+{\mathcal T}_\mu\nabla(\zeta{\mathcal T}_\mu \tV) \big)
        =\epsi^2 r^1_\epsi\vspace{1mm},\\
	\dt \tV+\nabla\zeta +\epsi\frac{1}{2}\big(\nabla \vert \tV\vert^2
        -\nabla({\mathcal T}_\mu\tV)^2\big)=\epsi^2\nabla r^2_\epsi,
	\end{array}\right.
\end{equation}
and where the exact expression of $R_\epsi:=(r^1_\epsi,\nabla r^2_\epsi)_\epsi$
 is of
no importance. Now, let $\dt+L$ denote the linear part of the above system
and $S(t)$ its evolution operator:
$L:=\left(	
	\begin{array}{cc}
	0 & -\cT_\mu\\
	\nabla & 0
	\end{array}\right)
$, and for all $U=(\zeta,V)$, $S(t)U:=u(t)$,  
where  $u$ solves $(\dt +L)u=0$, with initial condition $u_\init=U$. Since
$\cT_\mu$ is a Fourier multiplier, one can find an explicit expression
for $S(t)$, but we only need the following property: $S(t)$ is unitary
on $Z^r$ ($r\in\R$) defined as
\begin{eqnarray*}
	Z^r &:=& \{U=(\zeta,V)\in H^r(\R^2)^3, \\
        & & \vert U\vert_{Z^r}:=
	\vert \zeta\vert_{H^r}+
	\big\vert \big(\frac{\tanh(\sqrt{\mu}\vert\xi\vert)}{\vert\xi\vert}
	\big)^{1/2}V\big\vert_{H^r}<\infty\}.
\end{eqnarray*}
Writing $\widetilde{u}:=(\zeta,\tV)$, we define $w:=(\tz,W)$
as
$$
	w:=\widetilde{u}-\epsi^2\int_0^t S(t-t')R_\epsi(t')dt';
$$
remarking that 
$\vert V\vert_{H^{r-1/2}}\lesssim
	\big\vert 
	\big(\frac{\tanh(\sqrt{\mu}\vert\xi\vert)}{\vert\xi\vert}\big)^{1/2}
	V\big\vert_{H^r}
$ 
and that
$$
	\forall f\in H^{r+1}(\R^2),\qquad
	\vert \big(
	\frac{\tanh(\sqrt{\mu}\vert\xi\vert)}{\vert\xi\vert}\big)^{1/2}
	\nabla f\vert_{H^{r+1/2}}\lesssim \vert f\vert_{H^{r+1}},
$$
uniformly with respect to $\mu\geq 1$, and since $S(t)$ is unitary on $Z^r$, 
one gets
\begin{equation}\label{estdiff}
	\forall r\geq 0,\qquad
	\sup_{[0,\frac{\underline{T}}{\epsi}]}
	\vert w(t)-\widetilde{u}(t)\vert_{H^r}\lesssim
	\epsi \underline{T}\sup_{[0,\frac{\underline{T}}{\epsi}]}
	(\vert r^1_\epsi\vert_{H^{r+1/2}}+\vert r^2_\epsi\vert_{H^{r+1}}).
\end{equation}
Furthermore, one immediately checks that $w$ solves
$$
	\left\lbrace
	\begin{array}{l}
	\dt \tz -{\mathcal T}_\mu W+\epsi
	f^1(\tz,W)=\epsi^2 k^1_\epsi\vspace{1mm},\\
	\dt W+\nabla\tz+\epsi \nabla f^2(\tz,W)=\epsi^2\nabla k_\epsi^2,
	\end{array}\right.
$$
with initial condition $w_\init=(\zeta^0,\nabla\psi^0)^T$, 
and where
$k^1_\epsi:=\frac{1}{\epsi}(f^1(w)-f^1(\widetilde{u}))$, $k^2_\epsi:=\frac{1}{\epsi}(f^2(w)-f^2(\widetilde{u}))$, and
\begin{eqnarray*}
	f^1(\zeta,V):=\nabla\cdot(\zeta V)+{\mathcal T}_\mu\nabla(\zeta{\mathcal T}_\mu V),\qquad
	f^2(\zeta,V):=\frac{1}{2}
	\big(\vert V\vert^2-({\mathcal T}_\mu V)^2\big);
\end{eqnarray*}
from (\ref{estdiff}), one gets in particular that
\begin{equation}\label{bound}
	\vert (k_\epsi^1,k_\epsi^2)\vert_{X^{s+P}_{\underline{T}}}
	+\vert (\dt k_\epsi^1,\dt k_\epsi^2)
	\vert_{X^{s+P-5/2}_{\underline{T}}}
	\lesssim 
	C(\underline{T}, \vert(\zeta,V)\vert_{X^{s+Q}_{\underline{T}}}),
\end{equation}
provided that $Q$ is large enough.\\
Using the fact that all the terms in the equation on $W$ are a gradient
of a scalar expression, as well as $W_\init$, it is possible to write
$w=(\tz,\nabla\psi)^T$, and $(\tz,\psi)$
solves
\begin{equation}\label{numerique}
	\left\lbrace
	\begin{array}{l}
	\dt \tz -{\mathcal T}_\mu \nabla\psi+\epsi
	f^1(\tz,\nabla\psi)=\epsi^2 k^1_\epsi\vspace{1mm},\\
	\dt \psi+\tz+\epsi f^2(\tz,\nabla\psi)=\epsi^2k_\epsi^2,
	\end{array}\right.
\end{equation}
with initial condition $(\tz,\psi)_\init=(\zeta^0,\psi^0)$.\\
Remarking now
that 
$\cG[0]\psi=\sqrt{\mu}\cT_\mu\nabla \psi$
and writing $U=(\tz,\psi)$, one can check
that (\ref{numerique}) can be written
\begin{equation}\label{num}
	\dt U+\cL U+\epsi \cA^{(1)}[U]
	=
	\epsi^2(k^1_\epsi,k^2_\epsi)^T,
\end{equation}
where $\cA^{(1)}[U]$ is given by the 
same formula (\ref{A}) as $\cA[U]$, but with the Dirichlet-Neumann operator
$\cG[\eps\tz]\psi$ replaced by the first order expansion
given in Proposition \ref{proptrunc}, and with the $O(\epsi^2)$ terms 
neglected. One thus gets 
$$
	\dt U+\cL U+\epsi \cA[U]
	=
	\epsi^2 H_\epsi,
$$
with $H_\epsi=(k^1_\epsi,k^2_\epsi)^T +\frac{1}{\epsi}(\cA^{(1)}[U]-\cA[U])$.
From (\ref{estdiff}), Proposition \ref{proptrunc}, and (\ref{bound}), one has
$(H_\epsi)_\epsi$  is uniformly bounded in $C([0,\frac{\underline T}{\epsi}],X^{s+P})\cap
C^1([0,\frac{\underline T}{\epsi}],X^{s+P-5/2})$ and we can therefore conclude
with Remark \ref{remstab} and Proposition \ref{proptaylor}. \qed
\end{proof}

\subsubsection{A remark on a model used for numerical computations}
\label{sectjustifnum}

The Dirichlet-Neumann operator is one of the main difficulties in the 
numerical computation of solutions to 
the water-waves equations (\ref{nondimww}) because it requires to solve
a $d+1$ ($d=1,2$ is the surface dimension) Laplace equation on a 
domain which changes at each time step. A common strategy is to replace
the full Dirichlet-Neumann operator by an approximation which requires
less computations. An efficient method, set forth in \cite{CGHHS}, consists
in replacing the Dirichlet-Neumann operator by its $n$-th order expansion
with respect to the surface elevation $\zeta$. When $n=1$, 
it turns out that the model
thus obtained is exactly the same as the system (\ref{num}) used in the
proof of Theorem \ref{theodeep}.\\
We can therefore use this theorem to state that: \emph{ the precision
of the modelization in the numerical computations of \cite{CGHHS} is of the
same order as the \underline{steepness} of the wave}.\\
One will easily check that when the $n$-th order expansion is used, then
the precision is of the same order as the $n$-th power of the steepness.

\appendix
\section{Nondimensionalization(s) of the equations}\label{appnd}

Depending on the value
of $\mu$, two distinct nondimensionalizations are commonly used in
oceanography (see for instance \cite{Dingemans2}). Namely, with
dimensionless quantities denoted with a prime:
\begin{itemize}
\item Shallow-water, ie $\mu\ll1$, one writes
\begin{equation}\label{nondimshallow}
	\begin{array}{llll}
	\dsp x=\lambda x',&
	y=\frac{\lambda}{\gamma}y',&
	z=dz',&
	t=\frac{\lambda}{\sqrt{gd}}t',\\
	\dsp\zeta=a\zeta',&
	\Phi=\frac{a}{d}\lambda\sqrt{gd}\Phi',&
	b=Bb'.&
	\end{array}
\end{equation}
\item Deep-water, ie $\mu\gg1$, one writes
\begin{equation}\label{nondimdeep}
	\begin{array}{llll}
	\dsp x=\lambda x',&
	y=\frac{\lambda}{\gamma}y',&
	z=\lambda z',&
	t=\frac{\lambda}{\sqrt{g\lambda}}t',\\
	\dsp \zeta=a\zeta',&
	\Phi=a\sqrt{g\lambda}\Phi',&
	b=Bb'.
	\end{array}
\end{equation}
\end{itemize}
Remark that when $\mu\sim 1$, that is when $\lambda\sim d$, both 
nondimensionalizations are equivalent, we introduce the following 
general nondimensionalization, which is valid for all $\mu>0$:
$$
	\begin{array}{llll}
	x=\lambda x',&
	y=\frac{\lambda}{\gamma}y',&
	z=d\nu z',&
	t=\frac{\lambda}{\sqrt{gd\nu}}t',\\
	\zeta=a\zeta',&
	\Phi=\frac{a}{d}\lambda\sqrt{\frac{gd}{\nu}}\Phi',&
	b=Bb',
	\end{array}
$$
where $\nu$ is a smooth function of $\mu$ such that $\nu\sim 1$ when $\mu\ll1$
and $\nu\sim \mu^{-1/2}(=\lambda/d)$ when $\mu\gg1$ 
(say, $\nu=(1+\sqrt{\mu})^{-1}$).

The equations of motion (\ref{eqbern}) then become (after dropping the
primes for the sake of clarity):
\begin{equation}
	\label{app1}
	\left\lbrace
	\begin{array}{l}
	\displaystyle \nu^2\mu\partial_{x }^2\Phi 
	+\nu^2\gamma^2\mu\partial_{y }^2\Phi 
	+\partial_{z }^2\Phi =0,
	\qquad \frac{1}{\nu}(-1+\beta b )\leq z \leq \frac{\eps}{\nu}\zeta ,\vspace{1mm}\\
	\displaystyle -\nu^2\mu\nag(\frac{\beta}{\nu}b )\cdot\nag\Phi +\partial_{z }\Phi =0,
        \qquad z =\frac{1}{\nu}(-1+\beta b ),\vspace{1mm}\\
	\displaystyle \partial_{t }\zeta -\frac{1}{\mu\nu^2}\big(-\nu^2\mu\nag(\frac{\eps}{\nu}\zeta )\cdot\nag\Phi
	+\partial_{z }\Phi \big)
	=0,
	\qquad  z = \frac{\eps}{\nu}\zeta ,\vspace{1mm}\\
	\displaystyle \partial_{t}\Phi+\frac{1}{2}\big(\frac{\eps}{\nu}\vert\nag\Phi\vert^2
        +\frac{\eps}{\mu\nu^3}(\partial_{z }\Phi )^2\big)+\zeta =0,
	\qquad z =\frac{\eps}{\nu}\zeta .
	\end{array}
	\right.
\end{equation}

In order to reduce this set of equations into a system of two evolution
equations, define the Dirichlet-Neumann operator 
$\cG^\nu_{\mu,\gamma}[\frac{\eps}{\nu}\zeta, \beta b]\cdot$ as
$$
  \cG^\nu_{\mu,\gamma}[\frac{\eps}{\nu}\zeta, \beta b]\psi
  =\sqrt{1+\vert\nabla(\frac{\eps}{\nu}\zeta)\vert^2}\partial_n\Phi_{\vert_{z=\frac{\eps}{\nu}\zeta}},
$$
with $\Phi$ solving the boundary value problem
$$
	\left\lbrace
	\begin{array}{l}
	\dsp \nu^2\mu\partial_{x }^2\Phi 
	+\nu^2\gamma^2\mu\partial_{y }^2\Phi 
	+\partial_{z }^2\Phi =0,
	\qquad \frac{1}{\nu}(-1+\beta b )\leq z \leq \frac{\eps}{\nu}\zeta,
	\vspace{1mm}\\
	\dsp \Phi_{\vert_{z=\frac{\eps}{\nu}\zeta}}=\psi,\qquad
	\partial_n\Phi_{\vert_{z=\frac{1}{\nu}(-1+\beta b)}}=0,	
	\end{array}\right.
$$
(as always in this paper, $\partial_n\Phi$ stands for the upwards conormal
derivative associated to the elliptic equation). As remarked in
\cite{Zakharov,CSS1,CSS2}, the equations (\ref{app1}) are
equivalent to a set of two equations on the free surface parameterization
$\zeta$ and the trace of the velocity potential at the surface $\psi=\Phi_{\vert_{z=\eps/\nu\zeta}}$ involving the Dirichlet-Neumann operator. Namely,
\begin{equation}
	\label{app2}
	\left\lbrace
	\begin{array}{l}
	\dsp \dt \zeta-\frac{1}{\mu\nu^2}\cG^\nu_{\mu,\gamma}[\frac{\eps}{\nu}\zeta, \beta b]\psi=0,\\	
	\dsp \dt\psi+\zeta+\frac{\eps}{2\nu}\vert\nag\psi\vert^2-
        \frac{\eps\mu}{\nu^3}
        \frac{(\frac{1}{\mu}\cG^\nu_{\mu,\gamma}[\frac{\eps}{\nu}\zeta, \beta b]\psi
        +\nu\nag(\eps\zeta)\cdot\nag\psi)^2}
        {2(1+\eps^2\mu\vert\nag\zeta\vert^2)}=0.
	\end{array}\right.
\end{equation}
In order to derive the system (\ref{nondimww}), let $\cG_{\mu,\gamma}[\eps\zeta, \beta b]\cdot$ be
the Dirichlet-Neumann operator $\cG_{\mu,\gamma}^\nu[\frac{\eps}{\nu}\zeta, \beta b]\cdot$ corresponding to the case $\nu=1$. One will easily check that
$$
	\forall\nu>0,\qquad
	\cG_{\mu,\gamma}[\eps\zeta, \beta b]
	=\frac{1}{\nu}\cG_{\mu,\gamma}^\nu[\frac{\eps}{\nu}\zeta, \beta b],
$$
so that plugging this relation into (\ref{app2}) yields
$$
	\left\lbrace
	\begin{array}{l}
	\dsp \dt \zeta-\frac{1}{\mu\nu}\cG_{\mu,\gamma}[\eps\zeta, \beta b]\psi=0,\\
	\dsp \dt \psi+\zeta+\frac{\eps}{2\nu}\vert\nag\psi\vert^2-\frac{\eps\mu}{\nu}
        \frac{(\frac{1}{\mu}\cG_{\mu,\gamma}[\eps\zeta, \beta b]\psi+\nag(\eps\zeta)\cdot\nag\psi)^2}
        {2(1+\eps^2\mu\vert\nag\zeta\vert^2)}=0.
	\end{array}\right.
$$

%\noindent
%{\bf Acknowledgment.} 
\begin{acknowledgement}
This work was supported by the ACI Jeunes Chercheuses et
Jeunes Chercheurs ``Dispersion et nonlin\'earit\'e''.
\end{acknowledgement}

\end{document}